\documentclass[preprint,12pt]{elsarticle}
\usepackage{amsfonts,amsmath,amssymb}
\usepackage{tikz}

\newtheorem{thm}{Theorem}[section]
\newtheorem{lem}[thm]{Lemma}
\newtheorem{co}[thm]{Conjecture}
\newdefinition{df}{Definition}[section]
\newdefinition{rem}{Remark}[section]
\newdefinition{ex}{Example}[section]
\newproof{pf}{Proof}
\newproof{pot}{Proof of Theorem}
\numberwithin{equation}{section}

\journal{
}

\begin{document}

\begin{frontmatter}

\title{The Ismagilov conjecture over a finite field ${\mathbb F}_p$
}
\author{A.V.~Kosyak
}
\ead{kosyak02@gmail.com}
\address{Max-Planck-Institut f\"ur Mathematik, Vivatsgasse 7, D-53111 Bonn, Germany}
\address{Institute of Mathematics, Ukrainian National Academy of Sciences,
3 Tereshchenkivs'ka Str., Kyiv, 01601, Ukraine}

\begin{abstract}
 We construct the so-called quasiregular representations of the
group $B_0^{\mathbb N}({\mathbb F}_p)$ of infinite upper triangular matrices with coefficients in a
finite field and give the criteria of theirs irreducibility
in terms of the initial measure. These representations are particular case of the Koopman representation
hence,  we find new conditions of its irreducibility.
Since the field ${\mathbb F}_p$ is compact some new operators in the commutant
emerges. Therefore, the Ismagilov conjecture in the case of the finite field
should be corrected.
\end{abstract}


\end{frontmatter}
\newpage
\tableofcontents
\section{ Introduction }
Let $\hat G$ be a {\it dual space of a group} $G$, i.e., the set
of all equivalence classes of unitary irreducible representations
of the group $G$. Our {\it far reching goal} is to describe $\hat{G}$ for
$G=\varinjlim_{n}G_n$ where $G_n=B(n,{\bf k})$ is the group of all
upper triangular matrices with units on the principal  diagonal
with natural inclusion $G_n\subset G_{n+1}$, where ${\bf
k}={\mathbb R}$ or ${\bf k}={\mathbb F}_p$ is a finite field
${\mathbb Z}/p{\mathbb Z}$, $p$ is prime.

We mention here only some results concerning representations of algebraic groups over a finite field. The book by G,~Lusztig~\cite{Lusz84} presents a classification of all (complex) irreducible representations of a
reductive group with connected centre, over a finite field. To achieve this, the author uses etale intersection
cohomology, and detailed information on representations of Weyl groups.

From the article by P.~Deligne and G.~Lusztig~\cite{DelLusz76}:
``Let us consider a connected reductive algebraic group $G$,
defined over a finite field ${\mathbb F}_q$, with {\it Frobenius map}
$F$. We shall be concerned with the representation theory of the
finite group $G^F$, over field of characteristic $0$.
In 1968, Macdonalds conjectured, on the basis of the character table known at that time for
$({\rm GL}_4, {\rm Sp}_4 )$ that should be a well defined correspondence which, to any $F$-stable maximal
torus $T$ of $G$ and a character $\theta$ of $T^F$ in general position, associate an irreducible
representation of $G^F$; moreover, if $T$ modulo the centre of $G$ is {\it anisotropic} over ${\mathbb F}_q$, the
corresponding representation of $G^F$ should be {\it caspidal}. In this paper we prove Macdonald's conjecture.
 More precisely, for $T$ as above and $\theta$ an arbitrary character of $T^F$ we construct virtual
 representations $R^\theta_T$ which have all the required properties.''
 \index{torus $T$!anisotropic}
 \index{representation!caspidal}

The group $G_n=B(n,{\mathbb F}_p)$ is finite, hence the set
$\hat{G}_n$, in principal, is known (it is numerated by the set of all conjugacy
classes) and all irreducible representations are contained in the
regular representation. ``Nevertheless the complete classification
of the complex irreducible representations of this group has long
been known to be a difficult task'' \cite{Yan06} . Recently in 2006,  Ning Yan,
have introduced  in
\cite{Yan06} ``\,a natural variant of the orbit method, in which the
central role is played by certain {\it clusters of coadjoint
orbits}. This method of clusters leads to the construction of a
subring in the representation ring of $B(n,{\mathbb F}_p)$ that is
``\,rich in structure but pleasantly comprehensible''.

The article by V.~Gorin, A.~Vershik and S.~Kerov~\cite{GorKerVer14} is devoted to the representation
theory of locally compact infinite-dimensional group $GLB$ of {\it almost upper-triangular infinite
matrices}
\index{matrices!upper-triangular!almost}
over the finite field with $q$ elements.
From~\cite{GorKerVer14}:`` The group $GLB$ consist of all almost
triangular matrices of infinite order. An infinite matrix
$g=(g_{ij}),\,i,j=1,2,\dots,$ is said to be almost triangular if
the number of its subdiagonal elements $g_{ij}\not=0,\,i>j,$ is
finite. This group was defined by S.~Kerov, A.~Vershik, and A.~Zelevinsky
in 1982 as an adequate for  $n=\infty$ analogue of general linear
groups $GL(n,q)$. It serves as an alternative to $GL(\infty,q)$,
whose representation theory is poor. Our most important results
are the description of semi-finite unipotent traces (characters)
of the group $GLB$ via certain probability measures on the Borel
subgroup $B$ and the construction of the corresponding von Neumann
factor representations of type $II_\infty$.''

Coming back to our group $B_0^{\mathbb N}({\mathbb F}_p)=\varinjlim_{n}G_n$ where
$G_n=B(n,{\mathbb F}_p)$, we mention that with the natural
homomorphism $ p^{n}_{n+1}:G_{n+1}\to G_{n}$ (see (\ref{df.p})) we can associate an inclusion
$\widehat{G}_n\to \widehat{G}_{n+1}$  therefore,
$\widehat{G}\supset\bigcup_{n}\widehat{G_n}$. In the case ${\bf
k}={\mathbb R}$ one may use Kirillov's orbit method to describe
$\widehat{G_n}$.  We define $ p^{n}_{n+1}$ as follows:
\begin{equation}\label{df.p}
B(n+1,{\bf
k})\ni x=x^{n+1}x_n\mapsto p^{n}_{n+1}(x)=x_n\in
B(n,{\bf k}),
\end{equation}
where for
$$
x=I+\sum_{1\leq k<m\leq n+1}x_{km}E_{km}\in B(n+1,{\bf
k}),\quad\text{we set}
$$
\begin{equation}
x^{n+1}=I+\sum_{k=1}^nx_{kn+1}E_{kn+1},\quad x_n=I+\sum_{1\leq
k<m\leq n}x_{km}E_{km}.
\end{equation}
Obviously, $x=x^{n+1}x_n$ and $B(n+1,{\bf
k})$ is a semi-direct
product
\begin{equation}
\label{G(n+1)=A-ltimes-G(n)} B(n+1,{\bf
k})={\bf
k}^{n}\ltimes B(n,{\bf
k}).
\end{equation}
\begin{rem}
\label{B(F(p))-compact}
The group $B^{\mathbb N}({\mathbb F}_p)$ is compact, the corresponding {\it Haar measure} on this group is infinite tensor product of the normalised invariant measures on ${\mathbb F}_p$, where $\mu_{inv}^{kn}$ is defined by   (\ref{mu^kn(inv)})
\begin{equation}
\label{Haar-mes}
h=\mu_{inv}=\otimes _{1\leq k<n}\mu_{inv}^{kn}.
\end{equation}
\index{Haar measure on $B^{\mathbb N}({\mathbb F}_p)$}
Therefore, all irreducible representations of the group $B^{\mathbb N}({\mathbb F}_p)$ are {\it finite-dimensional} and are contained in the decomposition of the the regular representation of $B^{\mathbb N}({\mathbb F}_p)$. Moreover we have
\begin{equation}
\label{dual(projB,(nF(p)))}
\widehat{B^{\mathbb N}({\mathbb F}_p)}=\bigcup_{n\geq 1}\widehat{B(n,{\mathbb F}_p)}.
\end{equation}
The group  $B_0^{\mathbb N}({\mathbb F}_p)=\varinjlim_{n}B(n,{\mathbb F}_p)$ is subgroup of $B^{\mathbb N}({\mathbb F}_p)$
therefore, $\widehat{B_0^{\mathbb N}({\mathbb F}_p)}\supset \widehat{B^{\mathbb N}({\mathbb F}_p)}$. We construct the {\it infinite-dimensional} irreducible representations of the group $B_0^{\mathbb N}({\mathbb F}_p)$ as {\it quasiregular representations}. The most important observation is that
the measure on the homogeneous space $X^m$ (see (\ref{X^m}) below), we use for this, has the property that its projection and the projection of the Haar measure on subspace $X^{(k)}$ (see (\ref{X^(k)}) below) are {\it orthogonal}.
\end{rem}
\begin{rem}
\label{Ism-corr}
 We show that  $\widehat{G}\setminus\left(\bigcup_{n}\widehat{G_n}\right)\not=\emptyset.$
Namely, $\widehat{G}\setminus\left(\bigcup_{n}\widehat{G_n}\right)$
contains ``regular''  and ``quasiregular'' irreducible representations of the group $G$.
The criteria of the irreducibility of quasiregular
representations in the case of the field ${\mathbb F}_p$ is
established.  Some new conditions of the irreducibility are found, if we compare
with the Ismagilov conjecture~\ref{8.co.Ism} in the case $k={\mathbb R}$.
\end{rem}

Recall the definition of the {\it Koopman representation}. {\it Quasiregular representation} is a particular case of it.
Let we have the {\it measurable
action} $\alpha:G\rightarrow {\rm Aut}(X)$ of  the group $G$  on a
measurable space $(X,\mu)$ with $G$-{\it
quasi-invariant measure} $\mu$ having the following property:
$\mu^{\alpha_t}\sim\mu\quad\forall t\in G$. With these dates one
can associate a {\it representation} $ \pi^{\alpha
,\mu,X}:G\rightarrow U(L^2(X,d\mu))$ defined as follows:
\begin{equation}\label{8.Rep-pi}
(\pi^{\alpha
,\mu,X}_tf)(x)=(d\mu(\alpha_{t^{-1}}(x))/d\mu(x))^{1/2}f(\alpha_{t^{-1}}(x)),\quad
f\in L^2(X,\mu).
\end{equation}
In the case of invariant measure this representations is called {\it Koopman's representation}, see  \cite{Koo31}.
\index{represenation!Koopman's}
\index{subgroup!centarlizer}
Consider the {\it centralizer}  $Z_{{\rm Aut}(X)}(\alpha(G))$ of
the subgroup $\alpha(G)=\{\alpha_t\in{\rm Aut}(X)\mid t\in G\}$ in
the group ${\rm Aut}(X)$
$$ Z_{{\rm Aut}(X)}(\alpha(G))=\{g\in {\rm Aut}(X)\mid
\{g,\alpha_t\}=g\alpha_tg^{-1}\alpha_t^{-1}=e\,\,\forall t\in
G\}.$$ The following conjecture was proved for some
infinite-dimensional groups.
\begin{co} [\rm Kosyak, \cite{Kos02.3,Kos03}]
\label{8.co.G-Ism} The representation $\pi^{\alpha
,\mu,X}:G\rightarrow $\\ $U(L^2(X,\mu))$ is irreducible if and only if
\par
1) $\mu^g\perp \mu\,\,\forall g\in Z_{{\rm
Aut}(X)}(\alpha(G))\backslash{\{e\}},\,\,$ (where $\perp$ stands
for singular),
\par
2) the measure $\mu$ is $G$-ergodic.
\end{co}
We recall that a measure $\mu$ is $G$-ergodic if
$f(\alpha_t(x))=f(x)\,\,\mu-$a.e. for all $t\in G$ implies
$f(x)=const\,\,\mu-$a.e. (almost everywhere) for all functions
$f\in L^1(X,\mu)$.
\begin{rem}
\label{8.r.modify}
We show that Conjectures~\ref{8.co.G-Ism} dos not hold over the finite field  ${\mathbb F}_p$ and indicate how it  should be modified (see Conjecture~\ref{co.Qreg(F_p)-irr}). Nevertheless, we hope that the Ismagilov conjecture (Conjecture~\ref{8.co.Ism}) still holds in the case
of the finite field ${\mathbb F}_p$ (see Conjecture~\ref{co.Reg(F_p)-irr} and Remark~\ref{r.cond(3)-}).
\end{rem}

\section{Regular and guasiregular representations of the group $B_0^{\mathbb N}({\mathbb F}_p)$}
\subsection{Regular and quasiregular representations of infinite-dimensional groups, the case of ${\bf k}={\mathbb R}$}
Let the  group $G$ {\it be a locally compact}, $X=G$ and $h$ be
the Haar measure. If $\alpha$ is right or left action of the group
$G$ on itself  then $\rho=\pi^{R ,h,G}$ and $\lambda=\pi^{L,h,G}$
are well known {\it right} and {\it left regular representations}.
{\it Quasiregular representation} is a particular case of the
representation $\pi^{\alpha ,\mu,X}$ defined by (\ref{8.Rep-pi}), where $X=H\backslash
G,\,\,H$ is some closed subgroup of $G$ and $\mu$ is some
$G$-quasi-invariant measure on $X$.

Recall the notions of the regular and quasiregular representations for infinite-dimensional groups and the Ismagilov conjecture~\ref{8.co.Ism}.
To define a "{\it regular representation}" for {\it
infinite-dimensional group} $G$ the initial group $G$ as a
candidate  for $X$ is not suitable since on $G$ there is no Haar
(invariant) measure (Weil, \cite{Weil53}) no a $G$-quasi-invariant
measure (Xia Dao-Xing, \cite{{Xia78}}).
It is natural to consider some larger topological group $\widetilde G$
containing the initial group $G$ as the dense subgroup and a
$G$-quasi-invariant measure $\mu$ on $\widetilde G$.

\begin{df}
\label{8.def.reg.rep}
Representations $T^{\alpha,\mu}=\pi^{\alpha
,\mu,\widetilde{G}}$ where $\alpha=R$ (resp.  $\alpha=L$) we call
the {\it right (resp. the left) regular representation} of the
group $G$.
\end{df}

\begin{co}[\rm Ismagilov, 1985]
\label{8.co.Ism} The right regular representation
$T^{R,\mu}:G\rightarrow U(L^2(\tilde G,\mu))$ is irreducible if
and only if
\par
1) $\mu^{L_t}\perp \mu\quad\text{for all}\quad t\in G\backslash{\{e\}},\,\,$
\par 2) the measure $\mu$ is $G$-ergodic.
\end{co}
Similarly, we can generalize  the notion of  {\it quasiregular
representation} of a group $G$ associated with some subgroup $H$
using  a suitable completion $\tilde X\!=\!\!\widetilde{H\backslash G}$
of the homogeneous space $X=H\backslash G$ and constructing  some
$G$-right quasi-invariant measure $\mu$ on $\tilde X$.

Consider the group $G=B_0^{\mathbb N}(\mathbb
R)=\varinjlim_{n}B(n,{\mathbb R})$.
 Let us fix the space $X$ and the measure
$\mu$ on $X$ as follows, where $E_{kn}$ are matrix units of
infinite order:
$$
X=\widetilde{G}=B^{\mathbb N}=\{I+x\mid x=\sum_{1\leq
k<n}x_{kn}E_{kn},\,x_{kn}\in{\mathbb R}\},
$$
$$
d\mu_b(x)=\otimes_{1\leq k<n}(b_{kn}/\pi)^{1/2}
\exp(-b_{kn}x_{kn}^2)dx_{kn},\,\,b=(b_{kn})_{k<n}.
$$
\begin{thm}[\cite{Kos90,Kos92}]
Ismagilov's conjecture  holds, i.e., $T^{R,\mu}\in \hat{G}$ if and
only if $\mu^{L_t}\perp \mu\quad\text{for all}\quad t\in
G\backslash{\{e\}},\,\,$ and the measure $\mu$ is $G$-ergodic.
Moreover
$$
T^{R,\mu_1}\sim T^{R,\mu_1}\quad\text{if and only if}\quad \mu_1\sim \mu_2.
$$
\end{thm}
{\it Quasiregular representations}  for the group $G=B_0^{\mathbb
N}$. Let us consider two subgroups of the group $B^{\mathbb
N}({\mathbb R})$:
\begin{align*}
B_m=&\,\{I+x\in B^{\mathbb N}\mid
x=\sum_{m<k<n}x_{kn}E_{kn},\,\,x_{kn}\in{\mathbb R}\},\\
B^m=&\,\{I+x\in B^{\mathbb N}\mid  x=\sum_{1\leq k\leq
m,k<n}x_{kn}E_{kn}\}.
\end{align*}
The group $B^{\mathbb N}(\mathbb R)$ is a semi-direct product
$B^{\mathbb N}=B_m\rtimes B^m$. Fix the corresponding
decomposition $x=x_m\cdot x^m$.  Define $X^m=B_m\setminus
B^{\mathbb N}\simeq B^m.$ Right action $R$ of $G$ on $X^m$ is well
defined $R_t(x):=(xt^{-1})^m,\,\,x \in X^m,\,\,t\in G$. The
measure on $X^m$ is defined by
$$
d\mu_{(b,a)}^m(x)=\otimes_{1\leq k\leq m,\,k<n}(b_{kn}/\pi)^{1/2}
\exp(-b_{kn}(x_{kn}-a_{kn})^2)dx_{kn}.
$$
Quasiregular representation is defined by
$T^{R,\mu_{(b,a)},m}=\pi^{R,\mu_{(b,a)},X^m}$.
\begin{thm}[\cite {Kos02.3}]
Conjecture~\ref{8.co.G-Ism} holds, i.e., quasiregular
representation $T^{R,\mu_{(b,a)},m}$ is irreducible if and only
if conditions 1) and 2) of Conjecture~\ref{8.co.G-Ism} holds.
Moreover,
$$
T^{R,\mu_1,m}\sim T^{R,\mu_2,n}\Leftrightarrow m=n\text{\, and
\,}\mu_1\sim \mu_2.
$$
\end{thm}
Conjecture~\ref{8.co.G-Ism} for quasiregular representations of the group
$B_0^{\mathbb N}({\mathbb R})$ is proved  by A.~Kosyak and
S.~Albeverio in \cite{KosAlb05C} for a tensor product of arbitrary
one-dimensional measures and for more general Gaussian measures by
A.~Kosyak and S.~Albeverio in \cite{KosAlb06J}.

\subsection{Regular and quasiregular representations and criteria
of irreducibility, the case ${\bf k}={\mathbb F}_p$}
{\it We show that in the case ${\bf k}={\mathbb F}_p$  Conjecture
\ref{8.co.G-Ism} does not hold but may be corrected easily}. More
precisely, two conditions of the irreducibility 1) and 2) of the
Conjecture \ref{8.co.G-Ism} are not sufficient, since the commutant
of the right quasiregular representation may be generated not only
by the operators of the left representations, as in the case of
${\bf k}={\mathbb R}$, but also by some operators acting in $L^2$
on infinite rows and existing only in the case when the
corresponding measures are equivalent with the infinite tensor
product of the invariant measures (see conditions 3)  of  Conjecture
\ref{co.Qreg(F_p)-irr}).

Define a {\it quasiregular representations}  for the group
$B_0^{\mathbb N }({\mathbb F}_p)=\varinjlim_{n}B(n,{\mathbb
F}_p)$. Let us consider two subgroups of the group of all upper
triangular matrices $B^{\mathbb N }({\mathbb F}_p)$:
\begin{align*}
B_m({\mathbb F}_p)=&\{I+x\in B^{\mathbb N}({\mathbb F}_p)\mid
x=\sum_{m<k<n}x_{kn}E_{kn},\,\,x_{kn}\in{\mathbb F}_p\},\\
B^m({\mathbb F}_p)=&\{I+x\in B^{\mathbb N}({\mathbb F}_p)\mid
x=\sum_{1\leq k\leq m,k<n}x_{kn}E_{kn},\,\,x_{kn}\in{\mathbb
F}_p\}.
\end{align*}
The group $B^{\mathbb N}({\mathbb F}_p)$ is semi-direct product
$B^{\mathbb N}({\mathbb F}_p)=B_m({\mathbb F}_p)\ltimes
B^m({\mathbb F}_p)$. Fix the corresponding decomposition
$x=x_m\cdot x^m$. Define the homogeneous space
\begin{equation}
\label{X^m}
X^m:=B_m({\mathbb F}_p)\setminus
B^{\mathbb N}({\mathbb F}_p)\simeq B^m({\mathbb F}_p).
\end{equation}

The measure $\mu_\alpha=\mu_\alpha^{(m)}$ on the space $X^m$ is defined as infinite
tensor product:
\begin{equation}
\label{mu^m-alpha} \mu_\alpha^{(m)}:=\otimes_{1\leq k\leq
m,k<n}\mu_{\alpha_{kn}}=\otimes_{k=1}^m\mu^k_\alpha,\quad {\rm
where}\quad \mu^k_\alpha:=\otimes_{n=k+1}^\infty\mu_{\alpha_{kn}},
\end{equation}
of the probability measures $\mu_{\alpha_{kn}}$  on ${\mathbb
F}_p$ defined as follows:
$${\mathbb F}_p\ni r\mapsto
\mu_{\alpha_{kn}}(r)=\alpha_{kn}(r)>0\quad\text{and}\quad
\sum_{r\in{\mathbb F}_p}\alpha_{kn}(r)=1.
$$
{\it The  right action $R:B_0^{\mathbb N }({\mathbb
F}_p)\rightarrow {\rm Aut}(X^m)$ of the group $B_0^{\mathbb N
}({\mathbb F}_p)$ on the factor space $X^m\!=\!B_m({\mathbb
F}_p)\setminus B^{\mathbb N}({\mathbb F}_p)$ is well defined by} $
R_t(x)\!=\!(xt^{-1})^m$.
\begin{lem}
\label{l.Qreg-admis(p)}
The right action  of the group $B_0^{\mathbb N }({\mathbb F}_p)$  on the space $X^m$ is {\rm admissible}, i.e., $(\mu_\alpha)^{R_t}\sim \mu_\alpha\,\,\forall
t\in B_0^{\mathbb N }({\mathbb F}_p).$
\end{lem}
{\it Quasiregular representation} $T^{R,\mu_\alpha,m}$ is defined in the space
$L^2(X^m,\mu_\alpha)$ by (\ref{Rep-q-fin}).
For $1\leq k\leq m$ define the measures $\mu_\alpha^k$ and
$\mu_{inv}^k$ as follows:
\begin{equation}
\label{mu^k(alpha),mu^k(inv)}
\mu^k_\alpha:=\otimes_{n=k+1}^\infty\mu_{\alpha_{kn}},\quad
\mu_{inv}^k=\otimes_{n=k+1}^\infty\mu_{inv}^{kn},\quad
\mu_{inv}=\otimes_{k=1}^m\mu_{inv}^k,
\end{equation}
where $\mu_{inv}^{kn}$ is the normalized invariant measure on ${\mathbb
F}_p$, i.e.,
\begin{equation}
\label{mu^kn(inv)}
\mu_{inv}^{kn}(r)\!=\!p^{-1},\quad r\in{\mathbb F}_p.
\end{equation}
\begin{co}
\label{co.Qreg(F_p)-irr}
 Let $m\in {\mathbb N}$. The quasiregular
representation $T^{R,\mu_\alpha,m}:B_0^{\mathbb N }({\mathbb
F}_p)\rightarrow U(L^2(X^m,\mu_\alpha))$ of the group
$B_0^{\mathbb N }({\mathbb F}_p)$ is irreducible if and only if
\par
1) $\mu^{L_t}_\alpha\perp \mu_\alpha\,\,\forall t\in B(m,{\mathbb
F}_p)\backslash{\{e\}},$
\par
2) the measure $\mu_\alpha$ on the space $X^m$ is $B_0^{\mathbb N }({\mathbb
F}_p)$-right-ergodic,
\par
3)  for the measure $\mu_\alpha=\otimes_{k=1}^m\mu^k_\alpha$ holds the following:
$\mu^k_\alpha\perp\mu_{inv}^k$ for all $1\leq k\leq m$ (By Remark~\ref{(3):only-k=m}
it is
sufficient to verify this condition only for  $k=m$),
\par
4) two irreducible representations $T^{R,\mu_\alpha,m}$ and $T^{R,\mu_\beta,n}$ are equivalent\\
$T^{R,\mu_\alpha,m}\!\sim \!T^{R,\mu_\beta,n}$ if and only if
$m\!=\!n$ and $\mu_\alpha\!\sim\!\mu_\beta.$
\end{co}
We would like to mention here an interesting {\it problem to solve}.
Define a {\it regular representation of the group } $B_0^{\mathbb N }({\mathbb
F}_p)$ as before. On the group $B^{\mathbb N }({\mathbb F}_p)$   of all upper
triangular matrices define the measure $\mu_\alpha$ as follows:
\begin{equation}
\label{mu-alpha} \mu_\alpha:=\otimes_{k<n}\mu_{\alpha_{kn}}=\otimes_{k=1}^\infty\mu^k_\alpha.
\end{equation}
\begin{lem}
\label{l.Reg-admis(p)}
The right action  of the group $B_0^{\mathbb N }({\mathbb F}_p)$  on the group $B^{\mathbb N }({\mathbb F}_p)$ is
admissible, i.e., $(\mu_\alpha)^{R_t}\sim \mu_\alpha\,\,\forall
t\in B_0^{\mathbb N }({\mathbb F}_p).$
\end{lem}
\begin{co}
\label{co.Reg(F_p)-irr} The regular representation $T^{R,\mu_\alpha}:B_0^{\mathbb N }({\mathbb
F}_p)\rightarrow $\\$U(L^2(B^{\mathbb N }({\mathbb F}_p),\mu_\alpha))$ of the group
$B_0^{\mathbb N }({\mathbb F}_p)$ is irreducible if and only if
\par
1) $\mu^{L_t}_\alpha\perp \mu_\alpha\,\,\forall t\in B_0^{\mathbb N }({\mathbb
F}_p)\backslash{\{e\}},$
\par
2) the measure $\mu_\alpha$ on the group $B^{\mathbb N }({\mathbb F}_p)$ is $B_0^{\mathbb N }({\mathbb
F}_p)$-right-ergodic,
\par
3) two irreducible representations $T^{R,\mu_\alpha}$ and $T^{R,\mu_\beta}$ are equivalent
$T^{R,\mu_\alpha}\!\sim \!T^{R,\mu_\beta}$ if and only if  $\mu_\alpha\!\sim\!\mu_\beta.$
\end{co}
\begin{rem}
\label{r.cond(3)-}
By Remark~\ref{(3):only-k=m}, the condition 3) of Conjecture~\ref{co.Qreg(F_p)-irr} will disappear in Conjecture~\ref{co.Reg(F_p)-irr}!
\end{rem}
\subsection{Idea of the proof of the irreducibility of the regular and quasiregular representations}
\label{8.Idea-irr}
Below we show that conditions 1)--3) of
Conjecture~\ref{co.Qreg(F_p)-irr} and  1)--2) of
Conjecture~\ref{co.Reg(F_p)-irr}) are necessary for the
irreducibility of the representation $T^{R,\mu_\alpha,m}$ (resp.
of $T^{R,\mu_\alpha}$). {\it The remaining part of the chapter is
devoted to the proof of the fact that these conditions are
sufficient for the irreducibility of the representations
$T^{R,\mu_\alpha,m}$}.

{\bf The conditions 1) -- 3) of the conjecture are necessary for
the irreducibility of $T^{R,\mu_\alpha,m}$ and $T^{R,\mu_\alpha}$}.
Indeed,  let conditions 1) does not hold, then
$\mu^{L_s}_\alpha\sim\mu_\alpha$ for some $s\in B(m,{\mathbb
F}_p)\backslash{\{e\}}$ therefore, the operator
$T^{L,\mu_\alpha,m}_s\!=\!\pi^{L,\mu_\alpha,X^m}_s$ is well
defined and commutes with the representation $T^{R,\mu_\alpha,m}$.
Similarly, the operator $T^{L,\mu_\alpha}_s$ commutes with $T^{R,\mu_\alpha}$.

The necessity of the condition 2) is evident. Indeed, if the
measure $\mu_\alpha$ is not $G$-ergodic on the space $X^m$
(resp. $X\simeq B^{\mathbb N }({\mathbb F}_p)$) then $X^m=X_1\cup
X_2$ (resp. $X=X_1\cup X_2$)  where $X_k$ are $G-$invariant and
$\mu_\alpha(X_k)>0,\,\,k=1,2.$ In this case
$L^2(X^m,\mu_\alpha)=H_1\oplus H_2$  (resp.
$L^2(X,\mu_\alpha)=H_1\oplus H_2$) is a direct sum of two
nontrivial $G$-invariant subspaces.

To explain  the condition 3) we define the {\it elementary
representations} $T^{R,\mu_\alpha^{k},(k)}$ of the group $G$ as
follows. Consider the subspace
\begin{equation}
\label{X^(k)}
X^{(k)}\!=\!\{I\!+\!\sum_{n=k+1}^\infty x_{kn}E_{kn}\}
\end{equation}
 of the
space $X^m$ and the projection $\mu^k_\alpha$ of the measure
$\mu_\alpha$ on the subspace $X^{(k)}$,  then
$$
X^m=X^{(m)}X^{(m-1)}\dots\,X^{(1)},\quad \mu_\alpha=\otimes_{k=1}^m\mu^k_\alpha,\quad
\text{where}\quad \mu^k_\alpha:=\otimes_{n=k+1}^\infty\mu_{\alpha_{kn}}.
$$
In this case the following decomposition of the representation
$T^{R,\mu_\alpha,m}$ holds:
$$
T_t^{R,\mu_\alpha,m}
=\otimes_{k=1}^mT_t^{R,\mu_\alpha^{k},(k)}\quad\text{in}\quad
L^2(X^m,\mu_\alpha)=\otimes_{k=1}^mL^2(X^{(k)},\mu_\alpha^{k}).
$$
We shall use the following {\it notations}
\begin{equation}
\label{T_{kn},T_{kn}(r)}
T_{kn}:=T^{R,\mu_\alpha,m}_{I-E_{kn}},\quad
T_{kn}(r):=T^{R,\mu_\alpha^{r},(r)}_{I-E_{kn}}.
\end{equation}
The following decomposition holds for the quasiregular representation $T^{R,\mu_\alpha,m}:$
\begin{equation}
\label{T(kn)=otimes}
T_{kn}=\otimes_{r=1}^kT_{kn}(r),\,\,1\leq k\leq m,\quad
T_{kn}=\otimes_{r=1}^m T_{kn}(r),\,k> m.
\end{equation}
For the regular representations $T^{R,\mu_\alpha}$  in
$L^2(B^{\mathbb N }({\mathbb F}_p),\mu_\alpha)=\otimes_{k=1}^\infty$ $L^2(X^{(k)},\mu_\alpha^{k})$
we have:
\begin{equation}
\label{Reg=otimes(r)}
\quad T^{R,\mu_\alpha,m}_t
=\otimes_{k=1}^\infty T^{R,\mu_\alpha^{k},(k)}_t,\quad
T_{kn}=\otimes_{r=1}^kT_{kn}(r),\,\,k<n.
\end{equation}

In Section~\ref{8.4.2Centre} we describe the commutant
$({\mathfrak A}^m)'$ of the von Neumann algebra
\begin{equation}
\label{A^m:=}
{\mathfrak
A}^m:=\big(T^{R,\mu_\alpha,m}_t\mid t\in G\big)''.
\end{equation}
 To be more precise, define the {\it Laplace operators}
$\Delta^{(m)}$ and $\Delta_k$ where
\begin{equation}
\label{Del,Frob} \Delta^{(m)}=\prod_{k=1}^m\Delta_k,\quad
\Delta_k:=\prod_{n=k+1}^\infty
p^{-1}C(T_{kn}(k))\quad\text{and}\quad C(T):=\sum_{r\in{\mathbb
F}_p}T^r.
\end{equation}
\index{operator!Laplace}
By Lemma~\ref{Lap-exist} we conclude that
{\it   the operator $\Delta_k$ is well defined and belongs to the commutant $({\mathfrak A}^m)'$
of the corresponding von Neumann algebra  ${\mathfrak
A}^m$ if $\mu_{\alpha}^k\sim\mu_{inv}^k$ for some $1\leq k\leq m$.}
This shows that condition 3) of Conjecture~\ref{co.Qreg(F_p)-irr}
are necessary conditions of the irreducibility of the
representation $T^{R,\mu_\alpha,m}$.
\begin{rem}
\label{r.only m=1,2}
 We were able to prove Conjecture~\ref{co.Qreg(F_p)-irr} only
in the case $m=1,\,\,p$ is arbitrary and $m=2,\,\,p=2$. The
general case of $m$ and $p$ is open. We shall try to study these
cases later.
\end{rem}
\begin{rem}
\label{r.Idea-irr} {\rm Idea to prove the irreducibility.} Roughly speaking, to
prove that conditions 1) -- 3) are sufficient for the
irreducibility, it is sufficient to show that in this case
operators $T_{\alpha_{kn}}$ defined by (\ref{T-alpha}) and
\begin{equation}
\label{X(kn)}
x_{kn}:={\rm diag}(0,1,\dots,p-1)=\sum_{r\in{\mathbb F}_p}rE_{rr},\quad 1\leq k\leq m,\quad k<n,
\end{equation}
acting on the Hilbert spaces $H_{kn}$
\begin{equation}
\label{H(kn)}
H_{kn}:=H_{\alpha_{kn}}:=L^2({\mathbb F}_p,\mu_{\alpha_{kn}}),
\end{equation}
belong to the von Neumann algebra ${\mathfrak A}^m$ generated by the representation $T^{R,\mu_\alpha,m}$. To be more precise, consider two infinite families of operators $X_k$ and $T_r$ defined as follows: $X_k=(x_{kn}\mid k<n)$ and $T_r=(T_{\alpha_{rn}}\mid r<n)$ for $1\leq k,r\leq m$.
For $m=1$ we prove that
$X_1\subset {\mathfrak A}^1$ therefore, $({\mathfrak A}^1)'\subset (L^\infty(X_1))'=L^\infty(X_1)$ since the von Neumann algebra $L^\infty(X_1)$ is maximal abelian (see Definition~\ref{max.ab.sub}).
For $m=2$ we prove that, depending on the measure, one of the families $(X_1,X_2),\,\,(X_1,T_2),\,\,(T_1,X_2),\,\,(T_1,T_2)$ belong to ${\mathfrak A}^2$.
For an arbitrary $m$ it is sufficient to prove that one of the following families
$(F_1,F_2,\dots,F_m)$ belongs to the von Neumann algebra ${\mathfrak A}^m$ where $F_k$ is $X_k$ or $T_k$ for $1\!\leq\! k\!\leq\! m$. To prove the irreducibility it is sufficient to prove that the von Neumann algebra $L^\infty(F_1,F_2,\dots,F_m)$ is maximal abelian therefore, $({\mathfrak A}^m)'\!\subset\! L^\infty(F_1,F_2,\dots,F_m)$ and use the ergodicity of the measure $\mu_\alpha$.
\index{algebra!von Neumann!maximal abelian}
\end{rem}
\begin{rem}
\label{otimes=1} For shortness we shall use the same notations
$A_k$ for the operator $A_k$ acting on the Hilbert space $H_k$ and
the operator ${\mathcal A}_k=I\otimes \dots\otimes I \otimes A_k
\otimes I\otimes \dots$ acting on the finite ${\mathcal
H}_r=\otimes_{n=1}^r H_n$ or infinite tensor product ${\mathcal
H}=\otimes_{n=1}^\infty H_n$.
\end{rem}
\section{The space $X$ and the measure}
Let us consider the finite field ${\mathbb F}_p\!=\!{\mathbb
Z}/p{\mathbb Z}$ of $p$ elements ${\mathbb F}_p\!=\!\{0,1,...,p-\!1\}$.
The group $B_0^{\mathbb N}({\mathbb F}_p)$ is defined as the
inductive limit (with natural inclusion) $B_0^{\mathbb N}({\mathbb
F}_p)= \varinjlim_{n} B(n,{\mathbb F}_p)$, where $B(n,{\mathbb
F}_p)$ is the group of $n$-by-$n$ upper-triangular matrices with
unities on the principal diagonal with entries from ${\mathbb
F}_p$. For the group $B_0^{\mathbb N}({\mathbb F}_p)$ we have the
following description
$$
B_0^{\mathbb N}({\mathbb F}_p)=\{I+\sum_{1\leq
k<n}x_{kn}E_{kn}\mid\,\,x_{kn}\in {\mathbb F}_p,\quad
x_{kn}=0\,\,\text{for large}\,\,n\}.
$$
Let $B^{\mathbb N}({\mathbb F}_p)$ be the group of all upper-triangular matrices:
 \begin{equation*}
 B^{\mathbb N}({\mathbb F}_p)=\{I+\sum_{1\leq
k<n}x_{kn}E_{kn}\mid\,\,x_{kn}\in {\mathbb F}_p\}.
\end{equation*}
We have the
following semi-direct product  $B^{\mathbb N}({\mathbb
F}_p)=B_m({\mathbb F}_p)\ltimes B^m({\mathbb F}_p)$, where
$B^m({\mathbb F}_p)$ is normal subgroup in $B^{\mathbb N}({\mathbb
F}_p)$ and
$$
B_m({\mathbb F}_p)=\{I+x\in B^{\mathbb N}({\mathbb F}_p) \mid
x=\sum_{m<k< n}x_{kn}E_{kn}\},
$$
$$
B^m({\mathbb F}_p)=\{I+x\in B^{\mathbb N}({\mathbb F}_p)\mid
x=\sum_{1\leq k\leq m,k<n}x_{kn}E_{kn}\},
$$
and we shall write $B^{\mathbb N}({\mathbb F}_p)\ni x=x_m\cdot
x^m\in B_m({\mathbb F}_p)\cdot B^m({\mathbb F}_p).$ We define the
space $X^m$ as the factor-space $X^m=B_m({\mathbb F}_p)\backslash
B^{\mathbb N}({\mathbb F}_p)\simeq B^m({\mathbb F}_p)$. The right
action $R_t$ of the group $B^{\mathbb N}({\mathbb F}_p)$ is
correctly defined on the factor-space $X^m$ by the formula $
R_t(x)=(xt^{-1})^m,\,\, t\in B^{\mathbb N}({\mathbb F}_p),\,\,x\in
B^m({\mathbb F}_p).$ We have
\begin{equation}\label{R_t(x)}
R_t(x)=xt^{-1},\,\,\text{if}\,\,t\in B^m({\mathbb
F}_p),\,\,\text{and}\,\,
R_t(x)=t_mxt^{-1},\,\,\text{if}\,\,t=t_mt^m\not\in B^m({\mathbb
F}_p).
\end{equation}
To prove (\ref{R_t(x)}) we get
$$
B_m({\mathbb F}_p)\cdot B^m({\mathbb F}_p)\ni x_m\cdot x^m
\stackrel{x\mapsto xt}{\mapsto } x_m\cdot x^mt_m\cdot
t^m=x_mt_m(t_m^{-1}x^mt_m)t^m
$$
hence, $(x^mt_mt^m)^m=t_m^{-1}x^mt$ and $(xt^{-1})^m=t_mxt^{-1}$
for $x\in B^m({\mathbb F}_p)$. We use relation
$(t^{-1})_m=(t_m)^{-1}$.
 The measure
$\mu_\alpha$ on the space $X^m$ is defined as infinite tensor
product
\begin{equation*}
\mu_\alpha=\otimes_{1\leq k\leq
m,k<n}\mu_{\alpha_{kn}}=\otimes_{k=1}^m\mu^k_\alpha,\,\,{\rm
where}\,\, \mu^k_\alpha:=\otimes_{n=k+1}^\infty\mu_{\alpha_{kn}}
\end{equation*}
of the probability measures $\mu_{\alpha_{kn}}$  on ${\mathbb
F}_p$ defined as follows:

${\mathbb F}_p\ni r\mapsto
\mu_{\alpha_{kn}}(r)=\alpha_{kn}(r)>0$ and $\sum_{r\in{\mathbb
F}_p}\alpha_{kn}(r)=1.$
\begin{lem}
\label{simFp} We have $\mu_\alpha^{R_t}\sim \mu_\alpha$ for all
$t\in B_0^{\mathbb N}({\mathbb F}_p).$
\end{lem}
Define the unitary representation
$T^{R,\mu_\alpha,m}:B_0^{\mathbb N}({\mathbb F}_p)\mapsto
U(L^2(X^m,\mu_\alpha))$ in a natural way, i.e., for $f\in
L^2(X^m,\mu_\alpha)$ set
\begin{equation}
\label{Rep-q-fin}
(T^{R,\mu_\alpha,m}_tf)(x)=(d\mu_\alpha(R_t^{-1}(x))/d\mu_\alpha(x))^{1/2}f(R_t^{-1}(x)),\,
t\in B_0^{\mathbb N}({\mathbb F}_p).
\end{equation}

\begin{co}
\label{t.irr,F_p} The quasiregular representation
$T^{R,\mu_\alpha,m}$ of the group $B_0^{\mathbb N }({\mathbb
F}_p)$ is irreducible if and only if conditions  1)--3) holds:
\par
1) $\mu_\alpha^{L_t}\perp \mu_\alpha\,\,\forall t\in B(m,{\mathbb
F}_p)\backslash{\{e\}},$
\par
2) the measure $\mu_\alpha$ is $G$-ergodic,
\par
3)  for the measure $\mu_\alpha=\otimes_{k=1}^m\mu^k_\alpha$ holds
$\mu^m_\alpha\perp\mu_{inv}^m$.
\par
4) Moreover, $T^{R,\mu_\alpha,m}\sim
T^{R,\mu_\beta,n}$
if and only if $m=n$ and $\mu_\alpha\sim\mu_\beta.$
\end{co}
\begin{rem}
In the case of the field $k={\mathbb R}$ and the measure being a
Gaussian product-measure, the irreducibilty holds if and only if
the condition 1) and 2) are valid (see \cite{Kos02.3,Kos03})
hence, the case $k={\mathbb F}_p$ is richer.
\end{rem}
The right action $R$ of the group $B_0^{\mathbb N}({\mathbb F}_p)$
on the space $X^m$ is given by the formula (\ref{R_t(x)}). The left
action $L$ of the group $B(m,{\mathbb F}_p)$ on the space $X^m$ is
as follows: $L_t(x)=tx,\,\,t\in B(m,{\mathbb F}_p),\,x\in X^m$.
Let us consider the case $p=2$ and $m=2$, i.e., the space $X^2$.
Set $E_{kn}(d)=I+dE_{kn}\in G,\,\,d\in {\mathbb F}_p$. We have
$$
E_{12}(d)x=\left(\begin{smallmatrix}
1&d\\
0&1
\end{smallmatrix}\right)
\left(\begin{smallmatrix}
1&x_{12}&x_{13}&\dots&x_{1n}&\dots \\
0& 1    &x_{23}&\dots&x_{2n}&\dots
\end{smallmatrix}\right)=
\left(\begin{smallmatrix}
1&x_{12}+d&x_{13}+dx_{23}&\dots &x_{1n}+dx_{2n}&\dots \\
0& 1    &x_{23}&\dots &x_{2n}&\dots
\end{smallmatrix}\right).
$$
For $t=I+E_{kn},\,\,1\leq k\leq n$ the right action is
$R_t(x)=xt^{-1}$ (see (\ref{R_t(x)}))
$$
\left(\begin{array}{cccccc}
1&x_{12}&x_{13}&\dots &x_{1n}&\dots \\
0& 1    &x_{23}&\dots &x_{2n}&\dots
\end{array}\right)
\stackrel{R^{-1}_{I+E_{1n}}}{\mapsto }
 \left(\begin{array}{cccccc}
1&x_{12}&x_{13}&\dots &x_{1n}+1&\dots \\
0& 1    &x_{23}&\dots &x_{2n}&\dots
\end{array}\right),
$$
$$
\left(\begin{array}{cccccc}
1&x_{12}&x_{13}&\dots &x_{1n}&\dots \\
0& 1    &x_{23}&\dots &x_{2n}&\dots
\end{array}\right)
\stackrel{R^{-1}_{I+E_{2n}}}{\mapsto }
 \left(\begin{array}{cccccc}
1&x_{12}&x_{13}&\dots &x_{1n}+x_{12}&\dots \\
0& 1    &x_{23}&\dots &x_{2n}+1&\dots
\end{array}\right),
$$
$$
\left(\begin{smallmatrix}
1&x_{12}&x_{13}&\dots &x_{1k}&\dots &x_{1n}&\dots \\
0& 1    &x_{23}&\dots &x_{2k}&\dots &x_{2n}&\dots
\end{smallmatrix}\right)
\stackrel{R^{-1}_{I+E_{kn}}}{\mapsto }
 \left(\begin{smallmatrix}
1&x_{12}&x_{13}&\dots &x_{1k}&\dots &x_{1n}+x_{1k}&\dots \\
0& 1    &x_{23}&\dots &x_{2k}&\dots &x_{2n}+x_{2k}&\dots
\end{smallmatrix}\right).
$$
Therefore, we have four actions to study:
\begin{equation}
\label{right.act1}
 R^{-1}_{I+E_{1n}}\,:\,x_{1n}\mapsto x_{1n}+1,\quad
(x_{1k},x_{1n})\!\mapsto\!(x_{1k},x_{1n}+x_{1k}),
\end{equation}
\begin{equation}
\label{right.act2}
 \left(\!\begin{array}{cc}
x_{12}&x_{1n}\\
1     &x_{2n}
\end{array}\!\right)\!\!
\stackrel{R^{-1}_{I+E_{2n}}}{\mapsto}\!\! \left(\!\begin{array}{cc}
x_{12}&x_{1n}+x_{12}\\
1     &x_{2n}+1
\end{array}\!\right),\,
 \left(\!\begin{array}{cc}
x_{1k}&x_{1n}\\
x_{2k}&x_{2n}
\end{array}\!\right)\!\!
\stackrel{R^{-1}_{I+E_{kn}}}{\mapsto}\!\!\left(\!\begin{array}{cc}
x_{1k}&x_{1n}+x_{1k}\\
x_{2k} &x_{2n}+x_{2k}
\end{array}\!\right),
\end{equation}
and
\begin{equation}
\label{left.ac.} L_{I+dE_{12}}:\,\,
 \left(\begin{array}{c}
x_{1n}\\
x_{2n}
\end{array}\right)\mapsto
\left(\begin{array}{c}
x_{1n}+dx_{2n}\\
x_{2n}
\end{array}\right),\quad d\in{\mathbb F}_p.
\end{equation}
Set
\begin{equation}
\label{H(inv),H(alpha)}
H_{inv}=L^2({\mathbb F}_p,\mu_{inv})\quad\text{and}\quad
H_\alpha=L^2({\mathbb F}_p,\mu_\alpha)
\end{equation}
where the normalized Haar
measure $\mu_{inv}$ on the additive group ${\mathbb F}_p$ is
defined by
\begin{equation}
\label{haar(F-p)} \mu_{inv}(r)=p^{-1},\,\, r\in{\mathbb F}_p,\quad
\text{and}\quad \mu_\alpha(r)=\alpha(r),\quad
\text{with}\quad\sum_{r\in {\mathbb F}_p}\alpha(r)=1.
\end{equation}
The operator $T_{inv}$ on the Hilbert space $H_{inv}$ associated
with the action $x\mapsto x-1$ on ${\mathbb F}_p$ is defined by
the following formula
{\small
$$
(T_{inv}f)(x)\!=\!\left(\frac{d\mu_{inv}(x-1)}{d\mu_{inv}(x)}\right)^{1/2}\!\!f(x-1)=f(x-1),\,
\,f(x)\!=\!(f_0,f_1,\dots ,f_{p-1})\!\in\!{\mathbb C}^p.
$$
}
Take the orthonormal basis (o.n.b.) in the space $H_{\alpha}$ as
follows:
\begin{equation}
\label{e^al_k-o.n.b.}( e^\alpha_k)_{k\in {\mathbb F}_p},\,\,\text{
where}\,\,e^\alpha_k=(e^\alpha_k(r))_{r\in {\mathbb
F}_p},\,\,e^\alpha_k(r)=(\alpha(r))^{-1/2}\delta_{k,r},\,\,\,\,\,k,r\in
{\mathbb F}_p.
\end{equation}
For $e_k(r)=(p)^{-1/2}\delta_{kr},\,k,r\in {\mathbb F}_p$ we get
$(Te_k)(r)=e_k(r-1)=e_{k+1}(r)$, so
$$
T(\sum_kf_ke_k)= \sum_{k\in {\mathbb F}_p}f_ke_{k+1}=\sum_{k\in
{\mathbb F}_p}f_{k-1}e_{k}\quad\text{hence,}\quad T=\sum_{r\in
{\mathbb F}_p}E_{r+1,r}.
$$
To define the corresponding operator $T_{\alpha}$ on the Hilbert
space $H_\alpha$ we use the following commutative diagram:
$$
\left.
\begin{array}{ccccc}
&H_\alpha &\stackrel{T_{\alpha}}{\rightarrow}&H_\alpha&\\
U_\alpha&\downarrow&&\downarrow&
U_\alpha\\
 &H_{inv}&\stackrel{T_{inv}}{\rightarrow} &H_{inv}&
\end{array}
\right.
$$
where $U_\alpha : H_\alpha\rightarrow H_{inv}$ is the isomorphisms
defined by
$$
U_\alpha= \left( d\mu_\alpha(x)/d\mu_{inv}(x) \right)^{1/2}= {\rm
diag}((p\alpha(0))^{1/2},(p\alpha(1))^{1/2},\dots
,(p\alpha(p-1))^{1/2}).
$$
Finally, the operator $T_{\alpha}$ is equal to
$T_{\alpha}=U_\alpha^{-1} T_{inv}U_\alpha$ hence, we have for $p=
2$
\begin{equation}
\label{T-alpha.p=2}
T_{\alpha}=\left(\begin{smallmatrix} \frac{1}{\sqrt{2\alpha(0)}}
&0\\
0& \frac{1}{\sqrt{2\alpha(1)}}
\end{smallmatrix}\right)
\left(\begin{smallmatrix}
0&1\\
1&0
\end{smallmatrix}\right)
\left(\begin{smallmatrix} \sqrt{2\alpha(0)}
&0\\
0& \sqrt{2\alpha(1)}
\end{smallmatrix}\right)=
\left(\begin{smallmatrix} 0& \sqrt{\frac{\alpha(1)}{\alpha(0)}}
\\
 \sqrt{\frac{\alpha(0)}{\alpha(1)}} &0
\end{smallmatrix}\right).
\end{equation}
For general $p$ we have in the basis $( e^\alpha_k)_{k\in {\mathbb
F}_p}$
\begin{equation}
\label{T-alpha}
 T_\alpha=
\left(\begin{smallmatrix}
 0&0&0&\dots &0& \sqrt{
\frac{\alpha(p-1)}{\alpha(0)}}\\
\sqrt{
\frac{\alpha(0)}{\alpha(1)}}&0&0&\dots &0&0\\
0& \sqrt{
\frac{\alpha(1)}{\alpha(2)}}&0&\dots &0&0\\
&&&\dots &\\
0&0&0&\dots & \sqrt{ \frac{\alpha(p-2)}{\alpha(p-1)}}&0
\end{smallmatrix}\right)
,\quad T_{inv}= \left(\begin{smallmatrix}
0&0&0&\dots &0& 1\\
1&0&0&\dots &0&0\\
0& 1&0&\dots &0&0\\
&&&\dots &\\
0&0&0&\dots & 1&0
\end{smallmatrix}\right).
\end{equation}
\subsection{ The Kakutani criterion}  We find the condition of orthogonality
$\mu_{\alpha}^{L_{I+dE_{12}}}\perp\mu_{\alpha},\,\,d\in{\mathbb
F}_p\setminus \{0\}$, using the {\it Kakutani criterion}
\cite{Kak48}. The {\it Hellinger  integral} $H(\mu,\nu)$ for two
measures $\mu$ and $\nu$ on the space $X$ is defined \cite{Kuo75}
as follows:
$$
H(\mu,\nu)=
\int_X\sqrt{\frac{d\mu(x)}{d\rho(x)}\frac{d\nu(x)}{d\rho(x)}}d\rho(x),
$$
where $\rho$ is some measure on $X$ such that both measures $\mu$
and $\nu$ are {\it absolutely continuous} with respect to the
measure $\rho$. For example, one can take
$\rho=\frac{1}{2}(\mu+\nu)$.

Let we have two probability measures $\mu_\alpha$ and $\mu_\beta$
 on the group ${\mathbb F}_p$ defined as  follows:
$\mu_\alpha(r)=\alpha(r),\,\,\sum_{r\in {\mathbb F}_p}\alpha(r)=1$
and $\mu_\beta(r)=\beta(r),\,\,\sum_{r\in {\mathbb
F}_p}\beta(r)=1$. The Hellinger  integral
$H(\mu_\alpha,\mu_\beta)$ for two measures $\mu_\alpha$ and $
\mu_\beta$ is given in this case by
$$
H(\mu_\alpha,\mu_\beta)= \int_{{\mathbb
F}_p}\sqrt{\frac{d\mu_\alpha(x)}{d\mu_{inv}(x)}\frac{d\mu_\beta(x)}{d\mu_{inv}(x)}}d\mu_{inv}(x)=
\sum_{r\in {\mathbb F}_p} \sqrt{\alpha(r) \beta(r)}.
$$
Let us consider two probability measures
$\mu_\alpha=\otimes_{n\in{\mathbb N}}\mu_{\alpha_n}$ and
$\mu_\beta=\otimes_{n\in{\mathbb N}}\mu_{\beta_n}$ defined on the
space $({\mathbb F}_p)^\infty= {\mathbb F}_p\times{\mathbb
F}_p\times...$ as the infinite tensor product, where
$\mu_{\alpha_n}$ and $\mu_{\beta_n},\,\,\, n\in{\mathbb N}$ are
probability measures defined on the space ${\mathbb F}_p$, as
before. The Hellinger integral $H(\mu_\alpha,\mu_\beta)$ for two
measures $\mu_\alpha$ and $ \mu_\beta$ is given in this case by
$$
H(\mu_\alpha,\mu_\beta)=\prod_{n\in{\mathbb
N}}H(\mu_{\alpha_n},\mu_{\beta_n})= \prod_{n\in{\mathbb
N}}\sum_{r\in {\mathbb F}_p} \sqrt{\alpha_n(r) \beta_n(r)}.
$$
We use the {\it notation} $\mu^f(\Delta)=\mu(f^{-1}(\Delta))$ for
a measure $\mu$ on the space $X$ and a measurable bijection
$f:X\to X$. For two measures $\mu_\alpha\otimes\mu_\beta$ and
$(\mu_\alpha\otimes\mu_\beta)^{L_{I+dE_{12}}}$ on ${\mathbb
F}_p\times{\mathbb F}_p$ where
$L_{I+dE_{12}}:(x,y)\mapsto(x+dy,y)$ (see (\ref{left.ac.})) we
have
\begin{equation}
\label{(k,r)^L}
 (\mu_\alpha\otimes\mu_\beta)^{L^{-1}_{I+dE_{12}}}(k,r)=
(\alpha(k)\beta(r))^{L^{-1}_{I+dE_{12}}}= \alpha(k+dr)\beta(r).
\end{equation}
Hence, we have for the Hellinger integral the following expression:
$$
H_{12}^d:=H\left((\mu_\alpha\otimes\mu_\beta)^{L^{-1}_{I+dE_{12}}},\mu_\alpha\otimes\mu_\beta\right)=
\sum_{r,k\in{\mathbb
F}_p}\sqrt{(\alpha(k)\beta(r))^{L^{-1}_{I+dE_{12}}}
\alpha(k)\beta(r)}
$$
$$
=\sum_{r,k\in{\mathbb F}_p}\sqrt{\alpha(k+dr)\beta(r)
\alpha(k)\beta(r)}= \sum_{r\in{\mathbb
F}_p}\beta(r)\sum_{k\in{\mathbb F}_p}\sqrt{\alpha(k+dr)\alpha(k)}.
$$
\begin{lem}
\label{l.perp2} For the measure $\mu_\alpha=\otimes_{1\leq k\leq m,k<n}
\mu_{\alpha_{kn}}$ on the space $X^m$ five following conditions
are equivalent:
\begin{align*}
1)&\,\,\mu_\alpha^{L_t}\perp \mu_\alpha, \forall\,t\in
B(m,{\mathbb F}_p)\setminus\{e\},\\
2)&\,\,(\mu_\alpha)^{L_{I+dE_{ls}}}\perp \mu_\alpha,\,\,\forall\,d
\in{\mathbb F}_p\setminus \{0\},\,\,1\leq l<s\leq m,\\
3)&\,\,(\mu_\alpha)^{L_{I+E_{ls}}}\perp \mu_\alpha,\quad 1\leq l<s\leq m,\\
4)&\,\,\Pi^{L,d}_{ls}(\mu_\alpha)=\prod_{n=s+1}^\infty H_{n,ls}=
\prod_{n=s+1}^\infty  \sum_{r\in{\mathbb
F}_p}\alpha_{sn}(r)\sum_{k\in{\mathbb
F}_p}\sqrt{\alpha_{ln}(k+dr)\alpha_{ln}(k)}=0,\\
5)&\,\,S^{L,d}_{ls}(\mu_\alpha) =\sum_{n=s+1}^\infty
 \sum_{r\in{\mathbb F}_p}
\alpha_{sn}(r)\Big(1-\sum_{k\in{\mathbb F}_p
\setminus\{0\}}\sqrt{\alpha_{ln}(k+dr)\alpha_{ln}(k)}\Big)=\infty.
\end{align*}
\end{lem}
\begin{pf}
Obviously $1)\Rightarrow 2)\Rightarrow 3)\Rightarrow 4)$. We show
 that   $4)\Leftrightarrow 5)$.  The implication
$5)\Rightarrow 1)$ will follow from the irreducibility that we
prove later.

We show that $2)\Leftrightarrow 3)$. Indeed, since
$(\mu_\alpha)^{L_{I+dE_{ls}}}$ and $\mu_\alpha$ are product
measures, by Kakutani criterion, we conclude  that
$(\mu_\alpha)^{L_{I+dE_{ls}}}$ and $\mu_\alpha$ are orthogonal or
equivalent. It is sufficient to show that $3)$ implies $2)$ for
all $d\in{\mathbb F}_p^*={\mathbb F}_p\setminus\{0\}$. Let us suppose the opposite, i.e., that
for some  $d\in{\mathbb F}_p^*:={\mathbb F}_p\setminus\{0\}$ holds
$(\mu_\alpha)^{L_{I+dE_{ls}}}\sim \mu_\alpha$. Since ${\mathbb
F}_p^*$ is a multiplicative group there exists an inverse
$a=d^{-1}\in{\mathbb F}_p^*$. For this element $a$ we then get
$$
\mu_\alpha \sim
(\mu_\alpha)^{L^a_{I+dE_{ls}}}=(\mu_\alpha)^{L_{(I+dE_{ls})^a}}=
(\mu_\alpha)^{L_{I+adE_{ls}}}= (\mu_\alpha)^{L_{I+E_{ls}}}.$$
This contradicts with $3)$.
We have $4)\Leftrightarrow 5)$ since
\begin{align*}
\sum_{r\in{\mathbb
F}_p}\alpha_{sn}(r)\sum_{k\in{\mathbb
F}_p}\sqrt{\alpha_{ln}(k+r)\alpha_{ln}(k)}=&\,\alpha_{sn}(0)\\
\!+\!\!\!\sum_{r\in{\mathbb
F}_p\setminus\{0\}}\alpha_{sn}(r)\sum_{k\in{\mathbb
F}_p}\sqrt{\alpha_{ln}(k+r)\alpha_{ln}(k)}=&\,
1-\!\!\!\!\sum_{r\in{\mathbb F}_p\setminus\{0\}}\alpha_{sn}(r)\!+\!\!\!\!
\sum_{r\in{\mathbb F}_p\setminus\{0\}}\!\!\alpha_{sn}(r)\\
\times\sum_{k\in{\mathbb
F}_p}\sqrt{\alpha_{ln}(k+r)\alpha_{ln}(k)}=1-\sum_{r\in{\mathbb
F}_p\setminus\{0\}} \alpha_{sn}(r)&
\Big(1-\sum_{k\in{\mathbb F}_p}\sqrt{\alpha_{ln}(k+r)\alpha_{ln}(k)}\Big).
\end{align*}
\qed\end{pf}
\begin{rem}
\label{(3):only-k=m}
If $\mu^l_\alpha\sim \mu^l_{inv}$ for some $l,\,\,1\leq l<m$ then by Lemma~\ref{l.perp2}, 4)  we get
$(\mu_\alpha)^{L_{I+dE_{ls}}}\sim \mu_\alpha$ for $l<s\leq m$ hence,  the representation is reducible.
Therefore, only one condition from the list of conditions 3) in Conjecture~\ref{t.irr,F_p} is independent,
namely: $\mu^m_\alpha\perp \mu^m_{inv}$!
\end{rem}
\begin{pf}
If we replace the factor $\mu^l_\alpha$ in the expression for the
measure $\mu_\alpha=\otimes_{k=1}^m\mu^k_\alpha$ by $\mu^l_{inv}$
we get the equivalent measure $\mu_{\alpha'}$ and the
representation $T^{R,\mu_{\alpha'},m}$ equivalent with the initial
one $T^{R,\mu_{\alpha},m}$. For this measure we have
$(\mu_{\alpha'})^{L_{I+dE_{ls}}}\sim \mu_{\alpha'}$ for $s:l<s\leq
m$. Indeed, in this case we have\\
$\sum_{k\in{\mathbb F}_p}\sqrt{\alpha_{tn}(k+dr)\alpha_{tn}(k)}\!=\!1$
hence, $\Pi^{L,d}_{ls}(\mu)\!=\!1$. The representation $T^{R,\mu_{\alpha'},m}$ is reducible in this case,
since the operator $T^{L,\mu_{\alpha'},m}_{I+dE_{ts}}$ generated by the
transformation $L_{I+dE_{ls}}$ is well defined and commutes with the
representation $T^{R,\mu_{\alpha'},m}$.
\qed\end{pf}
%
{\bf Examples.} 1) In the particular case $p=2$ we have
$$
H_{n,12}=
\beta(0)+2\beta(1)\sqrt{ \alpha(0)\alpha(1)}=
\alpha_{2n}(0)+2\alpha_{2n}(1)\sqrt{ \alpha_{1n}(0)\alpha_{1n}(1)}
$$
$$
=1-\alpha_{2n}(1)(1-2\sqrt{ \alpha_{1n}(0)\alpha_{1n}(1)}).
$$
Hence, for $X=X^2$ and ${\mathbb F}_2$ we have
\begin{equation}
\label{Pi^L12F2} \Pi^{L,1}_{12}(\mu_\alpha)=\prod_{n=3}^\infty
H_{n,12}= \prod_{n=3}^\infty\left( 1-\alpha_{2n}(1)(1-2\sqrt{
\alpha_{1n}(0)\alpha_{1n}(1)}
\right).
\end{equation}
We see that $\Pi^{L,1}_{12}(\mu_\alpha)=0$ if and only if
$S^{L,1}_{12}(\mu_\alpha)=\infty$ where
$$
S^{L}_{12}(\mu_\alpha):=
S^{L,1}_{12}(\mu_\alpha)=\sum_{n=3}^\infty
\alpha_{2n}(1)\left(1-2\sqrt{\alpha_{1n}(0)\alpha_{1n}(1)}\right).
$$
2) For $X=X^3$ and ${\mathbb F}_2$ we have
\begin{gather*}
\Pi^{L,1}_{12}(\mu_\alpha)=\prod_{n=3}^\infty H_{n,12}=
\prod_{n=3}^\infty\left( \alpha_{2n}(0)+2\alpha_{2n}(1)\sqrt{
\alpha_{1n}(0)\alpha_{1n}(1)}\right),\\
\Pi^{L,1}_{13}(\mu_\alpha)=\prod_{n=4}^\infty H_{n,13}=
\prod_{n=4}^\infty\left( \alpha_{3n}(0)+2\alpha_{3n}(1)\sqrt{
\alpha_{1n}(0)\alpha_{1n}(1)}\right),\\
\Pi^{L,1}_{23}(\mu_\alpha)=\prod_{n=4}^\infty H_{n,23}=
\prod_{n=4}^\infty\left( \alpha_{3n}(0)+2\alpha_{3n}(1)\sqrt{
\alpha_{2n}(0)\alpha_{2n}(1)}\right).
\end{gather*}
3) For ${\mathbb F}_3$ and $X^2$ we have
\begin{align*}
H^1_{12}=&\beta(0)\left(\alpha(0)+\alpha(1)+\alpha(2)\right)\\
+&\beta(1)\left(\sqrt{\alpha(1)\alpha(0)}+\sqrt{\alpha(2)\alpha(1)}+
\sqrt{\alpha(0)\alpha(2)}\right)\\
+&\beta(2)\left(\sqrt{\alpha(2)\alpha(0)}+\sqrt{\alpha(0)\alpha(1)}+
\sqrt{\alpha(1)\alpha(2)}\right)\\
=&\beta(0)+(\beta(1)+\beta(2))\left(\sqrt{\alpha(2)\alpha(0)}+\sqrt{\alpha(0)\alpha(1)}+
\sqrt{\alpha(1)\alpha(2)}\right),
\end{align*}
hence, for ${\mathbb F}_3$ and $X^2$ we have
\begin{equation}
\label{Pi^L12F3}
 \Pi^{L,1}_{12}(\mu_\alpha)=\prod_{n=2}^\infty H_{n,12}= \prod_{n\in{\mathbb N},n>2}\left(
\alpha_{2n}(0)+(\alpha_{2n}(1)+\alpha_{2n}(2))\right.
\end{equation}
$$
\times\left.
\left(\sqrt{\alpha_{1n}(2)\alpha_{1n}(0)}+\sqrt{\alpha_{1n}(0)\alpha_{1n}(1)}+
\sqrt{\alpha_{1n}(1)\alpha_{1n}(2)}\right) \right).
$$
We study  first the condition 1) of Lemma~\ref{l.perp2}.
\begin{lem}
\label{mu^k-per-mu^i} The following three conditions are
equivalent:
\begin{align*}
1)\, &\mu^l_{\alpha}\perp\mu_{inv}^l,\quad 1\leq l\leq m,\\
2)\,
&\Pi_{ll}(\mu_\alpha)=\prod_{n=l+1}^\infty\frac{1}{p}\Big(1+\sum_{r\in{\mathbb
F}_p}\sum_{k\in{\mathbb
F}_p\setminus\{0\}}\sqrt{\alpha_{ln}(k)\alpha_{ln}(k+r)} \Big)=0,\\
3)\,&S_{ll}^L(\mu_\alpha)= \sum_{n=l+1}^\infty\sum_{r\in{\mathbb
F}_p\setminus\{0\}} \Big(1-\sum_{k\in{\mathbb
F}_p}\sqrt{\alpha_{ln}(k+r)\alpha_{ln}(k)}\Big)=\infty.
\end{align*}
\end{lem}
{\bf Particular cases}. 1) $p=2$ and $m=1.$ We have only one condition:
$$
S_{11}^L(\mu_\alpha)=
\sum_{n=2}^\infty\left(1-2\sqrt{\alpha_{1n}(0)\alpha_{1n}(1)}\right)=\infty.
$$
2) The case $p=2$ and $m\in {\mathbb N}.$ We have the following conditions for $1\leq k\leq n\leq m$:
\begin{equation}
\label{S^L(kk)}
S_{kn}^L(\mu_\alpha)=\infty,\quad
\text{where}\quad S_{kk}^L(\mu_\alpha)=
\sum_{r=k+1}^\infty\Big(1-2\sqrt{\alpha_{kr}(0)\alpha_{kr}(1)}\Big),
\end{equation}
\begin{equation}
\label{S^L(kn)}
 S_{kn}^L(\mu_\alpha)=
\sum_{r=n+1}^\infty\alpha_{nr}(1)\Big(1-2\sqrt{\alpha_{kr}(0)\alpha_{kr}(1)}\Big),\,\,\,\,k<n.
\end{equation}
\begin{rem}
\label{Haar(Z_p)} The conditions 3) of the Conjecture~\ref{co.Qreg(F_p)-irr} mean
the following. The space $X^{(k)}\!=\!\prod_{n=k+1}^\infty ({\mathbb
F}_p)_n$ is isomorphic to the set ${\mathbb Z}_p=\{x\in {\mathbb
Q}_p:\Vert x\Vert_p\leq 1\}$ of {\rm entire $p$-adic numbers of the
field ${\mathbb Q}_p$ of all $p$-adic numbers}  since ${\mathbb
Z}_p$ has the following description:
 ${\mathbb Z}_p=\{\sum_{n=0}^\infty a_np^n\mid a_n\in {\mathbb F}_p\}$. The
measure $\mu_{inv}^k$ on $X^{(k)}$ is the Haar measure on
${\mathbb Z}_p$ under this identification.
\end{rem}
\begin{rem} The lemma analogous to Lemma \ref{l.perp2}  holds in the case when we replace
the field $k={\mathbb F}_p$ by the ring $k={\mathbb Z}$. The
measure $\mu_\alpha$ on ${\mathbb Z}$ is defined by
$\mu_\alpha(r)=\alpha(r)>0,\,\,r\in {\mathbb Z}$ such that
$\sum_{r\in{\mathbb Z}}\alpha(r)=1$. The corresponding conditions
are the following:
\begin{align*}
1)&\,\,\mu_\alpha^{L_t}\perp \mu_\alpha, \forall\,t\in
B(m,{\mathbb Z})\setminus\{e\},\\
2)&\,\,(\mu_\alpha)^{L_{I+dE_{ls}}}\perp \mu_\alpha,\,\,\forall\,d
\in{\mathbb Z}\setminus \{0\},\,\,1\leq l<s\leq m,\\
3)&\,\,\Pi^{L,d}_{ls}(\mu_\alpha)=\prod_{n=s+1}^\infty H_{n,ls}=
\prod_{n=s+1}^\infty \sum_{r\in{\mathbb
Z}}\alpha_{sn}(r)\sum_{k\in{\mathbb
Z}}\sqrt{\alpha_{ln}(k+dr)\alpha_{ln}(k)}=0.
\end{align*}
\end{rem}

\subsection{ Fourier transform.}
 Let us consider an additive group of the
field ${\mathbb F}_p$. The Haar measure $\mu_{inv}$ on ${\mathbb
F}_p$ is defined by $\mu_{inv}(r)=1/p,\,\,r\in{\mathbb F}_p$. The
set of unitary characters $\chi_R(r),\,\,R\in {\mathbb F}_p$,
 are defined as follows:
\begin{equation}
{\mathbb F}_p\ni r\mapsto \chi_R(r)=\exp \frac{2\pi iRr}{p}\in
S^1.
\end{equation}
The Fourier transform $F$ is defined on the space
$H_{inv}=L^2({\mathbb F}_p,\mu_{inv})$ by the formula
\begin{equation}
(Ff)(R):=\tilde f(R):=\sqrt{p}\int_{{\mathbb F}_p}f(x){\overline
\chi}_R(x)d\mu_{inv}(x)= \frac{1}{\sqrt{p}}\sum_{r\in{\mathbb
F}_p}f(r)\exp\Big( -\frac{2\pi iRr}{p}\Big).
\end{equation}
The operator $F$ is a unitary operator on the space $L^2({\mathbb
F}_p,\mu_{inv})$.
\begin{lem}
\label{l.FTF^-1}
 The image ${\tilde T}_{inv}=FT_{inv}F^{-1}$ of the operator $T_{inv}$ with
respect to the Fourier transform is defined by
\begin{equation}
\label{FTF^-1}
({\tilde T}_{inv}\tilde f)(R)=\exp\Big(- \frac{2\pi
iR}{p}\Big)\tilde f(R),
\end{equation}
i.e., ${\tilde T}_{inv}\!=\!{\rm diag}(\bar{\chi}_1(0),\bar{\chi}_1(1),\dots
,\bar{\chi}_1(p-1))\!=\!{\rm diag}(1,\lambda,\dots
,\lambda^{p-1}),\,\,\lambda\!=\!\bar{\chi}_1(1)$.
\end{lem}
\begin{pf} Indeed, by (\ref{T-alpha})we have $T_{inv}:f(x)\mapsto f(x-1)$ hence,
$$
\tilde f(R)=\frac{1}{\sqrt{p}}\sum_{r\in{\mathbb
F}_p}f(r)\exp\Big(- \frac{2\pi i Rr}{p}\Big)
\stackrel{T_{inv}}\mapsto \frac{1}{\sqrt{p}}\sum_{r\in{\mathbb
F}_p}f(r-1)\exp\Big(-\frac{2\pi iRr}{p}\Big)=
$$
$$
\frac{1}{\sqrt{p}}\sum_{s\in{\mathbb F}_p}f(s)\exp\Big(-\frac{2\pi
iR(s+1)}{p}\Big)=\exp\Big(-\frac{2\pi iR}{p}\Big)\tilde f(R).
$$
\qed\end{pf}
To define the Fourier transform $F_{\alpha}:H_\alpha\rightarrow
H_{\tilde\alpha}$,  where the measure $\tilde\alpha$ on ${\mathbb F}_p$ is defined by (\ref{F(a)1=1}), we use the following commutative diagram:
$$
\left.
\begin{array}{ccccc}
&L^2({\mathbb F}_p,\mu_\alpha)
&\stackrel{F_{\alpha}}{\rightarrow}&L^2(
{\mathbb F}_p,\mu_{\tilde\alpha})&\\
U_\alpha&\downarrow&&\downarrow&
U_{\tilde\alpha}\\
 &L^2({\mathbb
F}_p,\mu_{inv})&\stackrel{F}{\rightarrow} &L^2({\mathbb
F}_p,\mu_{inv})&
\end{array}
\right.
$$
where $(U_\alpha f)(r)=(p\alpha(r))^{1/2}f(r)$. We have
$F_{\alpha}=U_{\tilde\alpha}^{-1}FU_{\alpha}$ hence (compare with the case
of the Fourier transform (\ref{Fourier.1}) in $L^2(\mathbb R,\mu)$ defined
below)
\begin{equation}
\label{Fourier.p}
\left(F_{\alpha}f\right)(R)= \frac{1}{\sqrt{\tilde\alpha(R)p}}
\sum_{r\in{\mathbb F}_p}\sqrt{\alpha(r)}f(r)\exp\left( -\frac{2\pi
iRr}{p}\right).
\end{equation}
\begin{rem}
\label{r.case-R}
Let us denote  by $F^\mu_{kn}$ the one-dimensional Fourier
transform corresponding to the measure $\mu_{kn}$ (see
\cite{Kos00m} formula (6) and (7))
$$
\left.
\begin{array}{ccccc}
&L^2({\mathbb
R},\mu_{kn})&\stackrel{F^\mu_{kn}}{\rightarrow}&L^2({\mathbb
R},\tilde\mu_{kn})&\\
U_{kn}^\mu=\left(\frac{d\mu_{kn}(x)}{dx}\right)^{1/2}
&\downarrow&&\downarrow&
U_{kn}^{\tilde\mu}=\left(\frac{\tilde\mu_{kn}(y)}{dy}\right)^{1/2}\\
 &L^2({\mathbb R},dx)&\stackrel{F}{\rightarrow} &L^2({\mathbb R},dy).&
\end{array}
\right.
$$
By definition, $F^\mu_{kn}=(U_{kn}^{\tilde\mu})^{-1}FU_{kn}^\mu$,
where
$$
(Ff)(y)=\frac{1}{\sqrt{2\pi}}\int_{\mathbb R}f(x)\exp(iyx)dx,
$$
so we have
\begin{equation}
\label{Fourier.1} (F^\mu_{kn}f)(y)=
\left(\frac{d\tilde\mu_{kn}(y)}{dy}\right)^{-1/2}
\frac{1}{\sqrt{2\pi}}\int_{\mathbb
R}f(x)\exp(iyx)\left(\frac{d\mu_{kn}(x)}{dx}\right)^{1/2}dx.
\end{equation}
In the case where the Fourier transform $F\mu_{kn}^{1/2}$ of the
function $\mu_{kn}^{1/2}$,
is positive, we define the density
\begin{equation}
\label{(F(mu^1/2))^2-mesure}
\tilde\mu_{kn}(y):=\mid(F\mu_{kn}^{1/2})(y)\mid^2
\end{equation}
of the corresponding measure
$d\tilde\mu_{kn}(y):=\tilde\mu_{kn}(y)dy.$
\end{rem}

\begin{rem}
\label{S^L(12):p} We compare the conditions $\mu_\alpha^{L_t}\perp
\mu_\alpha, \forall\,t\in B(m,{\bf k})\setminus\{e\}$ for ${\bf
k}={\mathbb R}$ and ${\bf k}={\mathbb F}_p$  when $m=2$.
\par (a) In
the case ${\bf k}={\mathbb R}$ we have
$$
X^2=\left(\begin{array}{cccccc}
1&x_{12}&x_{13}&\dots &x_{1n}&\dots \\
0& 1    &x_{23}&\dots &x_{2n}&\dots
\end{array}\right),
$$
$$
d\mu_{b}^m(x)=\otimes_{1\leq k\leq 2,\,k<n}\sqrt{\frac{b_{kn}}{\pi}}
\exp(-b_{kn}x_{kn}^2)dx_{kn}.
$$
For the operator $U_{12}^L(t):=T^{L,\mu_b,2}_{I+tE_{12}}=\exp(itA^L_{12}),\,\,t\in {\mathbb R}$ acting on $L^2(X^2,\mu_b)$,
where $A^L_{12}=D_{12}+\sum_{k=3}^\infty x_{2k}D_{1k}$ we have
$$
U_{12}^L(t)=
\exp\big(it(D_{12}+\sum_{k=3}^\infty
x_{2k}D_{1k})\big)=\exp(itD_{12})
\prod_{k=3}^\infty\exp(itx_{2k}D_{1k})= :\prod_{k=2}^\infty
U_{k}^L(t),
$$
$$
S^L_{12}(\mu)=\Vert A^L_{12}{\bf 1}\Vert^2=\Vert
\Big(D_{12}+\sum_{k=3}^\infty x_{2k}D_{1k}\Big){\bf 1}\Vert^2=
\Vert D_{12}{\bf 1}\Vert^2+\sum_{k=3}^\infty \Vert x_{2k}{\bf
1}\Vert^2\Vert D_{1k}{\bf 1}\Vert^2_{H_{1k}}
$$
$$=\!
\Vert (F_{12}D_{12}){\bf 1}\Vert^2+\sum_{k=3}^\infty \Vert
x_{2k}{\bf 1}\Vert^2\Vert (F_{1k}D_{1k}){\bf 1}\Vert^2_{{\tilde
H}_{1k}}\!=\!\Vert iy_{12}{\bf 1}\Vert^2+ \sum_{k=3}^\infty \Vert
x_{2k}{\bf 1}\Vert^2\Vert iy_{1k}{\bf 1}\Vert^2_{{\tilde H}_{1k}}
$$
$$
=
\frac{b_{12}}{2}+\sum_{k=3}^\infty\frac{1}{2b_{2k}}\frac{b_{1k}}{2}.
$$
\par (b)  In the case ${\bf k}={\mathbb F}_2$ we have $X^2$ and
$\mu_\alpha=\otimes_{1\leq k\leq 2,k<n}\mu_{\alpha_{kn}}$. Using
(\ref{left.ac.}) we get $\left(\begin{smallmatrix}
x_{12}\\
1
\end{smallmatrix}\right)\mapsto
\left(\begin{smallmatrix}
x_{1n}+1\\
1
\end{smallmatrix}\right)
$
and
 $\left(\begin{smallmatrix}
x_{1n}\\
x_{2n}
\end{smallmatrix}\right)\mapsto
\left(\begin{smallmatrix}
x_{1n}+x_{2n}\\
x_{2n}
\end{smallmatrix}\right).
$
Hence, the corresponding operator $U_{12}^L(t):=T^{L,\mu_b,2}_{I+tE_{12}},\,\,t\in {\mathbb F}_2$
acting on $L^2(X^2,\mu_\alpha)$, has the following form for $t=1$:
$$
U^L_{12}(1)=\otimes_{k=3}^\infty U_k(1),\,\text{where}\,\,U_2-1=(T_{\alpha_{12}}-1),\,\,
U_k-1\!=\!P_{1k}^{(1)}\otimes(T_{\alpha_{1k}}-1),\,k\geq 3.
$$
To get the two latter expressions we use (\ref{T-alpha}) and (\ref{T2n(1)-1.0}). Therefore, we get
$$
\sum_{k=3}^\infty\Vert (U_k-1){\bf
1}\Vert^2=\Vert
(T_{\alpha_{12}}-1){\bf 1}\Vert_{H_{12}}^2+\sum_{k=3}^\infty\Vert P_{2k}^{(1)}{\bf 1} \Vert^2\Vert
(T_{\alpha_{1k}}-1){\bf 1}\Vert_{H_{1k}}^2
$$
$$
=\Vert ({\tilde
T}_{\alpha_{12}}-1){\bf 1}\Vert_{{\tilde
H}_{12}}^2+
\sum_{k=2}^\infty\Vert P_{2k}^{(1)}{\bf 1} \Vert^2\Vert ({\tilde
T}_{\alpha_{1k}}-1){\bf 1}\Vert_{{\tilde
H}_{1k}}^2
$$
$$
=
\left(1-2\sqrt{
\alpha_{12}(0)\alpha_{12}(1)}\right)+
\sum_{n=3}^\infty \alpha_{2n}(1)\left(1-2\sqrt{
\alpha_{1n}(0)\alpha_{1n}(1)}\right).
$$
\end{rem}
\begin{rem}
\label{r.F(a)1=1}
Suppose that the square of the Fourir transform of the square root of the measure $\mu_\alpha$ on ${\mathbb F}_p$ is again a measure on ${\mathbb F}_p$. Compare with  the case of the field $\mathbb R$, Remark~\ref{r.case-R}.
The latter condition is equivalent with the following one: $F_\alpha {\bf 1}={\bf 1}$ that means by (\ref{Fourier.p} the following:
\begin{equation}
\label{F(a)1=1}
\left(F_{\alpha}{\bf 1}\right)(R)= \frac{1}{\sqrt{\tilde\alpha(R)p}}
\sum_{r\in{\mathbb F}_p}\sqrt{\alpha(r)}\exp\left( -\frac{2\pi
iRr}{p}\right)=1,\quad R\in {\mathbb F}_p.
\end{equation}
For $p=2$ we get $\left(F_{\alpha}{\bf 1}\right)(0)=\left(F_{\alpha}{\bf 1}\right)(1)=1$ or
$$
(\sqrt{\alpha(0)}+\sqrt{\alpha(1)})/
\sqrt{2\tilde\alpha(0)}
=1,\quad
(\sqrt{\alpha(0)}-\sqrt{\alpha(1)})/
\sqrt{2\tilde\alpha(1)}=1,
$$
hence, we get if $\alpha(0)>\alpha(1)$
$$
\tilde\alpha(0)=
(1+2\sqrt{\alpha(0)\alpha(1)})/2,\quad
\tilde\alpha(1)=
(1-2\sqrt{\alpha(0)\alpha(1)})/2.
$$
For $p=3$ we get $\left(F_{\alpha}{\bf 1}\right)(0)=\left(F_{\alpha}{\bf 1}\right)(1)=\left(F_{\alpha}{\bf 1}\right)(2)=1$ if
$\sqrt{\alpha(0)}>(\sqrt{\alpha(1)}+\sqrt{\alpha(2)})/2\,\,\,$ and $\,\,\,\alpha(1)=\alpha(2)$.
Indeed, we get
{\small
\begin{align*}
\left(F_{\alpha}{\bf 1}\right)(0)=
(\sqrt{\alpha(0)}+\sqrt{\alpha(1)}+\sqrt{\alpha(2)})/\sqrt{3\tilde\alpha(0)}=1,\\
\left(F_{\alpha}{\bf 1}\right)(1)=
(\sqrt{\alpha(0)}+\sqrt{\alpha(1)}\exp( -2\pi
i/3)+\sqrt{\alpha(2)}\exp( -4\pi i/3)
)\sqrt{3\tilde\alpha(1)}=1,\\
\left(F_{\alpha}{\bf 1}\right)(2)=
(\sqrt{\alpha(0)}+\sqrt{\alpha(1)}\exp( -4\pi i/3)+\sqrt{\alpha(2)}\exp( -8\pi
i/3))/\sqrt{3\tilde\alpha(2)}=1
\end{align*}
}
hence,
\begin{align*}
\sqrt{\alpha(0)}+\big(-\frac{1}{2}+i\frac{\sqrt{3}}{2}\big)\sqrt{\alpha(1)}+
\big(-\frac{1}{2}-i\frac{\sqrt{3}}{2}\big)\sqrt{\alpha(2)}=\\
\sqrt{\alpha(0)}-\frac{1}{2}(\sqrt{\alpha(1)}+\sqrt{\alpha(2)})+i\frac{\sqrt{3}}{2}
(\sqrt{\alpha(1)}-\sqrt{\alpha(2)})>0.
\end{align*}
For general prime $p$ we get $\sum_{R\in {\mathbb F}_p}\sqrt{\alpha(r)}\exp(-2\pi Rr/p)=1$ for $R\in {\mathbb F}_p$.
\end{rem}
\subsection{ Maximal abelian subalgebra and a simple spectrum}
\begin{df}
\label{max.ab.sub}
An abelian  subalgebra  of a von Neumann algebra ${\mathfrak A}$  is called  maximal if it is not properly included in any other such subalgebra of ${\mathfrak A}$.
\end{df}

Consider a finite-dimensional Hilbert space $H={\mathbb C}^n$ with the standard scalar product $(x,y)=\sum_{k=1}^nx_{k}\bar{y}_k.$
\begin{df}
\label{Sp-simple}
A spectrum $Sp(A)$ of an operator $A$ in an $n$-dimensional Hilbert space $H$ we call {\rm simple} if $Sp(A)$ consists of $n$ distinct eigenvalues.
\end{df}
\index{spectrum!simple}
\begin{lem}
\label{max-ab=simpl(Sp)}
A von Neumann algebara ${\mathfrak A}=L^\infty(A)$ generated by a diagonal operator $A={\rm diag}(\lambda_k)_{k=1}^n$ in $H={\mathbb C}^n$ is maximal abelian if and only if
the spectrum of $A$ is simple. In addition, $L^\infty(A)=\{P(A)\mid {\rm ord}\,P\leq n-1\}$ where $P(x)=\sum_{k=0}^{n-1}a_kx^k,\,\,a_k\in {\mathbb C}$.
\end{lem}
\begin{pf}
We know that for a von Neumann algebra ${\mathfrak A}=L^\infty(A)$ holds  $({\mathfrak A})'={\mathfrak A}$ therefore, $L^\infty(A)=(A)''$. We show that
$$
(L^\infty(A))'=(A)'=\{{\rm diag}(b_k)_{k=1}^n \mid b=(b_k)_{k=1}^n\in{\mathbb C}^n\}.
$$
Indeed, let $[A,B]=0$ where $B=(b_{km})_{k,m=1}^n$ then
$$
\lambda_kb_{km}=b_{km}\lambda_m,\quad\text{for all}\quad k\not=m.
$$
Therefore, $b_{km}=0$ for $k\not=m$ since $\lambda_k\not=\lambda_m$.
By the same arguments we show that
$$
L^\infty(A)=(A)''=\{{\rm diag}(a_k)_{k=1}^n \mid a=(a_k)_{k=1}^n\in{\mathbb C}^n\}.
$$
When, for example, $A={\rm diag}(\lambda_1,\lambda_2,\lambda_2)$ with $\lambda_1\not=\lambda_2$ then $L^\infty(A)\not=(L^\infty(A))'$ since
$$
L^\infty(A)=\left\{ \left(\begin{smallmatrix}
a_{11}&0&0\\
0&a_{22}&0\\
0&0&a_{33}
\end{smallmatrix}\right)\mid a_{kn}\in {\mathbb C}
\right\},
\quad
(L^\infty(A))'=\left\{ \left(\begin{smallmatrix}
a_{11}&0&0\\
0&a_{22}&a_{23}\\
0& a_{32}&a_{33}
\end{smallmatrix}\right)\mid a_{kn}\in {\mathbb C}
\right\}.
$$
\qed\end{pf}
Denote by $L^\infty(T_\alpha)$ the von Neumann algebra of operators acting on $H_\alpha\simeq {\mathbb C}^p$ generated by operator
$T_\alpha$, defined by (\ref{T-alpha}), i.e.,
$L^\infty(T_\alpha)=(T_\alpha)''$.
\begin{lem}
\label{t.(T)'} The von Neumann algebra
$L^\infty(T_\alpha)$ {\rm is a maximal abelian subalgebra} in $B(H_\alpha)$, i.e.,
$(L^\infty(T_\alpha))'=L^\infty(T_\alpha).$ In addition, $L^\infty(A)=\{P(A)\mid {\rm ord}\,P\leq p-1\}$.
\end{lem}
\index{subalgebra!von Neuman!maximal abelian}
%
\begin{pf}
 By Lemma~\ref{FTF^-1} $T_{inv}\sim {\tilde T}_{inv}$ and
${\tilde T}_{inv}$ is the diagonal operator with different
eigenvalues:
$
Sp({\tilde T}_{inv})=\{\exp(-ik\pi/p)\mid k\in {\mathbb F}_p\}.
$
\qed\end{pf}
\begin{lem}
\label{t.(x)'} The von Neumann algebra
$L^\infty(x)$ generated by the operator $x\!=\!{\rm diag}(k)_{k\in{\mathbb F}_p}$ {\rm is a maximal abelian subalgebra} in $B(H_\alpha)$, i.e.,
$(L^\infty(x))'\!=\!L^\infty(x).$
\end{lem}
\begin{pf}Since the spectrum $Sp(x)=\{k\mid k\in{\mathbb F}_p\}$ of the operator $x$ is simple  the proof follows from Lemma~\ref{max-ab=simpl(Sp)}.
\qed\end{pf}
%
\section{The Laplace operator and the commutant description}
\label{s8.irr}
\subsection{The Laplace operator and the irreducibility}
%
 For approximationof an operator $x_{kn}$ defined by (\ref{X(kn)}) we shall use the well known
result (see for example \cite{BecBel61}, Chap. I, \S 52)
\begin{equation*}
 \min_{x\in{\mathbb
R}^n}\Big(\sum_{k=1}^na_kx_k^2\mid \sum_{k=1}^n x_k=1\Big)=
\Big(\sum_{k=1}^n \frac{1}{a_k}\Big)^{-1},\,\,
a_k>0,\,\,k=1,2,\dots ,n.
\end{equation*}
We use the same result in a slightly different form with
$b_k\not=0,\,\,k=1,2,\dots ,n$
\begin{equation}
\label{8.min2}
 \min_{x\in{\mathbb
R}^n}\Big(\sum_{k=1}^na_kx_k^2\mid \sum_{k=1}^n x_kb_k=1\Big)=
\Big(\sum_{k=1}^n \frac{b_k^2}{a_k}\Big)^{-1}.
\end{equation}
The minimum is realized for
$x_k=\frac{b_k}{a_k}\Big(\sum_{k=1}^n\frac{b_k^2}{a_k}\Big)^{-1}$.

Let $p=2$ and $X=X^1$. Set
\begin{equation}
\label{A(kn)}
A_{\alpha_{1n}}:=T_{\alpha_{1n}}-1=
\left(\begin{smallmatrix}
 -1&
\sqrt{\frac{\alpha_{1n}(1)}{\alpha_{1n}(0)}}
\\
 \sqrt{\frac{\alpha_{1n}(0)}{\alpha_{1n}(1)}} &-1
\end{smallmatrix}\right).
\end{equation}
\begin{rem}
\label{approx-I} In order to approximate the operator
$x_{1k}:={\rm diag}(0,1)$ acting on $H_{1k}$ (see (\ref{H(kn)}))
by linear combinations of $T_{kn}-I=x_{1k}
\otimes(T_{\alpha_{1n}}-I)$ (see (\ref{T2n(1)-1.0})) it is
sufficient to approximate the identity operator $Id=I$ by linear
combinations
$\sum_{n=2}^{N}t_nA_{\alpha_{1n}}=\sum_{n=2}^{N}t_n(T_{\alpha_{1n}}-1)$
as $N\to\infty$ (see Lemma~\ref{1in<2>}  below).
\end{rem}
{\bf Notations}. Set ${\bf 1}={\bf 1}\otimes {\bf 1}\otimes {\bf
1}\otimes \dots \in L^2(X^m,\mu_\alpha)=\otimes_{1\leq k\leq
m,k<n}L^2({\mathbb F}_p,\mu_{\alpha_{kn}})$, where ${\bf
1}\!=\!(1,1,\dots ,1)\in\!L^2({\mathbb F}_p,\mu_{\alpha_{kn}})$
and let $c_{1n}=2\sqrt{\alpha_{1n}(0)\alpha_{1n}(1)}$. As before, for
$f\in L^1(X,\mu)$ we use the following notation $Mf=\int_X
f(x)d\mu(x)$ and  let $\langle f_n\!\mid\! n\in {\mathbb N}\rangle$
be a  closed subspace generated by the set of vectors
$(f_n)_{ n\in{\mathbb N}}$ in a space $H$.
\begin{rem}
Obviously, two series with positive $a_n,b_n$ are
{\it equivalent}:
\begin{equation}
\label{8.sim1}
\sum_n\frac{a_n}{a_n+b_n}\,\,\sim\,\,\sum_n\frac{a_n}{b_n},
\end{equation}
i.e., they are
simultaneously convergent or divergent.
\end{rem}
\begin{lem}
\label{1in<2>} Three following conditions are equivalent if $p=2$:
\begin{align*}
(i)&\,\,{\bf 1}\in \langle A_{\alpha_{1n}}{\bf 1}\mid n\geq
2\rangle,\\
(ii)&\,\, S_{11}^L(\mu_\alpha)=
\sum_{n=2}^\infty\left(1-2\sqrt{\alpha_{1n}(0)\alpha_{1n}(1)}\right)=:
\sum_{n=2}^\infty\left(1-c_{1n}\right)= \infty,\\
(iii)&\,\, \mu_\alpha\perp \mu_{inv}.
\end{align*}
\end{lem}
\begin{pf}
To prove $(i)\Leftrightarrow (ii)$  set $\xi_n=
A_{\alpha_{1n}}{\bf 1}$. We have
$$
(T_{\alpha_{1n}}{\bf 1},{\bf
1})\!=\!\sqrt{\frac{\alpha_{1n}(1)}{\alpha_{1n}(0)}}\alpha_{1n}(0)
+\sqrt{\frac{\alpha_{1n}(0)}{\alpha_{1n}(1)}}\alpha_{1n}(1)\!=\!2\sqrt{\alpha_{1n}(0)
\alpha_{1n}(1)}\!=\!c_{1n},
$$
$$
M\xi_n=MA_{\alpha_{1n}}{\bf 1}=( A_{\alpha_{1n}}{\bf 1},{\bf
1})=((T_{\alpha_{1n}}-1){\bf 1},{\bf 1}) =c_{1n}-1,
$$
$$
\Vert \xi_n\Vert^2=\Vert A_{\alpha_{1n}}{\bf
1}\Vert^2=(A_{\alpha_{1n}}{\bf 1},A_{\alpha_{1n}}{\bf 1})=
((T_{\alpha_{1n}}-1){\bf 1},(T_{\alpha_{1n}}-1){\bf 1})
$$
$$
=2-2(T_{\alpha_{1n}}{\bf 1},{\bf 1})=2(1-c_{1n}).
$$
Finally, we have
$$
M\xi_n=-\left(1-2\sqrt{\alpha_{1n}(0)\alpha_{1n}(1)}\right),\quad
\Vert \xi_n\Vert^2=
2\left(1-2\sqrt{\alpha_{1n}(0)\alpha_{1n}(1)}\right).
$$
If we take $(t_n)_n$ such that $\sum_{n=2}^{N+2}t_nM\xi_n\!=\!1$ we obtain (since
$\xi_n\!\!-\!\!M\xi_n\perp$\\ $\xi_m\!-\!M\xi_m$ for $n\not=m$ )
$$
\Vert \Big(\sum_{n=2}^{N+2}t_nA_{\alpha_{1n}}-1\Big){\bf
1}\Vert^2= \Vert \sum_{n=2}^{N+2}t_n\left(A_{\alpha_{1n}}-M\xi_n
\right){\bf 1}\Vert^2=
$$
$$
\sum_{n=2}^{N+2}t_n^2\Vert\left(A_{\alpha_{1n}}-M\xi_n
\right){\bf 1}\Vert^2=
\Vert \sum_{n=2}^{N+2}t_n (\xi_n\!-\!M\xi_n)\Vert^2
$$
$$
=\sum_{n=2}^{N+2}t_n^2\left(\Vert
\xi_n\Vert^2-\mid M\xi_n\mid^2\right).
$$
Using  (\ref{8.min2}) for $b_n=M\xi_n$ and $a_n=\Vert
\xi_n\Vert^2-\mid M\xi_n\mid^2$ we conclude that
$$
\min_{t\in{\mathbb R}^N}\Big(\Vert
\big(\sum_{n=2}^{N+2}t_nA_{\alpha_{1n}}-1\big){\bf 1}\Vert^2\mid \sum_{n=2}^{N+2}t_nM\xi_n=1\Big)=(
S^L_{11,N}(\mu_\alpha))^{-1}
$$
where
$$
S_{11,N}^L(\mu_\alpha)= \sum_{n=2}^{N+2}\frac{b_n^2}{a_n}=
\sum_{n=2}^{N+2}\frac{\mid M\xi_n\mid^2}{\Vert \xi_n\Vert^2-\mid
M\xi_n\mid^2}
\stackrel{(\ref{8.sim1})}{\sim}
\sum_{n=2}^{N+2}\frac{\mid M\xi_n\mid^2}{\Vert
\xi_n\Vert^2}=
$$
$$
\sum_{n=2}^{N+2}\frac{\left(1-2\sqrt{\alpha_{1n}(0)\alpha_{1n}(1)}\right)^2}{
2\left(1-2\sqrt{\alpha_{1n}(0)\alpha_{1n}(1)}\right)}=
\frac{1}{2}\sum_{n=2}^{N+2}
\left(1-2\sqrt{\alpha_{1n}(0)\alpha_{1n}(1)}\right).
$$
 To prove $(ii)\Leftrightarrow (iii)$ we have
$$
H(\mu_{\alpha_{1n}},\mu_{inv})=\int_{{\mathbb F}_2} \sqrt{
\frac{d\mu_{\alpha_{1n}}(x)} {d\mu_{inv}(x)}} d\mu_{inv}(x)=
\frac{\sqrt{\alpha_{1n}(0)}+\sqrt{\alpha_{1n}(1)}}{\sqrt{2}}.
$$
Hence, we get
$$
H(\mu_{\alpha},\mu_{inv})=\prod_{n=2}^\infty
H(\mu_{\alpha_{1n}},\mu_{inv})=\prod_{n=2}^\infty\frac{
\sqrt{\alpha_{1n}(0)}+\sqrt{\alpha_{1n}(1)}}{\sqrt{2}}.
$$
By the Kakutani criterion, we conclude that
$\mu_\alpha\perp\mu_{inv}$ if and only if
$$
H(\mu_{\alpha},\mu_{inv})=\prod_{n=2}^\infty\frac{
\sqrt{\alpha_{1n}(0)}+\sqrt{\alpha_{1n}(1)}}{\sqrt{2}}=0
\Leftrightarrow
$$
$$
\prod_{n=2}^\infty\Big(\frac{
\sqrt{\alpha_{1n}(0)}+\sqrt{\alpha_{1n}(1)}}{\sqrt{2}}\Big)^2=0
\Leftrightarrow \prod_{n=2}^\infty
\frac{1}{2}\left(1+2\sqrt{\alpha_{1n}(0)\alpha_{1n}(1)}\right)=0
$$
$$
\Leftrightarrow
S_{11}^L(\mu_\alpha)=
\sum_{n=2}^\infty\left(1-2\sqrt{\alpha_{1n}(0)\alpha_{1n}(1)}\right)=\infty.
$$
\qed\end{pf}
\begin{lem}
\label{l.cntr1} Denote by  $T_n=T_{\alpha_{1n}},\,\,n\geq 2$.
The following strong limit of operators $\Delta_1= {\rm
s.}\lim_{k\to\infty} \prod_{n=2}^k\frac{I+T_n}{2} $ is correctly defined if and only if the following equivalent conditions hold:
$$
\mu_\alpha\sim\mu_{inv} \Leftrightarrow \prod_{n=2}^\infty\frac{
\sqrt{\alpha_{1n}(0)}+\sqrt{\alpha_{1n}(1)}}{\sqrt{2}}>0\Leftrightarrow
S_{11}^L(\mu_\alpha)<\infty.
$$
\end{lem}
\begin{pf}
We have
$$
\Vert 1/2(I+T_n){\bf 1}\Vert^2=1/4\Vert(I+T_n){\bf
1}\Vert^2=1/4[({\bf 1},{\bf 1})+2(T_n{\bf 1},{\bf 1})+(T_n{\bf
1},T_n{\bf 1})]
$$
$$
=1/2(1+(T_n{\bf 1},{\bf 1}))
\stackrel{(\ref{T-alpha})}{=} 1/2(1+
2\sqrt{\alpha_{1n}(0)\alpha_{1n}(1)})=\Big( \frac{
\sqrt{\alpha_{1n}(0)}+\sqrt{\alpha_{1n}(1)}}{\sqrt{2}} \Big)^2.
$$
Hence,
\begin{equation}
\label{H^2(.,.),p=2}
\lim_{k\to\infty} \Vert\prod_{n=3}^k\frac{I+T_n}{2}{\bf 1}\Vert^2=
\Big( \prod_{n=3}^\infty\frac{
\sqrt{\alpha_{1n}(0)}+\sqrt{\alpha_{1n}(1)}}{\sqrt{2}} \Big)^2=H^2(\mu_{\alpha},\mu_{inv}).
\end{equation}
\qed\end{pf}
Consider the space $X^1$, i.e., the case $m=1$, and fix a prime $p>2$. For
the operator $T_{\alpha}$ we have (see (\ref{T-alpha}))
$$
T_{\alpha}=
\left(\begin{smallmatrix}
 0&0&0&\dots &0& \sqrt{
\frac{\alpha(p-1)}{\alpha(0)}}\\
\sqrt{
\frac{\alpha(0)}{\alpha(1)}}&0&0&\dots &0&0\\
0& \sqrt{
\frac{\alpha(1)}{\alpha(2)}}&0&\dots &0&0\\
&&&\dots &\\
0&0&0&\dots & \sqrt{ \frac{\alpha(p-2)}{\alpha(p-1)}}&0
\end{smallmatrix}\right)=
\sum_{k\in{\mathbb
F}_p}\sqrt{\frac{\alpha(k)}{\alpha(k+1)}}E_{k+1,k}.
$$
\begin{rem}
\label{r.gen:C(T)-p} To guess the expression for the right version of the operator
$A_{\alpha_{1n}}$ defined by (\ref{A(kn)}), if $p>2$, we observe that
$A_{\alpha_{1n}}=T_{\alpha_{1n}}-1= \sum_{r\in{\mathbb
F}_2}T_{\alpha_{1n}}^r-2$ for $p=2$. Hence, it is natural to
replace the operator $A_{\alpha_{1n}}$ in the case  $p=2$ by the
expression $C(T_{\alpha_{1n}})-p=\sum_{r\in{\mathbb F}_p}T^{r}_{\alpha_{1n}}-p$ in the case  of an arbitrary $p$.
\end{rem}
\index{$C(T)$}
\begin{lem}
\label{1in<p>} Three following conditions are equivalent for an
arbitrary $p$:
\par
$(i)\quad {\bf 1}\in \langle \Big(\sum_{r\in{\mathbb
F}_p}T^{p-r}_{\alpha_{1n}}-p\Big){\bf 1}\mid n\geq 2\rangle,$
\par
$(ii)\quad S_{11}^L(\mu_\alpha)=\infty,$ where
$$
S_{11}^L(\mu_\alpha):=\sum_{n=2}^\infty
 \sum_{r\in{\mathbb F}_p\setminus\{0\}}
\Big(1-\sum_{k\in{\mathbb
F}_p}\sqrt{\alpha_{1n}(k+r)\alpha_{1n}(k)}\Big)=
$$
$$
\sum_{n=2}^\infty \Big(p-\sum_{k,s\in{\mathbb F}_p}
\sqrt{\alpha_{1n}(s)\alpha_{1n}(k)}\Big)=p\sum_{n=2}^\infty\Big[
1-\Big(\sum_{k\in{\mathbb F}_p} \sqrt{\frac{\alpha_{1n}(k)}{p}}
\Big)^2 \Big],
$$
\par
$(iii)\quad \mu_\alpha\perp\mu_{inv}.$
\end{lem}
\begin{pf}
To prove $(i)\Leftrightarrow (ii)$ we have
\begin{equation}
\label{T,T^r=}
 T_{\alpha}=\sum_{k\in{\mathbb
F}_p}\sqrt{\frac{\alpha(k)}{\alpha(k+1)}}E_{k+1,k},\quad
T_{\alpha}^r=\sum_{k\in{\mathbb
F}_p}\sqrt{\frac{\alpha(k)}{\alpha(k+r)}}E_{k+r,k}.
\end{equation}
Using the change of variables $k+p-r=s,\,k=s-p+r=s+r\,\,{\rm
mod}\,p\,\,$ we get:
$$
T_{\alpha}^{p-r}=\sum_{k\in{\mathbb
F}_p}\sqrt{\frac{\alpha(k)}{\alpha(k+p-r)}}E_{k+p-r,k}=
\sum_{s\in{\mathbb
F}_p}\sqrt{\frac{\alpha(s+r)}{\alpha(s)}}E_{s,s+r}.
$$
Finally, we have
\begin{equation}
\label{xi_n}
C(T_{\alpha})= \sum_{r\in{\mathbb F}_p}T_{\alpha}^{p-r}=
\sum_{r\in{\mathbb F}_p}\sum_{s\in{\mathbb
F}_p}\sqrt{\frac{\alpha(s+r)}{\alpha(s)}}E_{ss+r}
=\sum_{k\in{\mathbb F}_p}\sum_{r\in{\mathbb F}_p}
\sqrt{\frac{\alpha(r)}{\alpha(k)}}E_{kr}.
\end{equation}
Set $\xi_n=\xi_{\alpha_{1n}}=(C(T_{\alpha_{1n}})-p){\bf 1}=\left(\sum_{r\in{\mathbb
F}_p}T^{p-r}_{\alpha_{1n}}-p\right){\bf 1}$. We get
$$
M\xi_{\alpha}=\big((C(T)-p){\bf 1}, {\bf 1}\big)=(C(T){\bf 1}, {\bf 1})-p=
$$
$$
\stackrel{(\ref{xi_n})}{=} \sum_{k\in{\mathbb
F}_p}\Big(\sum_{r\in{\mathbb F}_p}
\sqrt{\alpha(r)/\alpha(k)}
\Big)
\alpha(k)-p
=\sum_{r,k\in{\mathbb F}_p}\sqrt{\alpha(r)\alpha(k)}-p.
$$
To calculate $\Vert\xi_{\alpha}\Vert^2$ we get using (\ref{xi_n})
$$
\Vert\xi_{\alpha}\Vert^2= \sum_{k\in{\mathbb
F}_p}\Big|\sum_{r\in{\mathbb
F}_p}\sqrt{\frac{\alpha(r)}{\alpha(k)}}-p\,\Big|^2\alpha(k)=
$$
$$
= \sum_{k\in{\mathbb F}_p}\Big( \sum_{r\in{\mathbb
F}_p}\alpha(r)+p^2\alpha(k)+\sum_{r,s\in{\mathbb F}_p,\,r\not=s}
\sqrt{\alpha(r)\alpha(s)}-2p\sum_{r\in{\mathbb
F}_p}\sqrt{\alpha(r)\alpha(k)}\Big)
$$
$$
=\sum_{k\in{\mathbb F}_p}\Big( p^2\alpha(k)+\sum_{r,s\in{\mathbb
F}_p} \sqrt{\alpha(r)\alpha(s)}-2p\sum_{r\in{\mathbb
F}_p}\sqrt{\alpha(r)\alpha(k)}\Big)
$$
$$
= p^2+p\sum_{r,s\in{\mathbb F}_p} \sqrt{\alpha(r)\alpha(s)}-
2p\sum_{r,k\in{\mathbb F}_p} \sqrt{\alpha(r)\alpha(k)}
$$
$$
=p^2-p\sum_{r,k\in{\mathbb F}_p} \sqrt{\alpha(r)\alpha(k)}=
p\Big(p-\sum_{r,k\in{\mathbb F}_p}
\sqrt{\alpha(r)\alpha(k)}\Big).
$$
Finally, we get
\begin{equation*}
M\xi_{\alpha}\!=\!-\Big(p-\sum_{r,k\in{\mathbb
F}_p}\sqrt{\alpha(r)\alpha(k)}\Big),\quad
\Vert\xi_{\alpha}\Vert^2\!=\!p\Big(p\!-\!\sum_{r,k\in{\mathbb F}_p}
\sqrt{\alpha(r)\alpha(k)}\Big).
\end{equation*}
If we take $\sum_{n=2}^{N+2}t_nM\xi_n\!=\!1$ we obtain (since
$\xi_n\!-\!M\xi_n\!\perp\!\xi_m\!-\!M\xi_m$ for $n\not=m$)
$$
\Vert \Big(\sum_{n=2}^{N+2}t_n \Big(C(T_{\alpha_{1n}})-p\Big) -1\Big){\bf 1}\Vert^2= \Vert
\sum_{n=2}^{N+2}t_n\Big(\Big(C(T_{\alpha_{1n}})-p\Big)-M\xi_n \Big){\bf 1}\Vert^2
$$
$$
=\Vert \sum_{n=2}^{N+2}t_n (\xi_n\!-\!M\xi_n)\Vert^2
=\sum_{n=2}^{N+2}t_n^2\Big(\Vert \xi_n\Vert^2-\mid
M\xi_n\mid^2\Big).
$$
Using  (\ref{8.min2}) for $b_n=M\xi_n$ and $a_n=\Vert
\xi_n\Vert^2-\mid M\xi_n\mid^2$ we conclude that
$$
\min_{t\in{\mathbb R}^N}\Big(\Vert
\Big[\sum_{n=2}^{N+2}t_n\big(
C(T_{\alpha_{1n}})-p
\big)-I\Big]{\bf 1}\Vert^2\mid \sum_{n=2}^{N+2}t_nM\xi_n=1\Big)=
(S_{11,N}^L(\mu))^{-1},
$$
where
$$
S_{11,N}^L(\mu)= \sum_{n=2}^{N+2}\frac{\mid M\xi_n\mid^2}{\Vert
\xi_n\Vert^2-\mid M\xi_n\mid^2}\stackrel{(\ref{8.sim1})}\sim \sum_{n=2}^{N+2}\frac{\mid
M\xi_n\mid^2}{\Vert \xi_n\Vert^2}=
$$
$$
\sum_{n=2}^{N+2}\frac{ \Big(p-\sum_{r,k\in{\mathbb F}_p}
\sqrt{\alpha_{1n}(r)\alpha_{1n}(k)}\Big)^2}{
p\left(p-\sum_{r,k\in{\mathbb F}_p}
\sqrt{\alpha_{1n}(r)\alpha_{1n}(k)}\right)}=
\frac{1}{p}\sum_{n=2}^{N+2} \Big(p-\sum_{r,k\in{\mathbb F}_p}
\sqrt{\alpha_{1n}(r)\alpha_{1n}(k)}\Big).
$$
To prove $(ii)\Leftrightarrow (iii)$ we have $H(\mu_{\alpha},\mu_{inv})=$
\begin{equation}
\label{H(.,.),m=1}  \prod_{n=2}^\infty
H(\mu_{\alpha_{1n}},\mu_{inv})= \prod_{n=2}^\infty \int_{{\mathbb
F}_p} \sqrt{ \frac{d\mu_{\alpha_{1n}}(x)} {d\mu_{inv}(x)} } d\mu_{inv}(x)=
 \prod_{n=2}^\infty
 \sum_{k\in{\mathbb
F}_p}\sqrt{\frac{\alpha_{1n}(k)}{p}}.
\end{equation}
So, by Kakutani's criterion, we conclude that
$\mu_\alpha\perp\mu_{inv}$ if and only if
$$
H(\mu_{\alpha},\mu_{inv})=\prod_{n=2}^\infty \sum_{k\in{\mathbb
F}_p}\sqrt{\frac{\alpha_{1n}(k)}{p}} =0 \Leftrightarrow
\prod_{n=2}^\infty \Big( \sum_{k\in{\mathbb
F}_p}\sqrt{\frac{\alpha_{1n}(k)}{p}}\Big)^2=0.
$$
We note that
$$
\Big( \sum_{k\in{\mathbb F}_p}\sqrt{\alpha_{1n}(k)}\Big)^2=
\sum_{k\in{\mathbb F}_p}\sum_{r\in{\mathbb
F}_p}\sqrt{\alpha_{1n}(k)\alpha_{1n}(r)} = \sum_{r\in{\mathbb
F}_p}\sum_{k\in{\mathbb
F}_p}\sqrt{\alpha_{1n}(k)\alpha_{1n}(k+r)}=
$$
$$
1+\sum_{r\in{\mathbb F}_p\setminus\{0\}}\sum_{k\in{\mathbb
F}_p}\sqrt{\alpha_{1n}(k)\alpha_{1n}(k+r)}=: 1+\sum_{r\in{\mathbb
F}_p\setminus\{0\}}c_{1n}(r),
$$
where $c_{1n}(r):=\sum_{k\in{\mathbb
F}_p}\sqrt{\alpha_{1n}(k)\alpha_{1n}(k+r)}$. Finally,
$\mu_\alpha\perp\mu_{inv}$ if and only if
$$
\prod_{n=2}^\infty\frac{1}{p}\Big(1+\sum_{r\in{\mathbb
F}_p\setminus\{0\}}\sum_{k\in{\mathbb
F}_p}\sqrt{\alpha_{1n}(k)\alpha_{1n}(k+r)}
\Big)=0\Leftrightarrow
$$
$$
S_{11}^L(\mu_\alpha)= \sum_{n\in{\mathbb N},n>2}
 \sum_{r\in{\mathbb F}_p\setminus\{0\}}
\Big(1-\sum_{k\in{\mathbb
F}_p}\sqrt{\alpha_{1n}(k+r)\alpha_{1n}(k)}\Big)=\infty.
$$
\qed\end{pf}
\begin{lem}
\label{l.cntr1,p} Set $C(T)=\sum_{r\in {\mathbb F}_p}T^r$ and $T_n=T_{\alpha_{1n}}$.
The following strong limit $ \Delta_1={\rm
s.}\lim_{k\to\infty} \prod_{n=2}^kp^{-1}C(T_n), $ is correctly defined if and only if three equivalent conditions hold:
\begin{align*}
(i)&\quad\mu_\alpha\sim\mu_{inv},\\
(ii)&\,\,\,H(\mu_\alpha,\mu_{inv})=\prod_{n=2}^\infty\sum_{r\in {\mathbb F}_p}\sqrt{
\frac{\alpha_{1n}(r)}{p}}>0,\\
(iii)&\quad S_{11}^L(\mu_\alpha)<\infty.
\end{align*}
\end{lem}
\begin{pf}
Using (\ref{xi_n}) we have $\Vert C(T_n){\bf 1}\Vert^2\!=\!$
\begin{equation}
\label{normC(T)}
 \sum_{k\in {\mathbb
F}_p}\Big(\sum_{r\in {\mathbb
F}_p}\sqrt{\frac{\alpha_{1n}(r)}{\alpha_{1n}(k)}}\,\Big)^2
\alpha_{1n}(k)\!=\!p\Big(\sum_{r\in {\mathbb
F}_p}\sqrt{\alpha_{1n}(r)}\Big)^2\!=\!p^2
\Big(H(\mu_{\alpha_{1n}},\mu_{inv})\Big)^2,
\end{equation}
hence, by (\ref{H(.,.),m=1}) we get
\begin{equation}
\label{H^2(.,.)=norm(D),p} \lim\limits_{k\to\infty}\Vert \prod_{n=2}^k
p^{-1}C(T_n){\bf 1} \Vert^2= \prod_{n=2}^\infty
\Big(\sum_{r\in{\mathbb F}_p} \sqrt{\frac{\alpha_{1n}(r)}{p}}\,
\Big)^2=\Big(H(\mu_{\alpha},\mu_{inv})\Big)^2.
\end{equation}
\qed\end{pf}
Using Lemmas \ref{1in<p>} and \ref{l.cntr1,p}  we conclude that
\begin{lem}
\label{4cond.} The following four conditions are equivalent for
the measure $\mu_\alpha$ on the space $X^1$:
\begin{align*}
(i)&\,\,\, \mu_\alpha\sim\mu_{inv},\\
(ii)&\,\,\, S^L_{11}(\mu_\alpha)<\infty,\\
(iii)&\,\,\,{\bf 1}\not\in \langle \left(C(T_{1n})-p\right){\bf 1}\mid n\geq 2\rangle,\\
(iv)&\,\,\, \text{there exist a non trivial limit}\quad
\Delta_{1}:=\lim_{n\to\infty}\prod_{k=2}^np^{-1}C(T_{1k}).
\end{align*}
\end{lem}
\begin{rem}
\label{r.(3) -p-irr}
Using Lemma~\ref{4cond.} we conclude that {\rm condition 3) of
Conjecture~\ref{co.Qreg(F_p)-irr} are necessary for the irreducibility of the
representation} $T^{R,\mu_\alpha,m}$.
\end{rem}
Consider the measure $\mu_\alpha=\otimes_{k=1}^m\mu_\alpha^k$ on the space $X^m$ and the representation
$T^{R,\mu_\alpha,m}$.
\begin{lem}
\label{Lap-exist}
If $\mu_\alpha^k\sim \mu_{inv}^k$ for some
$1\leq k\leq m$ then the Laplace operator
$$
\Delta_k=s.\lim_{r\to\infty}\prod_{n=k+1}^rp^{-1}C(T_{kn}(k))
$$
 is  well defined  and commutes with the representation
$T^{R,\mu_\alpha,m}$. In particular, if
$\mu_\alpha\sim\mu_{inv}=\otimes_{k=1}^m\mu_{inv}^k$ then the
Laplace operator $\Delta^{(m)}:=\Delta_1\Delta_2\dots\Delta_m$ is
well defined and commutes with the representation.
\end{lem}
\index{operator!Laplace}
\begin{pf}
The operator $\Delta_l$ is well defined  by analogue of Lemma
\ref{l.cntr1,p}. To prove that $\Delta_l$ commutes with the
representation, i.e., $[\Delta_l,T^{R,\mu_\alpha,m}_t]=0$ for all $t\in
B_0^{\mathbb N}({\mathbb F}_p)$ it is sufficient to prove
commutation $[\Delta_l,T_{kk+1}]=0$ for all $k\in {\mathbb N}$
since the subgroups $E_{kk+1}(t)=I+tE_{kk+1},\,\,t\in {\mathbb
F}_p,\,\,k\in {\mathbb N}$ generate all the group $B_0^{\mathbb
N}({\mathbb F}_p)$.

In the case $m=1$ we prove the commutation relations
$[\Delta_{1},T_{kk+1}]=0$ for all $k\in {\mathbb N}$. The latter
relation follows from $[C(T_{12}),T_{12}]=0$ that is evident,
since $TC(T)=C(T)T=C(T)$ for $T$ such that $T^p=I$, and the
relation $[C(T_{1k})C(T_{1k+1}),T_{kk+1}]\!=\!0$ $k\!>\!2$. We
prove more general relations:
\begin{equation}
\label{[CC,T]=0}
 [C(T_{1k})C(T_{1m}),T_{km}]\!=\!0\quad\text{for}\quad 1<k<m.
\end{equation} We have
$$
C(T_{1k})C(T_{1m})T_{km}=C(T_{1k})T_{km}C(T_{1m})= \sum_{r\in{\mathbb
F}_p}T_{1k}^rT_{km}C(T_{1m})\stackrel{(\ref{[kr,rs]})}{=}
$$
$$
T_{km}\sum_{r\in{\mathbb
F}_p}T_{1k}^rT_{1m}^rC(T_{1m})=T_{km}\sum_{r\in{\mathbb
F}_p}T_{1k}^rC(T_{1m})=T_{km}C(T_{1k})C(T_{1m}).
$$
For the general $m$
we show that $[\Delta_l,T_{kk+1}]=0$ for
$1\leq l\leq m$ and $k\in {\mathbb N}$. First, using (\ref{T(kn)=otimes}) we conclude that
$[\Delta_l,T_{kk+1}]=0$ for $1\leq k<l$. Further, we conclude that $[\Delta_l,T_{kk+1}]=0$ for $k\geq l$ by analogy with the
relation $[\Delta_1,T_{kk+1}]=0$ for $k\geq 1$.
\qed\end{pf}
\subsection{Commutant of the von Neumann algebra ${\mathfrak A}^m$, case $m=1$}
\label{8.4.2Centre}
In this subsection we explain how the Laplace operator $\Delta_k$
(see (\ref{Del,Frob})) in the commutant $({\mathfrak A}^m)'$ was
found. Let
$${\mathfrak A}^{m}=\big(T^{R,\mu_\alpha,m}_t\mid t\in G\big)''=\left(T_{rk}\mid 1\leq r<k\right)''
$$
be the von-Neumann algebra generated by the representation
$T^{R,\mu_\alpha,m}$ acting in the space $H^m=L^2(X^m,\mu_\alpha)$
and let ${\mathfrak A}^{m,n}$ be its von-Neumann subalgebra,
\begin{equation}
\label{A^m,n}
{\mathfrak A}^{m,n}=\left(T_{rk}\mid 1\leq r\leq
m,\,r<k\leq n\right)'',
\end{equation}
 where $T_{kn}$ are defined by (\ref{T_{kn},T_{kn}(r)}). We have
${\mathfrak A}^{m}=$
$(\bigcup_{n>m}{\mathfrak A}^{m,n})''$. We would
like to describe the {\it commutant} $({\mathfrak
A}^{m})'=\bigcap_{n>m}({\mathfrak A}^{m,n})'$ of the von Neumann algebra
${\mathfrak A}^{m}$
First, we shall do this for $m=1.$

To describe $({\mathfrak A}^{m,n})'$ it is sufficient to consider
the invariant measures $\mu_{inv}$ since in the finite-dimensional
space $X^{m,n}$ (see below (\ref{X^(m,n)})) all considered measures are equivalent. Set
$H_{kn}=L^2({\mathbb F}_p,\mu_{inv}^{kn})$, where
$\mu_{inv}(r)=\mu_{inv}^{kn}(r)=p^{-1}$. For $m=1$ and $n=3$, we
denote
$$
X=\{(x_{12},x_{13})\mid x_{kn}\in{\mathbb F}_p\},\quad H=H_{12}\otimes H_{13}=L^2({\mathbb
F}_p,\mu_{inv})\otimes L^2({\mathbb F}_p,\mu_{inv}).
$$
We fix the basis $(e_k)_{k\in {\mathbb F}_p}$ in $L^2({\mathbb
F}_p,\mu_{inv})$, where
$e_k(x)\!=\!p^{1/2}\delta_{k,x},\,k,x\!\in\!{\mathbb F}_p$ (see
\ref{e^al_k-o.n.b.}), i.e.,
$$
\quad \quad \quad e_0=(p^{1/2},0,\dots
,0),\,\,\,\,e_2=(0,p^{1/2},0,\dots ,0),\,\,\dots
,\,\,e_{p-1}=(0,\dots 0,p^{1/2}).
$$
Fix $p=2$.  The operators $T_{12}$ and $T_{13}$ act as follows on the spaces $H_{12}$ and $H_{13}$:
$$
T_{12}=T_{13}=T:=T_{inv}=
\left(\begin{smallmatrix}
0&1\\
1&0
\end{smallmatrix}\right).
$$
In the space $H_{12}\otimes H_{13}$ we get $T_{12}=T_{inv}\otimes
I,\,\, T_{13}=I\otimes T_{inv}$ (see Remark \ref{otimes=1}), i.e.,
$$
T_{12}=
\left(\begin{smallmatrix}
0&1\\
1&0
\end{smallmatrix}\right)\otimes
\left(\begin{smallmatrix}
1&0\\
0&1
\end{smallmatrix}\right),\quad
T_{13}=
\left(\begin{smallmatrix}
1&0\\
0&1
\end{smallmatrix}\right)\otimes
\left(\begin{smallmatrix}
0&1\\
1&0
\end{smallmatrix}\right).
$$
The basis in the space $L^2({\mathbb F}_p,\mu_{inv})\otimes
L^2({\mathbb F}_p,\mu_{inv})$ is $(e_{kr}:=e_k\otimes
e_r)_{k,r\in{\mathbb F}_p}$. Let us fix the {\it lexicographic
order} on the set $(k,r)_{k,r\in{\mathbb F}_p}$. For $p=2$ the
basis $(e_{kr})_{k,r}$ in $H_{12}\otimes H_{13}$ is ordered as
follows:
$$
e_{00}=e_0\otimes e_0,\quad e_{01}=e_0\otimes
e_1,\quad e_{10}=e_1\otimes e_0,\quad e_{11}=e_1\otimes e_1.
$$
In this basis the operators $T_{12}$ and $T_{13}$ on the space $H_{12}\otimes H_{13}$ have the following form:
$$
T_{12}=
\left(\begin{smallmatrix}
0&0&1&0\\
0&0&0&1\\
1&0&0&0\\
0&1&0&0
\end{smallmatrix}\right)=
\left(\begin{smallmatrix}
0&I\\
I&0
\end{smallmatrix}\right)
,\,\,\,
T_{13}=
\left(\begin{smallmatrix}
0&1&0&0\\
1&0&0&0\\
0&0&0&1\\
0&0&1&0
\end{smallmatrix}\right)
=
\left(\begin{smallmatrix}
T&0\\
0&T
\end{smallmatrix}\right),
$$
where
$$
I=\left(\begin{smallmatrix}
1&0\\
0&1
\end{smallmatrix}\right),\quad
T=\left(\begin{smallmatrix}
0&1\\
1&0
\end{smallmatrix}\right).
$$
Consider the general case of $p$ and $\mu$. In the space $H=H_{\alpha_{12}}\otimes H_{\alpha_{13}}\otimes \dots\otimes
H_{\alpha_{1n}}$ the basis
$e_{i_2i_3\dots i_n},\,\,i_2,i_3,\dots,i_n\in {\mathbb F}_p$ is defined
by $ e_{i_2i_3\dots i_n}=e_{i_2}\otimes e_{i_3}\otimes \dots\otimes
e_{i_n}, $ and the scalar product for two elements $f$ and
$g$ in $H$
$$
f=\sum_{i_2,i_3,\dots,i_n\in {\mathbb
F}_p}f_{i_2i_3\dots i_n}e_{i_2i_3\dots i_n},\quad
g=\sum_{i_2,i_3,\dots,i_n\in {\mathbb
F}_p}g_{i_2i_3\dots i_n}e_{i_2i_3\dots i_n}
$$
is defined by the formula:
\begin{equation}
\label{(f,g)} (f,g)_{H_{12}\otimes H_{13}\otimes \dots\otimes
H_{1n}}=\sum_{i_2,i_3,\dots,i_n\in {\mathbb
F}_p}f_{i_2i_3\dots i_n}\overline{g}_{i_2i_3\dots i_n}.
\alpha_{12}(i_2)\alpha_{13}(i_3)\dots\alpha_{1n}(i_n).
\end{equation}
To describe $({\mathfrak A}^{1,2})'$ when $p=2$ we take any operator ${\rm
A}$ on the space $H_{12}\otimes H_{13}$
$$
{\rm A}=
\left(\begin{smallmatrix}
A&B\\
C&D
\end{smallmatrix}\right)
=
\left(\begin{smallmatrix}
a_{11}&a_{12}&b_{11}&b_{12}\\
a_{21}&a_{22}&b_{21}&b_{22}\\
c_{11}&c_{12}&d_{11}&d_{12}\\
c_{21}&c_{22}&d_{21}&d_{22}
\end{smallmatrix}\right).
$$
Since $T_{12}=\left(\begin{smallmatrix}
0&I\\
I&0
\end{smallmatrix}\right),\,\,\,
T_{13}=\left(\begin{smallmatrix}
T&0\\
0&T
\end{smallmatrix}\right)$ the relations $
[{\rm A},T_{13}]=0$ and $[{\rm A},T_{12}]=0$  gives us
$$
\left[
\left(\begin{smallmatrix}
A&B\\
C&D
\end{smallmatrix}\right),\left(\begin{smallmatrix}
T&0\\
0&T
\end{smallmatrix}\right)\right]=0\quad\text{and}\quad
\left[\left(\begin{smallmatrix}
A&B\\
C&D
\end{smallmatrix}\right),\left(\begin{smallmatrix}
0&I\\
I&0
\end{smallmatrix}\right)\right]=0,
$$
 or
 $$
\left(\begin{smallmatrix}
AT&BT\\
CT&DT
\end{smallmatrix}\right)=
\left(\begin{smallmatrix}
TA&TB\\
TC&TD
\end{smallmatrix}\right)\quad
{\rm and}\quad \left(\begin{smallmatrix}
B&A\\
D&C
\end{smallmatrix}\right)=
\left(\begin{smallmatrix}
C&D\\
A&B
\end{smallmatrix}\right).
$$
The second relation gives us $A=D$ and $B=C$, the first relation gives $TA=AT$ and $TB=BT$ hence, we get
 $$
\left(\begin{smallmatrix}
a_{21}&a_{22}\\
a_{11}&a_{12}
\end{smallmatrix}\right)=
\left(\begin{smallmatrix}
a_{12}&a_{11}\\
a_{22}&a_{21}
\end{smallmatrix}\right)
,\,\, {\rm and}\,\, \left(\begin{smallmatrix}
b_{21}&b_{22}\\
b_{11}&b_{12}
\end{smallmatrix}\right)=
\left(\begin{smallmatrix}
b_{12}&b_{11}\\
b_{22}&b_{21}
\end{smallmatrix}\right).
 $$
Finally, we get $ A=D,\,\,\,B=C,\,\,\,
a_{11}=a_{22},\,\,\,a_{12}=a_{21},\,\,\,b_{11}=b_{22},\,\,\,b_{12}=b_{21}
$ where $a_{11},a_{12},b_{11},b_{12}\in {\mathbb C}$ hence,
$$
{\rm A}=
\left(\begin{smallmatrix}
a_{11}&a_{12}&b_{11}&b_{12}\\
a_{12}&a_{11}&b_{12}&b_{11}\\
b_{11}&b_{12}&a_{11}&a_{12}\\
b_{12}&b_{11}&a_{12}&a_{11}
\end{smallmatrix}\right)
=a_{11} \left(\begin{smallmatrix}
I&0\\
0&I
\end{smallmatrix}\right)+
a_{12}\left(\begin{smallmatrix}
T&0\\
0&T
\end{smallmatrix}\right)+
b_{11} \left(\begin{smallmatrix}
0&I\\
I&0
\end{smallmatrix}\right)
+b_{12}\left(\begin{smallmatrix}
0&T\\
T&0
\end{smallmatrix}\right)
$$
$$
=a_{11}(I\otimes I)+a_{12}( I\otimes T)+b_{11}(T\otimes
I)+b_{12}(T\otimes T).
$$
Finally, in the case of $p=2$ the following statement is proved
\begin{lem}
\label{l.comA12} The von Neumann algebra $L^\infty(T_{12},T_{13})$
is maximal abelian, i.e.,
$(L^\infty(T_{12},T_{13}))'=L^\infty(T_{12},T_{13})$.
\end{lem}
In the case of an arbitrary ${\mathbb F}_p$ and the space
$X^{1,n}$ denote by $L^\infty(T_{12},\dots,T_{1n})$ the von Neumann algebra generated by operetors $T_{1k},\,2\leq k\leq n$.
\begin{lem}
\label{l.comA1n} The von Neumann algebra
$L^\infty(T_{12},\dots,T_{1n})$ is maximal abelian, i.e.,
$(L^\infty(T_{12},\!\dots,\!T_{1n}))'\!=\!L^\infty(T_{12},\!\dots,\!T_{1n}).$
In other words, any operator $A\!\in\!
(L^\infty(T_{12},\!\dots,\!T_{1n}))'$ has the following form:
\begin{equation}
\label{A-in-L^infty(n)}
A=
\sum_{i_2,i_3,\dots,i_n\in{\mathbb
F}_p}a_{i_2i_3\dots i_n}T_{12}^{i_2}T_{13}^{i_3}\dots T_{1n}^{i_n}.
\end{equation}
\end{lem}
\begin{pf} The proof follows from the fact that the von Neumann
algebra $L^\infty(T_{1k})$ is maximal abelian, i.e.,
$(L^\infty(T_{1k}))'=L^\infty(T_{1k})$,  by Theorem \ref{t.(T)'}.
Indeed, let us consider  the operator
$T=T_{inv}=\sum_{k\in{\mathbb F}_p}E_{k+1k}$ (see (\ref{T-alpha}))
acting on the space $L^2({\mathbb F}_p,\mu_{inv}),$ then
$T_{1k}={\underset{k}{\underbrace{I\otimes I\otimes \dots \otimes I\otimes
T}}}\otimes I\otimes \dots \otimes I$. Finally,
\begin{align*}
(L^\infty(T_{12},\dots,T_{1n}))'=&(T_{1k}\mid 2\leq k\leq n)'=\bigcap_{k=2}^n(T_{1k})'=\otimes_{k=2}^n(T_{1k})'\\
=&\otimes_{k=2}^n(L^\infty(T_{1k}))'=\otimes_{k=2}^nL^\infty(T_{1k})=L^\infty(T_{12},\dots,T_{1n}).
\end{align*}
 \qed\end{pf}
We calculate explicitly the commutant $({\mathfrak A}^{1,n})'$ for $p=2$
and small $n=3,4$  to guess  the general rule.
\begin{lem}
\label{l.1} In the case $p=2$ the commutant $({\mathfrak
A}^{1,3})'$ of the von Neumann algebra ${\mathfrak A}^{1,3}=
(T_{12},T_{13},T_{23})''$ is generated by operators
$T_{12}(I+T_{13})$ and $T_{13}$ or by $T_{12}C(T_{13})$ and
$T_{13}${\rm :}
$$
({\mathfrak A}^{1,3})'=\left(T_{12}C(T_{13}),\,T_{13}\right)''.
$$
\end{lem}
\begin{pf}
Let ${\rm A}\in ({\mathfrak A}^{1,3})',$ since $({\mathfrak
A}^{1,3})'=(T_{12},T_{13})'\bigcap(T_{23})'$ so, by Lemma
\ref{l.comA1n},
$$
({\mathfrak A}^{1,3})'=\left(A=aI+bT_{13}+cT_{12}+dT_{12}T_{13}\mid [{\rm A},T_{23}]=0\right).
$$
The operator $T_{23}$ has the following form (see (\ref{T1nT2n}))
$$
T_{23}=
\left(\begin{smallmatrix}
I&0\\
0&T
\end{smallmatrix}\right)
=
\left(\begin{smallmatrix}
1&0&0&0\\
0&1&0&0\\
0&0&0& 1\\
0&0&1&0
\end{smallmatrix}\right).
$$
We prove that for $k<r<s$ holds
\begin{equation}
\label{[kr,rs]} T_{kr}T_{rs}=T_{rs}T_{kr}T_{ks},\quad
T_{kr}^vT_{rs}=T_{rs}T_{kr}^vT_{ks}^v,\,v\in {\mathbb F}_p,
\end{equation}
in particular, for $(k,r,s)=(1,2,3)$ we have
$
T_{12}T_{23}=T_{23}T_{12}T_{13}.
$
Let us denote $E_{kr}(t)=I+tE_{kr},\,\,t\in{\mathbb F}_p,\,\,k<r$,
then we have by (\ref{T_{kn},T_{kn}(r)})
$T_{kn}^{-1}:=T^{R,\mu_\alpha,m}_{E_{kn}(1)}$. We calculate
$$
E_{12}(t)E_{23}(s)=\left(\begin{smallmatrix}
1&t&0\\
0&1&0\\
0&0& 1
\end{smallmatrix}\right)
\left(\begin{smallmatrix}
1&0&0\\
0&1&s\\
0&0& 1
\end{smallmatrix}\right)=
\left(\begin{smallmatrix}
1&t&ts\\
0&1&s\\
0&0& 1
\end{smallmatrix}\right)=
\left(\begin{smallmatrix}
1&t&0\\
0&1&s\\
0&0& 1
\end{smallmatrix}\right)
\left(\begin{smallmatrix}
1&0&ts\\
0&1&0\\
0&0& 1
\end{smallmatrix}\right)=
$$
$$
=\left(\begin{smallmatrix}
1&0&0\\
0&1&s\\
0&0& 1
\end{smallmatrix}\right)
\left(\begin{smallmatrix}
1&t&0\\
0&1&0\\
0&0& 1
\end{smallmatrix}\right)
\left(\begin{smallmatrix}
1&0&ts\\
0&1&0\\
0&0& 1
\end{smallmatrix}\right)=E_{23}(s)E_{12}(t)E_{13}(st),
$$
hence, $E_{12}(t)E_{23}(s)=E_{23}(s)E_{12}(t)E_{13}(st)$ or, if we
take $s=t=1$, we get (\ref{[kr,rs]}) for $k=1,\,r=2,\,s=3,\,v=1$.
The proof for an arbitrary $v\in {\mathbb F}_p$ is similar.

Using (\ref{[kr,rs]}) we get for
$A=a_{00}I+a_{01}T_{13}+a_{10}T_{12}+a_{11}T_{12}T_{13}$:
$$
A T_{23}=
a_{00}T_{23}+a_{01}T_{13}T_{23}+T_{23}T_{12}T_{13}(a_{10}+a_{11}T_{13})=
$$
$$
a_{00}T_{23}+a_{01}T_{13}T_{23}+T_{23}T_{12}(a_{10}T_{13}+a_{11}T_{13}^2)=
$$
$$
a_{00}T_{23}+a_{001}T_{13}T_{23}+T_{23}T_{12}(a_{10}T_{13}+a_{11}).
$$
Similarly, we have
$$
T_{23}A =a_{00}T_{23}+a_{01}T_{13}T_{23}+T_{23}T_{12}(a_{10}+a_{11}T_{13}).
$$
The condition $A T_{23}=T_{23}A$ gives us
$(a_{10}T_{13}+a_{11})=(a_{10}+a_{11}T_{13})$ or $a_{10}=a_{11}$.
Hence, $A=a_{00}I+a_{01}T_{13}+a_{10}T_{12}(I+T_{13})$ and lemma
is proved. \qed\end{pf}

\begin{lem}
\label{l.2} The commutant $({\mathfrak A}^{1,4})'$ of the von
Neumann algebra ${\mathfrak A}^{1,4}$ is generated by operators
$T_{12}(I+T_{13})(I+T_{14}),\,\,\,T_{13}(I+T_{14})$ and $T_{14}$
or by operators $T_{12}C(T_{13})C(T_{14}),\,\,T_{13}C(T_{14})$ and
$T_{14}$:
$$
({\mathfrak
A}^{1,4})'=\left(T_{12}C(T_{13})C(T_{14}),\,\,T_{13}C(T_{14}),\,\,T_{14}
\right)''.
$$
\end{lem}
\begin{pf}
Let ${\rm A}\in ({\mathfrak A}^{1,4})'$, since $ {\mathfrak
A}^{1,4}=(T_{12},T_{13},T_{14},T_{23},T_{24},T_{34})'' $ and
$T_{23}T_{34}=T_{34}T_{23}T_{24}$ or $\{T_{23},T_{34}\}=T_{24},$
where $\{a,b\}:=aba^{-1}b^{-1}$ we conclude that $({\mathfrak
A}^{1,4})'=(T_{12},T_{13},T_{14})'$ $\bigcap(T_{23},T_{34})' $ so,
$$({\mathfrak A}^{1,4})'=\left(
{\rm A}\in (T_{12},T_{13},T_{14})' \mid [{\rm A},T_{23}]=[{\rm
A},T_{34}]=0 \right).
$$
Using Lemma \ref{l.comA1n} we have for $ A\in
(T_{12},T_{13},T_{14})'$, $\quad A=$
$$
a_{000}\!+\!a_{100}T_{12}\!+\!a_{010}T_{13}\!+\!a_{001}T_{14}\!+\!a_{110}T_{12}T_{13}\!+\!
a_{101}T_{12}T_{14}\!+\!a_{011}T_{13}T_{14} \!+\!a_{111}T_{12}T_{13}T_{14}
$$
$$
=a_{000}+a_{010}T_{13}+a_{001}T_{14}+a_{011}T_{13}T_{14}
+T_{12}\left(a_{100}+a_{110}T_{13}+ a_{101}T_{14}
+a_{111}T_{13}T_{14}\right)
$$
$$
=a_{000}+a_{100}T_{12}+a_{001}T_{14}+a_{101}T_{12}T_{14}
+T_{13}\left(a_{010}+a_{110}T_{12}+a_{011}T_{14}
+a_{111}T_{12}T_{14}\right).
$$
The condition $[ A,T_{23}]=0$ gives us
$$  0=AT_{23}-T_{23}A=
T_{12}T_{23}\left(a_{100}+a_{110}T_{13}+ a_{101}T_{14}
+a_{111}T_{13}T_{14}\right)
$$
$$
-T_{23}T_{12}\left(a_{100}+a_{110}T_{13}+ a_{101}T_{14}
+a_{111}T_{13}T_{14}\right).
$$
Since $T_{12}T_{23}=T_{23}T_{12}T_{13}$, we have
$$
T_{12}T_{23}\left(a_{100}+a_{110}T_{13}+ a_{101}T_{14}
+a_{111}T_{13}T_{14}\right)$$
$$
=T_{23}T_{12}T_{13}\left(a_{100}+a_{110}T_{13}+ a_{101}T_{14}
+a_{111}T_{13}T_{14}\right)
$$
$$
=T_{23}T_{12}\left(a_{100}T_{13}+a_{110}+ a_{101}T_{13}T_{14}.
+a_{111}T_{14}\right).
$$
Therefore,
$$
a_{100}T_{13}+a_{110}+ a_{101}T_{13}T_{14} +a_{111}T_{14}=
a_{100}+a_{110}T_{13}+ a_{101}T_{14} +a_{111}T_{13}T_{14}
$$
hence,
\begin{equation}
\label{14.1}
 a_{100}=a_{110},\,\,a_{101}=a_{111}.
\end{equation}
Similarly, we get using condition $[ A,T_{34}]=0$,
$$  0=AT_{34}-T_{34}A=
T_{13}T_{34}\left(a_{010}+a_{110}T_{12}+a_{011}T_{14}
+a_{111}T_{12}T_{14}\right)-
$$
$$
T_{34}T_{13}\left(a_{010}+a_{110}T_{12}+a_{011}T_{14}
+a_{111}T_{12}T_{14}\right).
$$
Since $T_{23}T_{34}=T_{34}T_{23}T_{24}$, we have
$$
T_{13}T_{34}\left(a_{010}+a_{110}T_{12}+a_{011}T_{14}
+a_{111}T_{12}T_{14}\right)$$
$$
=T_{34}T_{13}T_{14}\left(a_{010}+a_{110}T_{12}+a_{011}T_{14}
+a_{111}T_{12}T_{14}\right)
$$
$$
=T_{34}T_{13}\left(a_{010}T_{14}+a_{110}T_{12}T_{14}+a_{011}
+a_{111}T_{12}\right)
$$
hence,
$$
a_{010}T_{14}+a_{110}T_{12}T_{14}+a_{011} +a_{111}T_{12}=
a_{010}+a_{110}T_{12}+a_{011}T_{14} +a_{111}T_{12}T_{14},
$$
and finally, we get
\begin{equation}
\label{14.2}
 a_{010}=a_{011},\,a_{110}=a_{111}.
\end{equation}
Using (\ref{14.1}) and (\ref{14.2}) we conclude
$a_{100}=a_{110}=a_{101}=a_{111},\,\,a_{010}=a_{011}$ hence,
$$
A=a_{000}I+a_{100}T_{12}(I+T_{13}+T_{14}+T_{13}T_{14})+a_{010}T_{13}(I+T_{14})+a_{001}T_{14}
$$
$$
=a_{000}I+a_{100}T_{12}(I+T_{13})(I+T_{14})+a_{010}T_{13}(I+T_{14})+a_{001}T_{14}.
$$
\qed\end{pf}
The previous lemmas were proved for ${\mathbb F}_2$ and von
Neumann algebras ${\mathfrak A}^{1,n}$ with $n=3,\,n=4.$ For the
general case ${\mathfrak A}^{1,n}$ and arbitrary ${\mathbb F}_p$
we can guess
\begin{lem}
\label{l.3} The commutant $({\mathfrak A}^{1,n})'$ of the von
Neumann algebra ${\mathfrak A}^{1,n}$ as the linear
space  is generated by  the following operators (we set $\Delta_{n,n}^{1,r}:=T_{1n}^r$){\rm:}
\begin{equation}
\label{centre1n}
({\mathfrak A}^{1,n})'=\Big(\Delta_{s,n}^{1,r}:=T_{1s}^r\prod_{k=s+1}^np^{-1}C(T_{1k})\mid 2\leq s\leq n,\,r\in
{\mathbb F}_p\setminus\{0\}\Big)''.
\end{equation}
The dimension of the von Neumann algebra $({\mathfrak A}^{1,n})'$ equals to $(n-1)(p-1)+1.$
\end{lem}
\begin{pf} We prove the statement first for $p=3$ and $n=3$. Any operator
$A\in L^\infty(T_{12},T_{13},T_{14})$ can be expressed as follows:
$$
A=\sum_{i_1,i_2,i_3\in{\mathbb
F}_3}a_{i_1,i_2,i_3}T_{12}^{i_1}T_{13}^{i_2}T_{14}^{i_3}.
$$
Rewrite the operator $A$ in the following form:
$$
A=\sum_{i_2,i_3\in{\mathbb F}_3}a_{0,i_2,i_3}T_{13}^{i_2}T_{14}^{i_3}+
T_{12}\sum_{i_2,i_3\in{\mathbb F}_3}a_{1,i_2,i_3}T_{13}^{i_2}T_{14}^{i_3}+
T_{12}^2\sum_{i_2,i_3\in{\mathbb F}_3}a_{2,i_2,i_3}T_{13}^{i_2}T_{14}^{i_3},
$$
$$
A=\sum_{i_1,i_3\in{\mathbb F}_3}a_{i_1,0,i_3}T_{12}^{i_2}T_{14}^{i_3}+
T_{13}\sum_{i_1,i_3\in{\mathbb F}_3}a_{i_1,1,i_3}T_{12}^{i_2}T_{14}^{i_3}+
T_{13}^2\sum_{i_1,i_3\in{\mathbb F}_3}a_{i_1,2,i_3}T_{12}^{i_2}T_{14}^{i_3}.
$$
 Using  (\ref{[kr,rs]}) we get
$$
T_{12}^rT_{23}=T_{23}T_{12}^rT_{13}^r,\quad T_{13}^rT_{34}=T_{34}T_{13}^rT_{14}^r.
$$
Using the relation $AT_{23}=T_{23}A$, we conclude that
$$
T_{13}\sum_{i_2,i_3\in{\mathbb F}_3}a_{1,i_2,i_3}T_{13}^{i_2}T_{14}^{i_3}=
\sum_{i_2,i_3\in{\mathbb F}_3}a_{1,i_2,i_3}T_{13}^{i_2}T_{14}^{i_3},
$$
$$
T_{13}^2\sum_{i_2,i_3\in{\mathbb
F}_3}a_{2,i_2,i_3}T_{13}^{i_2}T_{14}^{i_3}=
\sum_{i_2,i_3\in{\mathbb
F}_3}a_{2,i_2,i_3}T_{13}^{i_2}T_{14}^{i_3},
$$
therefore,  we get
\begin{equation}
\label{p3,n3,1}
a_{1,i_2,i_3}=a_{1,i_2+1,i_3},\,\,
a_{2,i_2,i_3}=a_{2,i_2+2,i_3},\,\,\forall i_2\in{\mathbb F}_3.
\end{equation}
Using the relation $AT_{34}=T_{34}A$, we conclude that
$$
T_{14}\sum_{i_1,i_3\in{\mathbb F}_3}a_{i_1,1,i_3}T_{12}^{i_1}T_{14}^{i_3}=
\sum_{i_1,i_3\in{\mathbb F}_3}a_{i_1,1,i_3}T_{12}^{i_1}T_{14}^{i_3},
$$
$$
T_{14}^2\sum_{i_1,i_3\in{\mathbb F}_3}a_{i_1,2,i_3}T_{12}^{i_1}T_{14}^{i_3}=
\sum_{i_1,i_3\in{\mathbb F}_3}a_{i_1,2,i_3}T_{12}^{i_1}T_{14}^{i_3},
$$
\begin{equation}
\label{p3,n3,2}
\text{so,}\quad a_{i_1,1,i_3}\!=\!a_{i_1,1,i_3+1},\,\,
a_{i_1,2,i_3}=a_{i_1,2,i_3+2},\,\,\forall i_3\in{\mathbb F}_3.
\end{equation}
Using (\ref{p3,n3,1}) and (\ref{p3,n3,2}) we conclude that
$$
a_{2i_2i_3}=a_{200},\,\,
a_{1i_2i_3}=a_{100},\,\,\,\forall i_2,i_3,\in{\mathbb F}_3,\,\,\,
a_{01i_3}=a_{010},\,\,a_{02i_3}=a_{020},\,\,\forall i_2\in{\mathbb F}_3.
$$
This implies that $A$ has the following form:
$$
A=a_{000}I+a_{001}T_{14}+a_{002}T_{14}^2+
a_{010}T_{13}C(T_{14})+a_{020}T_{13}^2C(T_{14})+
$$
$$
a_{100}T_{12}C(T_{13})C(T_{14})+a_{200}T_{12}^2C(T_{13})C(T_{14}).
$$
To prove the statement for general $p$ and $n$, set $T_k=T_{1k+1},\,\,k=1,2,\dots,n$. Let an operator $A$ has the following form:
$$
A=\sum_{i_1,i_2,\dots,i_n\in{\mathbb F}_p}a_{i_1,i_2,\dots,i_n}T_{1}^{i_1}T_{2}^{i_2}\dots T_{n}^{i_n}.
$$
Rewrite the operator $A$ in the following form for all $r$:
$$
A=\sum_{i_r\in{\mathbb F}_p}T_r^{i_r}\sum_{i_1,\dots,\hat{i}_r,\dots,i_n\in{\mathbb F}_p}a_{i_1,\dots,\hat{i}_r,\dots,i_n}
T_{1}^{i_1}\dots\hat{T}_r\dots T_{n}^{i_n},\quad 1\leq r\leq n
$$
where $\hat{T}_r,$ (resp. $\hat{i}_r$) means that factor $T_r$ (resp, index $\hat{i}_r$) is absent in the expression.
The commutation relations $AT_{kk+1}=T_{kk+1}A$ for $1\leq k\leq n$ imply, as before, the following relations:
$$
T_2^r\sum_{i_2,\dots,i_n\in{\mathbb F}_p}a_{r,i_2,\dots,i_n}
T_{2}^{i_2}\dots T_{n}^{i_n}=
\sum_{i_2,\dots,i_n\in{\mathbb F}_p}a_{r,i_2,\dots,i_n}
T_{2}^{i_2}\dots T_{n}^{i_n},
$$
$$
T_3^r\sum_{i_1,\hat{i}_2,i_3\dots,i_n\in{\mathbb F}_p}a_{i_1,r,i_3,\dots,i_n}
T_1^{i_1}\hat{T}_2T_{3}^{i_3}\dots T_{n}^{i_n}=
\sum_{i_1,\hat{i}_2,i_3\dots,i_n\in{\mathbb F}_p}a_{i_1,r,i_3,\dots,i_n}
T_1^{i_1}\hat{T}_2T_{3}^{i_3}\dots T_{n}^{i_n}.
$$
In the general case of $r,\,\,1\leq r\leq n$ we get
{\small
$$
T_r^s\!\!\!\!\!\!\!\sum_{i_1,\dots,\hat{i_r},\dots,i_n\in{\mathbb F}_p}\!\!\!\!\!\!a_{i_1,\dots,r,i_{r+1},\dots,i_n}
T_{1}^{i_1}\!\dots\!\hat{T}_r\!\dots\! T_{n}^{i_n}=\!\!\!\!\!\!\!\!
\sum_{i_1,\dots,\hat{i_r},\dots,i_n\in{\mathbb F}_p}a_{i_1,\dots,r,i_{r+1},\dots,i_n}
T_{1}^{i_1}\!\dots\!\hat{T}_r\!\dots\! T_{n}^{i_n}.
$$
}
The previous relations implies the analog of the relations (\ref{p3,n3,1}) and (\ref{p3,n3,2}):
\begin{equation}
\label{p,n,12}
a_{r,i_2,\dots,i_n}=a_{r,i_2+r,\dots,i_n}\,\,\,\forall r,\,i_2\in{\mathbb F}_p,\,\,\,
a_{i_1,r,i_3,\dots,i_n}=a_{i_1,r,i_3+r,\dots,i_n}\,\,\,\forall r,\,i_3\in{\mathbb F}_p,
\end{equation}
\begin{equation}
\label{p,n,r}
a_{i_1,\dots,r,i_{r+1},\dots,i_n}=a_{i_1,\dots,r,i_{r+1}+r,\dots,i_n}\,\,\,\forall r,\,i_{r+1}\in{\mathbb F}_p.
\end{equation}
Using (\ref{p,n,12}) and  (\ref{p,n,r}) we conclude that
$$
a_{r,i_2,\dots,i_n}=a_{r,0,\dots,0},\quad\forall i_2,\dots,i_n\in{\mathbb F}_p,\quad
a_{0,r,i_3,\dots,i_n}=a_{0,r,0,\dots,0},\,\,\forall i_3,\dots,i_n\in{\mathbb F}_p,
$$
$$
a_{0,\dots,0,r,i_{r+1},\dots,i_n}=a_{0,\dots,0,r,0,\dots,0}\,\,\,\forall i_{r+1},\dots,i_n\in{\mathbb F}_p.
$$
This implies that $A$ has the following form:
$$
A\!=\!\!\sum_{r\in{\mathbb F}_p}a_{0,\dots,0,r}T_n^r+\!\!\!\!
\sum_{r\in{\mathbb F}_p\setminus\{0\}}a_{0,\dots,0,r,0}T_{n-1}^rC(T_n)+\dots+
\!\!\!\!\sum_{r\in{\mathbb F}_p\setminus\{0\}}\!\!a_{r,0,\dots,0}T_1^rC(T_2)\dots C(T_n).
$$
\qed\end{pf}
\begin{rem}
\label{c(T)-C(T)} We have proved in the previous lemma that  the
von Neumann algebra $({\mathfrak A}^{1,n})'$  as the linear space
is generated by  the operators $\delta_{s,n}^{1,r}$:
$$
({\mathfrak A}^{1,n})'=\Big(\delta_{s,n}^{1,r}:=T_{1s}^r\prod_{k=s+1}^nC(T_{1k})\mid 2\leq s\leq n,\,r\in
{\mathbb F}_p\setminus\{0\}\Big)''.
$$
But the uniform limit $\lim_n\delta_{s,n}^{1,r}$ is divergent,
since
$$
\Vert\prod_{k=s+1}^nC(T_{1k})\Vert=\prod_{k=s+1}^n\Vert C(T_{1k})\Vert=p^{n-s}\to\infty,\quad
\text{as}\quad n\to\infty.
$$
We use the fact that $\Vert C(T)\Vert=p$.
{\rm Instead of the basis $\delta_{s,n}^{1,r}$, we choose the basis $\Delta_{s,n}^{1,r}$ of the algebra
$({\mathfrak A}^{1,n})'$ in Lemma~\ref{l.3} to be shur that the  limit
$\lim_n\Delta_{s,n}^{1,r}$ is correctly defined}.
Consider again the expression $C(T)=\sum_{k\in
{\mathbb F}_p}T^k$. Since $C(T)T=TC(T)=C(T)$ we get $C(T)^2=pC(T)$
so, $C(T)$ is ``almost projector'', i.e., $A^2=\lambda A$.
The operator $C(T)$ has two eigenvalues, $\lambda_1=0$ and $\lambda_2=p$. Indeed, if
 $C(T)f=\lambda f$, then $\lambda^2f=p\lambda f$ so, $\lambda(\lambda-p)=0$. Therefore,
$\Vert C(T)\Vert=\max\{0,p\}=p$.
\begin{equation}
\label{c(T)}
\text{Set}\quad c(T)=p^{-1}C(T).\quad\text{We have}\quad
c^2(T)=c(T)
\end{equation}
therefore, the eigenvalues of $c(T)$ are $0$ and $1$
hence, $\Vert c(T)\Vert=1$ so, the operator
$\Delta_{s,\infty}^{1,r}=\lim_n\Delta_{s,n}^{1,r}$, at least
formally, is correctly defined since
$$
\Vert \Delta_{s,\infty}^{1,r}\Vert=
\lim_n\Vert \Delta_{s,n}^{1,r}\Vert=\lim_n \Vert T_{1s}^r\prod_{k=s+1}^n c(T_{1k})\Vert=
\Vert T_{1s}^r\Vert\prod_{k=s+1}^n\Vert c(T_{1s})\Vert=1.
$$
In Lemma~\ref{l.Delta_1} below we prove that the operator $\lim_n \Delta_{s,n}^{1,r}$
is correctly defined when $\mu_\alpha\sim\mu_{inv}$.
\end{rem}
\begin{rem}
\label{Alg=()''}
The von Neuman algebra $({\mathfrak A}^{1,n})'$ as  an algebra is generated
by the following expressions:
$$
 ({\mathfrak A}^{1,n})'=\Big(\Delta_{s,n}^{1,1}\mid 2\leq s\leq n\Big)''.
$$
\end{rem}
\begin{pf} It is sufficient to use Lemma~\ref{l.3} and the following relations:
$$
\Delta_{s,n}^{1,r}\Delta_{s,n}^{1,t}=\Delta_{s,n}^{1,r+t},\quad
\Delta_{s_1,n}^{1,r_1}\Delta_{s_2,n}^{1,r_2}=\Delta_{s_1,n}^{1,r_1}\quad\text{for}\quad
3\leq s_1<s_2\leq n.
$$
Using the relation $c^2(T)=c(T)$ and $Tc(T)=c(T)$ we get $\Delta^{1,r}_{sn}\Delta^{1,l}_{sn}=
\Delta^{1,r+l}_{sn}$.
Indeed,
$$
\Delta^{1,r}_{sn}\Delta^{1,l}_{sn}=
T^r_{1s}\prod_{k=s+1}^n c(T_{1k})T^l_{1s}\prod_{t=s+1}^n c(T_{1t})=T^{r+l}_{1s}\prod_{k=s+1}^n c^2(T_{1k})=
\Delta^{1,r+l}_{sn}.
$$
Similarly, we prove the second relation $\Delta_{s_1,n}^{1,r_1}\Delta_{s_2,n}^{1,r_2}=\Delta_{s_1,n}^{1,r_1}$.
\qed\end{pf}

{\it Another description of  the commutant $({\mathfrak
A}^{1,n})'$}. Any operator\\ $A\in L^\infty(T_{12},\dots,T_{1n})$
has the following form by (\ref{A-in-L^infty(n)}):
$A=f(T_{12},\dots,T_{1n})$.
\begin{lem}
\label{l.[f,T(k,k+1)]=0(n)}
An operator $A=f(T_{12},\dots,T_{1n})\in L^\infty(T_{12},\dots,T_{1n})$  commutes with  $T_{kk+1}$ for all
$2\leq k\leq n-1$ if and only if for all $2\leq k\leq n-1$ holds:
\begin{equation}
\label{[f,T(k,k+1)]=0(n)}
f(T_{12},\dots,T_{1k},T_{1k+1},\dots,T_{1n})=f(T_{12},\dots,T_{1k}T_{1k+1},T_{1k+1},\dots,T_{1n}).
\end{equation}
\end{lem}
\begin{pf}
Consider first the space $H_{1k}\otimes H_{1k+1}$ and the von Neumann subalgebra
$L^\infty(T_{1k},T_{1k+1})$ in the algebra $B(H_{1k}\otimes H_{1k+1})=B(H_{1k})\otimes B( H_{1k+1})$.
Take the function $f(T_{1k},T_{1k+1})\in L^\infty(T_{1k},T_{1k+1})$ of the following form
$f(T_{1k},T_{1k+1})$\\$=T_{1k}^rT_{1k+1}^s$. We show that commutation relation $[f,T_{kk+1}]=0$ implies
the relation (\ref{[f,T(k,k+1)]=0(n)}), i.e., $f(T_{1k},T_{1k+1})=f(T_{1k}T_{1k+1},T_{1k+1})$. Indeed,
using  the relations (\ref{[kr,rs]}) $T_{1k}^rT_{km}=T_{km}T_{1k}^rT_{1m}^r$ for $r\in {\mathbb F}_p$
and $2\leq k<m$ and commutation relation
$f(T_{1k},T_{1k+1})T_{kk+1}=T_{kk+1}f(T_{1k},T_{1k+1})$
 we get:
\begin{gather*}
T_{kk+1}f(T_{1k},T_{1k+1})=f(T_{1k},T_{1k+1})T_{kk+1}=T_{1k}^rT_{1k+1}^sT_{kk+1}=\\
T_{kk+1}T_{1k}^rT_{1k+1}^rT_{1k+1}^s=T_{kk+1}f(T_{1k}T_{1k+1},T_{1k+1}).
\end{gather*}
Finally, we prove $f(T_{1k},T_{1k+1})=f(T_{1k}T_{1k+1},T_{1k+1})$ for the particalar case of\\
$f(T_{1k},T_{1k+1})\!=\!T_{1k}^rT_{1k+1}^s$.
For the general function $f(T_{1k},T_{1k+1})\!=\!\sum_{r,s\in {\mathbb F}_p}a_{r,s}\times $\\$T_{1k}^rT_{1k+1}^s$
the proof is the same. Similarly, we prove (\ref{[f,T(k,k+1)]=0(n)}) for any function
$\in L^\infty(T_{12},\dots,T_{1n})$:
$$
f(T_{12},\dots,T_{1k})=\sum_{i_2,\dots,i_n\in {\mathbb F}_p}a_{i_2,\dots,i_n}T_{12}^{i_2}\dots T_{1n}^{i_n}.
$$
\qed\end{pf}
\begin{lem}
\label{l.f=f(c(T))} If for the function $f
\in L^\infty(T_{12},\dots,T_{1n})$ holds relation (\ref{[f,T(k,k+1)]=0(n)}) then
\begin{equation}
\label{f=f(c(T))}
f=f_s(T_{1s})\prod_{k=s+1}^nc(T_{1k})\quad\text{for some}\quad
s,\,\,2\leq s\leq n,
\end{equation}
where $f_s(T_{1s})\in L^\infty(T_{1s})$.
\end{lem}
\begin{pf}
If for the function $f(T_{1k},T_{1k+1})=\sum_{r,s\in {\mathbb F}_p}a_{r,s}T_{1k}^rT_{1k+1}^s$ holds\\
$f(T_{1k},T_{1k+1})=f(T_{1k}T_{1k+1},T_{1k+1})$ then we get
\begin{align*}
f(T_{1k}T_{1k+1},T_{1k+1})=\sum_{r,s\in {\mathbb F}_p}a_{r,s}(T_{1k}T_{1k+1})^rT_{1k+1}^s=
\sum_{r,s\in {\mathbb F}_p}a_{r,s}T_{1k}^rT_{1k+1}^{r+s}=\\
\sum_{r,t\in {\mathbb F}_p}a_{r,t-r}T_{1k}^rT_{1k+1}^t=\sum_{r,t\in {\mathbb F}_p}a_{r,t}T_{1k}^rT_{1k+1}^t.
\end{align*}
Therefore, $a_{r,t-r}=a_{r,t}$ for all $r,t\in {\mathbb F}_p$ hence, $a_{r,t}=a_{r,t-kr}$ for all $r,t,k\in {\mathbb F}_p$.
Since ${\mathbb F}_p$ is a field, we conclude that
\begin{equation}
\label{f-inv1}
a_{r,t}=a_{r,0}\quad\text{for all}\quad r,t\in {\mathbb F}_p.
\end{equation}

Finally, if we set $f_k(T_{1k})=\sum_{r\in {\mathbb F}_p}a_{r,0}T_{1k}^r$ we get $f=$
\begin{align*}
\sum_{r,s\in {\mathbb F}_p}\!\!a_{r,s}T_{1k}^rT_{1k+1}^s\!\!=\!\!\sum_{r,s\in {\mathbb F}_p}\!\!a_{r,0}T_{1k}^rT_{1k+1}^s
\!\!=\!\!\sum_{r\in {\mathbb F}_p}\!\!a_{r,0}T_{1k}^r\sum_{s\in {\mathbb F}_p}T_{1k+1}^s\!\!=\!\!p^{-1}f_k(T_{1k})c(T_{1k+1}).
\end{align*}
If for the function $f(T_{1k},T_{1k+1},T_{1k+2})=\sum_{r,s,t\in
{\mathbb F}_p}a_{r,s,t}T_{1k}^rT_{1k+1}^sT_{1k+2}^t$ holds
$$
f(T_{1k},T_{1k+1},T_{1k+2})=f(T_{1k}T_{1k+1},T_{1k+1},T_{1k+2})=f(T_{1k},T_{1k+1}T_{1k+2},T_{1k+2})
$$
we conclude similarly, that $a_{r,s-r,t}=a_{r,s,t}=a_{r,s,t-s}$
for all $r,s,t\in {\mathbb F}_p$ hence,
\begin{equation}
\label{f-inv2}
a_{r,s,t}=a_{r,0,0}\quad\text{for all}\quad r,s,t\in {\mathbb F}_p.
\end{equation}
Finally, we get
\begin{align*}
f=\sum_{r,s,t\in {\mathbb F}_p}a_{r,s,t}T_{1k}^rT_{1k+1}^sT_{1k+2}^t=\sum_{r,s,t\in {\mathbb F}_p}a_{r,0,0}
T_{1k}^rT_{1k+1}^sT_{1k+2}^t=\sum_{r\in {\mathbb F}_p}a_{r,0,0}T_{1k}^r\times\\
\sum_{s\in {\mathbb F}_p}T_{1k+1}^s
\sum_{t\in {\mathbb F}_p}T_{1k+2}^t=f_k(T_{1k})C(T_{1k+1})C(T_{1k+2})=p^{-2}f_k(T_{1k})c(T_{1k+1})c(T_{1k+2}),
\end{align*}
where $f_k(T_{1k})=\sum_{r\in {\mathbb F}_p}a_{r,0,0}T_{1k}^r$.
\qed\end{pf}
\begin{lem}
\label{l.Delta_1}  When $\mu_\alpha\sim\mu_{inv}$ the commutant
$({\mathfrak A}^{1})'$ of the von Neumann algebra ${\mathfrak
A}^{1}$ is generated as a linear space by the following
expressions:
 \begin{equation}
\label{centre1(inf)}
({\mathfrak A}^{1})'=\Big(
\Delta_{s,\infty}^{1,r}=
\lim_{n}\Delta_{s,n}^{1,r}:=T_{1s}^r\prod_{k=s+1}^\infty p^{-1}C(T_{1k})\mid 2\leq s,\,r\in
{\mathbb F}_p\setminus\{0\}\Big)''.
 \end{equation}
When $\mu_\alpha\perp\mu_{inv}\Leftrightarrow S^L_{11}(\mu)=\infty$ the commutant $({\mathfrak A}^{1})'$ is trivial, i.e., $({\mathfrak A}^{1})'=(\lambda I\mid
 \lambda\in {\mathbb C})$.
\end{lem}
\begin{pf}

Denote by $L^\infty(T_1)=L^\infty(T_{1k}\mid 2\leq k)=(T_{1k}\mid
2\leq k)''$ the von Neumann algebra generated by the commuting
family of operators $T_1:=(T_{1k}\mid 2\leq k)$. By definition, we
have
\begin{equation}
({\mathfrak A}^{1})'=\Big(f\in L^\infty(T_1)\mid
[f,T_{kk+1}]=0,\quad 2\leq k\Big).
\end{equation}
Using the spectral theorem for the family $T_1=(T_{1k}\mid 2\leq
k)$ of commuting unitary operators $T_{1k}$ we conclude that any
element $f\in L^\infty(T_1)$ has the following form:
$f(T_1)=\int_{X^1}f(\lambda)dE(\lambda)$ or
\begin{equation}
\label{f=spectal-int}
f(T_1)=f(T_{1k}\mid 2\leq k)=\int_{X^1}f(\lambda)dE(\lambda)
\end{equation}
where $X^1=\prod_{k=2}^\infty({\mathbb F}_p)_k$, $f$ is essentially
bounded function on $X^1$ and $E$ is a common resolution of the
identity of the family of operators $T_1$ defined on cylindrical sets
$\Delta_2\times\cdots\times\Delta_k$ as follows:
$$
E(\Delta_2\times\cdots\times\Delta_k):=E_2(\Delta_2)\cdots E_k(\Delta_k)
$$
where $E_k$ is  resolution of the identity of the operators $T_{1k}$.
See details in \cite{Ber86}. Similarly to the proof of
Lemma~\ref{l.[f,T(k,k+1)]=0(n)}, we get
\begin{lem}
\label{l.[f,T(k,k+1)]=0(inf)}
An operator $A=f(T_{12},\dots,T_{1n},\dots)\in L^\infty(T_{1k}\mid 2\leq k)$ defined by (\ref{f=spectal-int})
 commutes  with  $T_{kk+1}$ for all
$2\leq k$ if and only if for all $2\leq k$ holds:
\begin{equation}
\label{[f,T(k,k+1)]=0(inf)}
f(T_{12},\dots,T_{1k},T_{1k+1},\dots,T_{1n},\dots)=f(T_{12},\dots,T_{1k}T_{1k+1},T_{1k+1},\dots,T_{1n},\dots).
\end{equation}
\end{lem}
\begin{lem}
\label{l.f=f(c(T))(inf)} If for the function $f=f(T_1)\in L^\infty(T_{12},\dots,T_{1n},\dots)$ holds
relation (\ref{[f,T(k,k+1)]=0(inf)}) then
\begin{equation}
\label{f=f(c(T))(inf)} f=f_s(T_{1s})\prod_{k=s+1}^\infty
c(T_{1k})\quad\text{for some}\quad  s\geq 2,
\end{equation}
where $f_s(T_{1s})\in L^\infty(T_{1s})$.
\end{lem}
\qed\end{pf}
\subsection{The commutant of the von Neumann algebra ${\mathfrak A}^{m}$, case $m>1$}
Let us consider  the restriction $T^{R,m,n}$ of the representation
$T^{R,\mu_\alpha,m}:B_0^{\mathbb N}({\mathbb F}_p)\mapsto
U(L^2(X^m,\mu_\alpha))$ to the subgroup $B(n,{\mathbb
F}_p),\,m\leq n$, of the group $B_0^{\mathbb N}({\mathbb F}_p)$
acting in the space $H^{m,n}=\otimes_{1\leq r\leq m,r<k\leq n}
L^2({\mathbb
F}_p,\mu_{\alpha_{rk}})=L^2(X^{m,n},\mu_{\alpha^{m,n}})$ where
$\mu_{\alpha^{m,n}}=\otimes_{1\leq r\leq m,r<k\leq
n}\mu_{\alpha_{rk}}$ and
\begin{equation}
\label{X^(m,n)}
X^{m,n}=
\left(\begin{smallmatrix}
1&x_{12}&x_{13}&\dots &x_{1m}&\dots &x_{1n}\\
0& 1    &x_{23}&\dots &x_{2m}&\dots &x_{2n}\\
0& 0    & 1    &\dots &x_{3m}&\dots &x_{3n}\\
 &      &      &\dots &      &\dots &      \\
0& 0    &      &\dots & x_{m-1m}    &\dots &x_{m-1n}\\
0& 0    &      &\dots &1     &\dots &x_{mn}
\end{smallmatrix}\right).
\end{equation}
Let us denote as before by ${\mathfrak A}^{m,n}$ and ${\mathfrak
A}^m$ the von-Neumann algebras generated by the representation
$T^{R,m,n}$ (respectively by $T^{R,\mu_\alpha,m}$)
$$
{\mathfrak A}^{m,n}\!=\Big( T^{R,m,n}_t\mid t\in B(n,{\mathbb F}_p)
 \Big)'',\,\,
{\mathfrak A}^m\!=\!\Big( T^{R,\mu_\alpha,m}_t\mid t\in B_0^{\mathbb
N}({\mathbb F}_p) \Big)''\!=\!\Big(\bigcup_{n\geq m}{\mathfrak
A}^{m,n}\Big)''.
$$
Obviously,  the commutant $\left({\mathfrak
A}^{m,n}\right)'$ contains the following operators:
\begin{equation*}
({\mathfrak A}^{m,n})'\supset \Big(\Delta_{s,n}^{k,s}:=T_s^r\prod_{k=s+1}^np^{-1}C(T_{kr}(k))\mid 1\leq k\leq m,\,\,k+1\leq s\leq n\Big)''.
\end{equation*}
Since ${\mathfrak A}^{m}=(\bigcup_{n\geq m+1}{\mathfrak
A}^{m,n})''$ we get $({\mathfrak A}^{m})'=\bigcap_{n\geq
m+1}({\mathfrak A}^{m,n})'$. Hence, the commutant $({\mathfrak
A}^{m})'$ is not trivial if there exist a non trivial limit
$\Delta_{k},\,\,1\leq k\leq m$
$$
\Delta_{k}:=\lim_{n\to\infty}\prod_{r=k+1}^np^{-1}C(T_{kr}(k))=\lim_{n\to\infty}\Delta_{k,n}^{k,0}.
$$
The latter limit exists,  if
$\mu_\alpha^k\sim\mu^k_{inv}\Leftrightarrow
S^L_{kk}(\mu_\alpha)<\infty$ (see  Lemma~\ref{Lap-exist}).

\subsection{Commutant of the von Neumann algebra $({\mathfrak A}^{m})'$, case $m=2$}
Set
$$
\Delta_{s,n}^{2,r}:=T_{\alpha_{2s}}^r\prod_{k=s+1}^np^{-1}C(T_{\alpha_{2k}}),\quad 3\leq s< n,\quad
\Delta_{n,n}^{2,r}:=T_{\alpha_{2n}}^r,\quad r\in {\mathbb F}_p.
$$
\begin{lem}
\label{l.x(12),(A^(2,n)'=}
Let $x_{12}\in {\mathfrak A}^{2,n}$, then the commutant $({\mathfrak A}^{2,n})'$ of the von
Neumann algebra ${\mathfrak A}^{2,n}$ as the linear
space  is generated by  the following operators:
\begin{equation}
\label{comm.m=2,n}
({\mathfrak A}^{2,n})'=\Big(\Delta_{s,n}^{l,r}:=T_{1s}^r\prod_{r=s+1}^np^{-1}C(T_{1r})\mid 1\leq l\leq2,\,\, 3\leq s\leq n,\,r\in {\mathbb F}_p\setminus\{0\}\Big)''.
\end{equation}
The dimension of the von Neumann algebra $({\mathfrak A}^{2,n})'$ equals to $[(n\!-\!1)(p\!-\!1)\!+\!1]^2.$
\end{lem}
\begin{pf}
Since $x_{12},\,\,T_{1k},\,\,T_{2k}\in {\mathfrak A}^{2,n}$ we conclude that
$T_{2k}(2)=T_{\alpha_{2k}}\in {\mathfrak A}^{2,n}$ (see  Remark~\ref{r.x(12),(A)'=}).
Since  the commutative family of operators with common simple spectrum belongs to   ${\mathfrak A}^{2,n}$, i.e.,
$x_{12},\,\,T_{1k},\,\,T_{\alpha_{2k}}\in {\mathfrak A}^{2,n}$ for $3\leq k\leq n$ we conclude that
\begin{equation}
\label{x(12),(A^(2,n)'=}
({\mathfrak A}^{2,n})'\subset L^\infty\left(\begin{smallmatrix}
x_{12}&T_{\alpha_{13}}&\cdots&T_{\alpha_{1n}}\\
          &T_{\alpha_{23}}&\cdots&T_{\alpha_{2n}}\\
\end{smallmatrix}\right).
\end{equation}
Commutation relation $[f,T_{12}]=0$ for $f\in L^\infty\left(\begin{smallmatrix}
x_{12}&T_{\alpha_{13}}&\cdots&T_{\alpha_{1n}}\\
          &T_{\alpha_{23}}&\cdots&T_{\alpha_{2n}}\\
\end{smallmatrix}\right)$
means that $f$ does not depend on $x_{12}.$
Fix $n=4$ and $p=2$,  let $f\left(\begin{smallmatrix}
T_{\alpha_{13}}&T_{\alpha_{14}}\\
 T_{\alpha_{23}}&T_{\alpha_{24}}\\
\end{smallmatrix}\right) \in L^\infty\left(\begin{smallmatrix}
T_{\alpha_{13}}&T_{\alpha_{14}}\\
 T_{\alpha_{23}}&T_{\alpha_{24}}\\
\end{smallmatrix}\right)$.
We use the following relations (see (\ref{T(kn)=otimes}))
$$
T_{34}\!=\!T_{34}(1)\otimes T_{34}(1),\,\,
T_{\alpha_{13}}T_{34}(1)\!=\!T_{34}(1)T_{\alpha_{13}}T_{\alpha_{14}},\,\,
T_{\alpha_{23}}T_{34}(2)\!=\!T_{34}(2)T_{\alpha_{23}}T_{\alpha_{24}}
$$
\begin{equation}
\label{[T,T]=(T-I)}
\text{or}\quad [T_{\alpha_{13}},T_{34}(1)]\!=\!T_{34}(1)T_{\alpha_{13}}(T_{\alpha_{14}}-I),\!\quad\!
[T_{\alpha_{23}},T_{34}(2)]\!=\!T_{34}(2)T_{\alpha_{23}}(T_{\alpha_{24}}-I).
\end{equation}
Any operator $A\in L^\infty\left(\begin{smallmatrix}
T_{\alpha_{13}}&T_{\alpha_{14}}\\
 T_{\alpha_{23}}&T_{\alpha_{24}}\\
\end{smallmatrix}\right)$ has the form
$$
A=\sum_{i_1,i_2,j_1,j_2\in {\mathbb F}_2}a^{i_1,i_2}_{j_1,j_2}T^{i_1}_{\alpha_{13}}
T^{i_2}_{\alpha_{14}}T^{j_1}_{\alpha_{23}}T^{j_2}_{\alpha_{24}}.
$$
We rewrite $A$ as follows:
\begin{align*}
A&=[a^{00}_{00}I+a^{01}_{00}T_{\alpha_{14}}+a^{00}_{01}T_{\alpha_{24}}
+a^{01}_{01}T_{\alpha_{14}}T_{\alpha_{24}}]\\
&+T_{\alpha_{23}}[(a^{00}_{10}+a^{00}_{11}T_{\alpha_{24}})
+T_{\alpha_{14}}(a^{01}_{10}+a^{01}_{11}T_{\alpha_{24}})]\\
&+T_{\alpha_{13}}[(a^{10}_{00}+a^{11}_{00}T_{\alpha_{14}})
+T_{\alpha_{24}}(a^{10}_{01}+a^{11}_{01}T_{\alpha_{14}})]\\
&+T_{\alpha_{13}}T_{\alpha_{23}}[(a^{10}_{10}+a^{10}_{11}T_{\alpha_{24}}+a^{11}_{10}T_{\alpha_{14}}+
a^{11}_{11}T_{\alpha_{14}}T_{\alpha_{24}})]\\
&=A_1+A_2+A_3+A_4\\
&=A_1+T_{\alpha_{23}}a_2+T_{\alpha_{13}}a_3+T_{\alpha_{13}}T_{\alpha_{23}}a_4
\end{align*}
where
\begin{align*}
A_1&=[a^{00}_{00}I+a^{01}_{00}T_{\alpha_{14}}+a^{00}_{01}T_{\alpha_{24}}
+a^{01}_{01}T_{\alpha_{14}}T_{\alpha_{24}}],\\
a_2&=[(a^{00}_{10}+a^{00}_{11}T_{\alpha_{24}})
+T_{\alpha_{14}}(a^{01}_{10}+a^{01}_{11}T_{\alpha_{24}})],\\
a_3&=[(a^{10}_{00}+a^{11}_{00}T_{\alpha_{14}})
+T_{\alpha_{24}}(a^{10}_{01}+a^{11}_{01}T_{\alpha_{14}})],\\
a_4&=[(a^{10}_{10}+a^{10}_{11}T_{\alpha_{24}}+a^{11}_{10}T_{\alpha_{14}}+
a^{11}_{11}T_{\alpha_{14}}T_{\alpha_{24}})].\\
\end{align*}
Since
$T_{34}(1)-I=x_{13}\otimes(T_{\alpha_{14}}-I),\,\,
T_{34}(2)-I=x_{23}\otimes(T_{\alpha_{24}}-I)$
(see (\ref{T2n(1)-1.0}) and (\ref{63})) we get
$$
[A_1,T_{34}(1)]=[A_2,T_{34}(1)]=[T_{34}(2),A_1]=[T_{34}(2),A_3]=0
$$
hence,
\begin{align*}
&[A,T_{34}(1)\otimes T_{34}(2)]=[A,T_{34}(1)]\otimes T_{34}(2)+T_{34}(1)\otimes [A,T_{34}(2)]\\
=&\big([A_3,T_{34}(1)]+[A_4,T_{34}(1)]\big)\otimes T_{34}(2)\\
+&T_{34}(1)\otimes \big([A_2,T_{34}(2)]+[A_4,T_{34}(2)]\big)\\
=&\big([T_{\alpha_{13}},T_{34}(1)]a_3+[T_{\alpha_{13}},T_{34}(1)]T_{\alpha_{23}}a_4\big)\otimes T_{34}(2)\\
+&T_{34}(1)\otimes\big(
[T_{\alpha_{23}},T_{34}(2)]a_2+T_{\alpha_{13}}[T_{\alpha_{23}},T_{34}(2)]a_4\big)\\
=&\big(T_{34}(1)T_{\alpha_{13}}(T_{\alpha_{14}}-I)a_3+T_{34}(1)T_{\alpha_{13}}T_{\alpha_{23}}
(T_{\alpha_{14}}-I)a_4\big)\otimes T_{34}(2)\\
+&T_{34}(1)\otimes\big(
T_{34}(2)T_{\alpha_{23}}(T_{\alpha_{24}}-I)a_2+T_{\alpha_{13}}T_{34}(1)T_{\alpha_{23}}
(T_{\alpha_{24}}-I)a_4
\big)=0.
\end{align*}
Therefore, we have
\begin{align*}
(T_{\alpha_{14}}-I)a_3=
(T_{\alpha_{14}}-I)a_4=
(T_{\alpha_{24}}-I)a_2=
(T_{\alpha_{24}}-I)a_4=0
\end{align*}
or
\begin{align*}
(T_{\alpha_{14}}-I)[(a^{10}_{00}+a^{11}_{00}T_{\alpha_{14}})
+T_{\alpha_{24}}(a^{10}_{01}+a^{11}_{01}T_{\alpha_{14}})]=0,\\
(T_{\alpha_{14}}-I)[(a^{10}_{10}+a^{10}_{11}T_{\alpha_{24}}+a^{11}_{10}T_{\alpha_{14}}+
a^{11}_{11}T_{\alpha_{14}}T_{\alpha_{24}})]=0,\\
(T_{\alpha_{24}}-I)[(a^{00}_{10}+a^{00}_{11}T_{\alpha_{24}})
+T_{\alpha_{14}}(a^{01}_{10}+a^{01}_{11}T_{\alpha_{24}})]=0,\\
(T_{\alpha_{24}}-I)[(a^{10}_{10}+a^{10}_{11}T_{\alpha_{24}}+a^{11}_{10}T_{\alpha_{14}}+
a^{11}_{11}T_{\alpha_{14}}T_{\alpha_{24}})]=0.
\end{align*}
Hence, we get respectively
\begin{align*}
\label{m=2,(34)}
a^{10}_{00}=a^{11}_{00},\quad a^{10}_{01}=a^{11}_{01},\quad
a^{10}_{10}=a^{11}_{10},\quad a^{10}_{11}=a^{11}_{11},\\\nonumber
a^{00}_{10}=a^{00}_{11},\quad a^{01}_{10}=a^{01}_{11},\quad
a^{10}_{10}=a^{10}_{11},\quad a^{11}_{10}=a^{11}_{11}.
\end{align*}
At last, using the latter equalities we get
\begin{align*}
A=&[a^{00}_{00}I+a^{01}_{00}T_{\alpha_{14}}+a^{00}_{01}T_{\alpha_{24}}
+a^{01}_{01}T_{\alpha_{14}}T_{\alpha_{24}}]\\
+&T_{\alpha_{23}}[a^{00}_{10}(I+T_{\alpha_{24}})+a^{01}_{10}T_{\alpha_{14}}(I+T_{\alpha_{24}})]\\
+&T_{\alpha_{13}}[a^{10}_{00}(I+T_{\alpha_{14}})+a^{10}_{01}T_{\alpha_{24}}(I+T_{\alpha_{14}})]\\
+&T_{\alpha_{13}}T_{\alpha_{23}}a^{10}_{10}[I+T_{\alpha_{24}}+T_{\alpha_{14}}+T_{\alpha_{14}}T_{\alpha_{24}}].
\end{align*}

Finally, the basis in the algebra $({\mathfrak A}^{2,4})'$ can be
obtained as a tensor product of two basis in two von Neumann algebras
$$
({\mathfrak A}^{1,3,4})'=L^\infty(T_{13},T_{14})\bigcap
(T_{34})'\quad\text{ and} \quad({\mathfrak
A}^{2,3,4})'=L^\infty(T_{\alpha_{23}},T_{\alpha_{24}})\bigcap
(T_{34})'
$$
due to the following relations:
{\small
$$
({\mathfrak A}^{2,4})'=({\mathfrak A}^{1,3,4})'\otimes
({\mathfrak A}^{2,3,4})'=
\Big(T_{\alpha_{13}}(I+T_{\alpha_{14}}),I, T_{\alpha_{14}})\Big)''\otimes
\Big(T_{\alpha_{23}}(I+T_{\alpha_{24}}),I, T_{\alpha_{24}})\Big)''.
$$
}
\qed\end{pf}
Set for $s\geq 3$
$$
\Delta^{1,r}_{s,\infty}:\!=\!T_{\alpha_{1s}}^r\!\prod_{k=s+1}^\infty
p^{-1}C(T_{\alpha_{1k}}),\quad
\Delta^{2,r}_{s,\infty}:\!=\!T_{\alpha_{2s}}^r\!\prod_{k=s+1}^\infty
p^{-1}C(T_{\alpha_{2k}}),\,\,r\in{\mathbb F}_p\setminus\{0\}.
$$
\begin{lem}
\label{l.x(12),(A^(2)'=}
Let $x_{12}\in {\mathfrak A}^{2}$, and $\mu^2_\alpha\sim\mu^2_{inv}$ or $S^L_{11}(\mu)<\infty,\,\,\,S^L_{22}(\mu)<\infty$, then the commutant $({\mathfrak A}^{2})'$ of the von
Neumann algebra ${\mathfrak A}^{2}$ as the linear
space  is generated by  the following operators:
\begin{equation}
\label{comm.m=2}
({\mathfrak A}^{2})'=\Big(\Delta_{s,\infty}^{l,r}
\mid 1\leq l\leq2,\,\, 3\leq s,\,r\in {\mathbb F}_p\setminus\{0\}\Big)''.
\end{equation}
When $\mu^2_\alpha\perp\mu^2_{inv}$ or
$S^L_{11}(\mu)=S^L_{22}(\mu)=\infty$ and  $x_{12}\in {\mathfrak
A}^{2}$, the commutant $({\mathfrak A}^{2})'$ is trivial.
\end{lem}
\begin{pf}
Since $x_{12},\,\,T_{1k},\,\,T_{2k}\in {\mathfrak A}^{2}$ for $k\geq 3$ we conclude that
$T_{2k}(2)=T_{\alpha_{2k}}\in {\mathfrak A}^{2}$ (see  Remark~\ref{r.x(12),(A)'=}).
Since  the commutative family of operators with common simple spectrum belongs to    ${\mathfrak A}^{2}$, i.e.,
$x_{12},\,\,T_{1k},\,\,T_{\alpha_{2k}}\in {\mathfrak A}^{2}$ for $3\leq k$ we conclude that
\begin{equation}
\label{x(12),(A^2)'=}
({\mathfrak A}^{2})'\subset L^\infty\left(\begin{smallmatrix}
T_{\alpha_{13}}&\cdots&T_{\alpha_{1n}}&\dots\\
 T_{\alpha_{23}}&\cdots&T_{\alpha_{2n}}&\dots\\
\end{smallmatrix}\right).
\end{equation}
Set $(T_1,T_2)=(T_{\alpha_{1k}},T_{\alpha_{2k}}\mid 3\leq k)$.
Denote by
$L^\infty(T_1,T_2)=L^\infty(T_{\alpha_{1k}},T_{\alpha_{2k}}\mid
3\leq k)= (T_{\alpha_{1k}},T_{\alpha_{2k}}\mid 3\leq k)''$ the von
Neumann algebra generated by the commuting family of operators
$(T_1,T_2)$. By definition, we have
\begin{equation}
({\mathfrak A}^{2})'=\Big(f\in L^\infty(T_1,T_2)\mid
[f,T_{kk+1}]=0,\,\,3\leq k \Big).
\end{equation}
Using the spectral theorem for the family $(T_1,T_2)$ of commuting unitary operators
$(T_{\alpha_{1k}},T_{\alpha_{2k}}\mid 3\leq k)$
we conclude that any element $f\in L^\infty(T_1,T_2)$ has the following form:
\begin{equation}
\label{f=spectal-int(2)}
f(T_1,T_2)=f(T_{\alpha_{1k}},T_{\alpha_{2k}}\mid 3\leq k)=
\int_{X^2_0}f(\lambda_1,\lambda_2)dE(\lambda_1,\lambda_2)
\end{equation}
where $X^2_0=\prod_{k=3}^\infty({\mathbb F}_p\times{\mathbb
F}_p)_k$, $f$ is essentially bounded function on $X^2_0$ and $E$ is
 common resolution of the identity of the family of operators
$(T_1,T_2)$ defined on cylindrical sets
$\Delta_{13}\times\Delta_{23}\times\cdots\times\Delta_{1n}\times\Delta_{2n}$
as follows:
$$
E(\Delta_{13}\times\Delta_{23}\times\cdots\times\Delta_{1n}\times\Delta_{2n}):=\prod_{k=3}^n
E_{1k}(\Delta_{1k})E_{2k}(\Delta_{2k})
$$
where $E_{rk}$ is  resolution of the identity of the operators
$T_{\alpha_{rk}}$ for $1\leq r\leq 2,\,\,3\leq k$, i.e.,
$T_{rk}=\int_{Sp(T_{rk})}f(\lambda_{rk})dE_{rk}(\lambda_{rk})$.
Similarly to the proof of Lemma~\ref{l.[f,T(k,k+1)]=0(n)}, we get
\begin{lem}
\label{l.[f,T(k,k+1)]=0(inf2)}
An operator $f(T_1,T_2)\in L^\infty(T_1,T_2)$ defined by (\ref{f=spectal-int(2)})
commutes  with  $T_{kk+1}$ for all
$2\leq k$ if and only if for all $2\leq k$ holds:
\begin{equation}
\label{[f,T(k,k+1)]=0(inf2)}
f\left(
 \begin{smallmatrix}
T_{\alpha_{13}},\dots,T_{\alpha_{1k}},T_{\alpha_{1k+1}},\dots,T_{\alpha_{1n}},\dots\\
T_{\alpha_{23}},\dots,T_{\alpha_{2k}},T_{\alpha_{2k+1}},\dots,T_{\alpha_{2n}},\dots
\end{smallmatrix}\right)
=f\left(
 \begin{smallmatrix}
T_{\alpha_{13}},\dots,T_{\alpha_{1k}}T_{\alpha_{1k+1}},T_{\alpha_{1k+1}},\dots,T_{\alpha_{1n}},\dots\\
T_{\alpha_{23}},\dots,T_{\alpha_{2k}}T_{\alpha_{2k+1}},T_{\alpha_{2k+1}},\dots,T_{\alpha_{2n}},\dots
\end{smallmatrix}\right).
\end{equation}
\end{lem}
\begin{lem}
\label{l.f=f(c(T))(inf2)} If for the function $f(T_1,T_2)\in L^\infty(T_1,T_2)$ holds
relation (\ref{[f,T(k,k+1)]=0(inf2)}) then  for some
$(s_1,s_2),\,\,3\leq s_1,\,\,3\leq s_2$ holds
\begin{equation}
\label{f=f(c(T))(inf2)}
f=f_{s_1,s_2}(T_{\alpha_{1s_1}},T_{\alpha_{2s_2}})\prod_{k=s_1+1}^\infty c(T_{\alpha_{1k}})
\prod_{k=s_2+1}^\infty c(T_{\alpha_{2k}}),
\end{equation}
where $f_{s_1,s_2}(T_{\alpha_{1s_1}},T_{\alpha_{2s_2}})\in L^\infty(T_{\alpha_{1s_1}},T_{\alpha_{2s_2}})$.
\end{lem}
\qed\end{pf}
%

\section{The proof of the irreducibilty, case $m=1$ }
\subsection{The irreducibilty, case $m=1,\,\,p=2$}
Let us consider two operators $T_{1n}:=T_{\alpha_{1n}}$ and
$T_{kn}$ on the space $H=H_{\alpha_{1k}}\otimes
H_{\alpha_{1n}}=L^2(X,\mu)$ where $\mu=\mu_{\alpha_{1k}}\otimes
\mu_{\alpha_{1n}}$ and $X={\mathbb F}_p\times {\mathbb F}_p,\,\,2\leq k<n,$
$$
X=
\left(\begin{smallmatrix}
1&x_{1k}&x_{1n}\\
0&1&0\\
0&0&1
\end{smallmatrix}\right).
$$
%
%
The basis in the space $H_{1t}:=H_{\alpha_{1t}}:=L^2({\mathbb
F}_p,\mu_{\alpha_{1t}})$ is $(e_s^{\alpha})_{s\in{\mathbb F}_p}$,
where $e_s^{\alpha}(r)=(\alpha_{1t}(r))^{-1/2}\delta_{sr},$ $s,r\in
{\mathbb F}_p$ hence, the basis in the space
$$
H_{\alpha_{1k}}\otimes H_{\alpha_{1n}}\quad\text{is}\quad
(e_{sr}:=e_s^{\alpha}\otimes e_r^{\alpha})_{s,r\in{\mathbb  F}_p}.
 $$
We fix the lexicographic order on the set $(s,r)_{s,r\in {\mathbb
F}_p}$. So, we have chosen the following basis
$$
e_{00},\quad e_{01},\quad e_{10},\quad e_{11}.
$$
In this basis the operators $T_{1n}$ and $T_{kn}$ act as follows
if the measures $\mu_{\alpha_{1t}}$ are invariant, recall that
$(T_{1n}f)(x_{1n})=f(x_{1n}-1)$ and
$(T_{kn}f)(x_{1k},x_{1n})=f(x_{1k},x_{1n}-x_{1k})$:
\begin{equation}
\label{Tinv1} T_{1n}:e_{ij}\to e_{ij+1},\quad T_{kn}:e_{ij}\to
e_{ij+i}.
\end{equation}
For an arbitrary measure $\mu_\alpha$ operators $T_{1n}$ and $T_{kn}$ act as follows:
\begin{equation}
\label{T1} T_{1n}:e_{ij}\to \sqrt{\frac{ \alpha_{1n}(j) }{
\alpha_{1n}(j+1)} } e_{ij+1},\quad T_{kn}:e_{ij}\to \sqrt{ \frac{
\alpha_{1n}(j)} {
 \alpha_{1n}(j+i)}} e_{ij+i}.
\end{equation}
Using (\ref{Tinv1}) we  have the following transformation of indices
of the basis $e_{ij}$ under the action of $T_{1n}$ and $T_{kn}$:
$$
\left|\begin{smallmatrix}
      &0 & 1& 2& 3\\
    ij&00&01&10&11\\
T_{1n}&01&00&11&10\\
T_{kn}&00&01&11&10
\end{smallmatrix} \right|.
$$
For the invariant measure and in the general case we have respectively
\begin{equation}
\label{T1nT2n-inv}
 T_{1n}=
 \left(\begin{smallmatrix}
 1&0\\
 0&1
\end{smallmatrix}\right)\otimes
\left(\begin{smallmatrix}
 0&1\\
1&0
\end{smallmatrix}\right)=
 \left(\begin{smallmatrix}
 0&1&0&0\\
 1&0&0&0\\
 0&0&0&1\\
 0&0&1&0
\end{smallmatrix}\right),\quad
T_{kn}= \left(\begin{smallmatrix}
1&0&0&0\\
0&1&0&0\\
0&0&0&1\\
0&0&1&0
\end{smallmatrix}\right),
\end{equation}
\begin{equation}
\label{T1nT2n}
 T_{1n}=
 \left(\begin{smallmatrix}
 1&0\\
 0&1
\end{smallmatrix}\right)\otimes
\left(\begin{smallmatrix}
 0&a_n^{-1}\\
a_n&0
\end{smallmatrix}\right)=
 \left(\begin{smallmatrix}
 0&a_n^{-1}&0&0\\
 a_n&0&0&0\\
 0&0&0&a_n^{-1}\\
 0&0& a_n&0
\end{smallmatrix}\right),\quad
T_{kn}= \left(\begin{smallmatrix}
1&0&0&0\\
0&1&0&0\\
0&0&0& a_n^{-1}
\\
0&0&a_n
&0
\end{smallmatrix}\right),
\end{equation}
where $a_n\!=\!\sqrt{\frac{ \alpha_{1n}(0)}{\alpha_{1n}(1)}}$.

Recall that $x_{1k}={\rm diag}(0,1)=
\left(\begin{smallmatrix}
0&0\\
0&1
\end{smallmatrix}\right)$  (see (\ref{X(kn)}) for notations $x_{kn}$).
We would like to approximate the operator $x_{1k} \cong
\dots \otimes x_{1k} \otimes\dots $ (see Remark
\ref{otimes=1}) on the space $H^1\!=\!\otimes_{n=2}^\infty
H_{\alpha_{1n}}$ by linear combinations of operators $T_{kn}$.
\begin{lem}
\label{l.P2} We have
$$
x_{1k} {\bf 1}=
\left(\begin{smallmatrix}
 0&0\\
 0&1
\end{smallmatrix} \right)
{\bf 1}\in \langle (T_{kn}-I){\bf 1} \mid n> k\rangle
\Leftrightarrow S_{11}^L(\mu_\alpha)=\infty.
$$
\end{lem}
\begin{pf} In the space $H_{\alpha_{1k}}\otimes H_{\alpha_{1n}}$ we
have
$$
T_{kn}-I=
\left(\begin{smallmatrix}
0&0&0&0\\
0&0&0&0\\
0&0&-1& a_n^{-1}\\
0&0& a_n &-1
\end{smallmatrix} \right)=
\left(\begin{smallmatrix}
 0&0\\
 0&1
\end{smallmatrix} \right)
\otimes
\left(\begin{smallmatrix}
-1& a_n^{-1}\\
  a_n&-1
\end{smallmatrix} \right)
=x_{1k}
\otimes
(T_{\alpha_{1k}}-I)
$$
so,
\begin{equation}
\label{T2n(1)-1.0}
T_{kn}-I=x_{1k}
\otimes(T_{\alpha_{1n}}-I).
\end{equation}
Hence, we get $\sum_{n=k+1}^{N+k+1}t_n(T_{kn}-I)\to x_{1k} $.
Indeed, in the space $\otimes_{n=k+1}^\infty H_{1n}$ we get
$$
\sum_{n=k+1}^{N+k+1}t_n(T_{kn}-I)= x_{1k} \otimes I
\otimes...\otimes I\otimes
\sum_{n=k+1}^{N+k+1}t_n\left(\begin{smallmatrix}
 -1& a_n^{-1}\\
  a_n&-1
\end{smallmatrix} \right)
\otimes I...
$$
$$
=x_{1k} \otimes I\otimes ...\otimes I\otimes \Big[
\sum_{n=k+1}^{N+k+1}t_n \Big(\begin{smallmatrix}
-1&a_n^{-1}\\
a_n&-1
\end{smallmatrix} \Big)
 \Big]\otimes I...\to x_{1k} \otimes I\otimes\dots\otimes I\dots\cong x_{1k} ,
$$
by Lemma \ref{1in<2>}, where
$a_n=\sqrt{\frac{\alpha_{1n}(0)}{\alpha_{1n}(1)}}$.
\qed\end{pf}
%
Since $x_{1k} ={\rm diag}(0,1)\in {\mathfrak A}^1$ the proof of
the irreducibility for $m=1$
and $p=2$ follows from Remark~\ref{r.Idea-irr}.

\subsection{The irreducibilty, case $m=1$, $p$ is arbitrary}
{\bf Notation.} For an arbitrary $p$ let us denote by $P_{kn}^{(r)}$ the operators
$E_{rr}=$\\${\rm diag}({\underset{r}{\underbrace{0,...,0,1}}},0,...,0)$, acting on the spaces $H_{\alpha_{kn}},\,\,r\in{\mathbb F}_p,\,\,1\leq k<n$.

Let
us suppose that we are able to approximate $P_{1k}^{(r)}$ by the
operators of the representation, i.e., that
$P_{1k}^{(r)}\in{\mathfrak A}^1,\,\,r\in{\mathbb F}_p$ hence, an
operator $x_{1k}$ acting in $H_{1k}$ (see (\ref{X(kn)}) and
(\ref{H(kn)}))
 belongs to ${\mathfrak A}^1$:
$$
x_{1k}={\rm diag}(0,1,\dots,p-1)=\sum_{r\in{\mathbb F}_p}rE_{rr}=
\sum_{r\in{\mathbb F}_p}rP_{1k}^{(r)}\in{\mathfrak A}^1.
$$
In this case the  proof follows from Remark~\ref{r.Idea-irr}.

In order to find an appropriate combinations to approximate the
operators $P_{1k}^{(r)},\,\,r\in{\mathbb F}_p$ {\it we study first
the case} $p=3$. Let us denote (see (\ref{T-alpha}))
\begin{equation}
\label{T,T^2,p=3}
T_\alpha=
\left(\begin{smallmatrix}
0     &0&t_{02}\\
t_{10}&0&0  \\
0&t_{21}&0
 \end{smallmatrix} \right),\,\,\text{ then}\,\,
T^2_\alpha=
\left(\begin{smallmatrix}
0&t_{01}&0\\
0&    0 &t_{12}     \\
t_{20}&0&0
\end{smallmatrix} \right),\,\,\text{ where}\,\,t_{ij}=\sqrt{\frac{
\alpha_{1n}(j)}{\alpha_{1n}(i)}},\,\,i,j\in{\mathbb F}_p.
\end{equation}
Let $e_{kr}:=(e_k^\alpha\otimes e_r^\alpha)_{k,r\in{\mathbb F}_p}$ be the basis in the space $H_{\alpha_{12}}\otimes
H_{\alpha_{1n}}$ (see (\ref{e^al_k-o.n.b.})). Using (\ref{Tinv1}) we  have the following transformation
of indices of the basis $e_{ij}$ under the action of $T_{1n}$ and $T_{2n}$:
$$
\left|
\begin{smallmatrix}
      &0 & 1& 2& 3&4&5&6&7&8\\
    ij&00&01&02&10&11&12&20&21&22\\
T_{1n}&01&02&00&11&12&10&21&22&20\\
T_{2n}&00&01&02&11&12&10&22&20&21
\end{smallmatrix}
\right|.
$$
So, the operators $T_{1n}$ and $T_{2n}$ have the following forms in
$H_{12}\otimes H_{1n}$:
$$
 T_{1n}=
\left(\begin{smallmatrix}
0     &0&t_{02}&0&0&0&0&0&0\\
t_{10}&0&0 &0&0&0&0&0&0    \\
0&t_{21}&0&0&0&0&0&0&0\\
0&0&0&0     &0&t_{02}&0&0&0\\
0&0&0&t_{10}&0&0&0&0&0     \\
0&0&0&0&t_{21}&0&0&0&0\\
0&0&0&0&0&0&0     &0     &t_{02}\\
0&0&0&0&0&0&t_{10}&0     &0     \\
0&0&0&0&0&0&0     &t_{21}&0
\end{smallmatrix} \right),\quad
 T_{2n}=
\left(\begin{smallmatrix}
1     &0&0&0&0&0&0&0&0\\
0&1&0 &0&0&0&0&0&0    \\
0&0&1&0&0&0&0&0&0\\
0&0&0&0     &0&t_{02}&0&0&0\\
0&0&0&t_{10}&0&0&0&0&0     \\
0&0&0&0&t_{21}&0&0&0&0\\
0&0&0&0&0&0&0&t_{01}&0\\
0&0&0&0&0&0&0&    0 &t_{12}     \\
0&0&0&0&0&0&t_{20}&0&0
\end{smallmatrix} \right),
$$
\begin{equation}
\label{T{2n}(1)}
 T_{1n}=
 \left(\begin{smallmatrix}
 1&0&0\\
 0&1&0\\
 0&0&1
\end{smallmatrix} \right)
\otimes
\left(\begin{smallmatrix}
0     &0&t_{02}\\
t_{10}&0&0     \\
0&t_{21}&0
\end{smallmatrix} \right)=
\left(\begin{smallmatrix}
 T_{\alpha} &0&0\\
 0&T_{\alpha} &0\\
 0&0&T_{\alpha}
\end{smallmatrix} \right),\quad
T_{2n}=
\left(\begin{smallmatrix}
 T_{\alpha}^0 &0&0\\
 0&T_{\alpha}^1 &0\\
 0&0&T_{\alpha}^2
\end{smallmatrix} \right),
\end{equation}
where $\alpha=\alpha_{1n}$. Note that
\begin{equation}
\label{Tkn=diag(.,.,.)}
T_{1n}={\rm
diag}(T_\alpha,T_\alpha,T_\alpha),\quad T_{2n}={\rm diag}(I,T_\alpha,T_\alpha^2).
\end{equation}
Since
$$
T_{2n}={\rm diag}(I,T_\alpha,T_\alpha^2),\,\,\, T_{2n}^2={\rm
diag}(I,T_\alpha^2,T_\alpha)\,\,\,\,\text{so},\,\,\,\,C(T_{2n})= {\rm
diag}(3,C(T_\alpha),C(T_\alpha)).
$$
Similarly, we get
$$
T_{1n}T_{2n}={\rm diag}(T_\alpha,T_\alpha^2,I),\quad (T_{1n}T_{2n})^2={\rm diag}(T_\alpha^2,T_\alpha,I),
$$
$$
C(T_{1n}T_{2n})= {\rm diag}(C(T_\alpha),C(T_\alpha),3),
$$
$$
T_{1n}^2T_{2n}={\rm diag}(T_\alpha^2,I,T_\alpha),\quad(T_{1n}^2T_{2n})^2={\rm diag}(T_\alpha,I,T_\alpha^2),\,
$$
$$
C(T_{1n}^2T_{2n})= {\rm diag}(C(T_\alpha),3,C(T_\alpha)).
$$
So, we can try to approximate
$$
{\rm diag}(0,I,I)=(P_{12}^{(1)}+P_{12}^{(2)})\otimes
I=(I-P_{12}^{(0)})\otimes I \quad\text{by combinations of
}\,\,C(T_{2n})-3,
$$
$$
{\rm diag}(I,I,0)=(P_{12}^{(0)}+P_{12}^{(1)})\otimes
I=(I-P_{12}^{(2)})\otimes I \quad\text{by combinations of
}\,\,C(T_{1n}T_{2n})-3,
$$
$$
{\rm diag}(I,0,I)=(P_{12}^{(0)}+P_{12}^{(2)})\otimes
I=(I-P_{12}^{(1)})\otimes I \quad\text{by combinations of
}\,\,C(T_{1n}^2T_{2n})-3.
$$
In the general case, we can try to approximate
$$
(I-P_{1k}^{(p-r)})\otimes I\quad\text{by combinations of
}\,\,\sum_{s\in{\mathbb
F}_p}(T_{1n}^{r}T_{kn})^s-p=C(T_{1n}^{r}T_{kn})-p.
$$
\begin{lem}
\label{l.P_r} We have for $r\in{\mathbb F}_p$ and $k>1$
$$
(I-P_{1k}^{(p-r)}){\bf 1}\in \langle \Big[C(T_{1n}^{r}T_{kn})-p\Big]{\bf 1}\mid n> k \rangle
$$
if and only if $S_{11}^L(\mu_\alpha)=\infty \Leftrightarrow
\mu^\alpha\perp \mu_{inv}.$
\end{lem}
\begin{pf} Since $T_{1n}\!=\!{\rm
diag}(T_{\alpha_{1n}},\dots,T_{\alpha_{1n}}),$
$ T_{2n}={\rm
diag}(I,T_{\alpha_{1n}},T_{\alpha_{1n}}^2,\dots,T_{\alpha_{1n}}^{p-1}),$\\$
T_{2n}^s={\rm
diag}(I,T_{\alpha_{1n}}^s,T_{\alpha_{1n}}^{2s},\dots,T_{\alpha_{1n}}^{s(p-1)}),\,\,s\in{\mathbb
F}_p,
$
 we get
 $$
(T_{1n}^{r}T_{2n})^s=\Big[
{\rm diag}(T_{1n},T_{1n},T_{1n},\dots,T_{1n}){\rm diag}(I,T_{1n},T^2_{1n},\dots,T^{p-1}_{1n})
\Big] ^s=
 $$
 $$
 \Big[{\rm diag}(T_{1n}^r,T_{1n}^{r+1},T^{r+2}_{1n},\dots,T^{r+p-1}_{1n})\Big] ^s=
 \Big[{\rm diag}(T_{1n}^{rs},T_{1n}^{(r+1)s},T^{r+2}_{1n},\dots,T^{(r+p-1)s}_{1n})\Big].
 $$
 Therefore,
 $$
 \sum_{s\in{\mathbb F}_p}(T_{1n}^{r}T_{2n})^s=\Big(
\sum_{s\in{\mathbb F}_p} T_{1n}^{rs},\sum_{s\in{\mathbb F}_p}T_{1n}^{(r+1)s},\dots,\sum_{s\in{\mathbb F}_p}T^{(r+p-1)s}_{1n})
 \Big)
 $$
 $$
 =\Big(
 {\underset{p-r}{\underbrace{C(T_{1n}),C(T_{1n}),\dots,p}}},\dots,C(T_{1n}).
 \Big)
 $$
 At last, we get $\sum_{s\in{\mathbb F}_p}(T_{1n}^{r}T_{2n})^s-p=$
 \begin{align*}
\Big(C(T_{1n})-p\Big)
( {\underset{p-r}{\underbrace{I,I,\dots,I,0}}},I,\dots,I)=(I-P_{12}^{(p-r)})\otimes \Big(C(T_{1n})-p\Big).
 \end{align*}
 Finally, we get when  $N\to\infty$
  \begin{align*}
\sum_{n=k+1}^{k+N}\!t_n\Big[
C(T_{1n}^{r}T_{2n})-p\Big]\!-\!(I\!-\!P_{12}^{(p-r)})\otimes I\!=\!\\
(\!I-\!P_{12}^{(p-r)})\otimes \Big[\!\!\sum_{n=k+1}^{k+N}\!t_n \Big(
C(T_{1n})\!-\!p\Big)\!-\!I\Big]\!\to\!0.
 \end{align*}
The proof of the latter statement is similar to the proof of Lemma
\ref{l.P2}.
\qed\end{pf}
Finally, we can approximate $P_{1k}^{(r)}$ therefore, $x_{1k}=\sum_{r\in {\mathbb F}_p}rP_{1k}^{(r)}=\sum_{r\in {\mathbb F}_p}rE_{rr}
\in{\mathfrak A}^1$.  Using the Remark~\ref{r.Idea-irr}, we conclude that the representation is irreducible.

\section{Irreducibility, case $m=2$}
Let us consider three operators $T_{1n},\,\,T_{2n}$ and
$T_{kn}$  on the space $H=
L^2(X,\mu)=H_{\alpha_{12}}\otimes H_{\alpha_{1n}}\otimes
H_{\alpha_{2n}}$ where $\mu=\mu_{\alpha_{12}}\otimes
\mu_{\alpha_{1n}}\otimes \mu_{\alpha_{2n}}$ and
$$
 X=
\left(\begin{smallmatrix}
1&x_{12}&x_{1n}\\
0& 1    &x_{2n}\\
0& 0    & 1
\end{smallmatrix} \right).
$$
The basis in the space $H_{st}=H_{\alpha_{st}}=L^2({\mathbb
F}_p,\mu_{\alpha_{st}})$ is $(e_k^\alpha)_{k\in{\mathbb F}_p}$
(see (\ref{e^al_k-o.n.b.})) hence, the basis in the space
$H_{\alpha_{12}}\otimes H_{\alpha_{1n}}\otimes H_{\alpha_{2n}}$ is
$(e_{krs}:=e_t^\alpha\otimes e_r^\alpha\otimes
e_s^\alpha)_{t,r,s\in{\mathbb F}_p}$. We fix the lexicographic order on the set $(t,r,s)_{k,r,s\in
{\mathbb F}_p}$. So, we have chosen the following basis
$$
e_{000},\,\,e_{001},\,\,e_{010},\,\,e_{011},
\,\,e_{100},\,\,e_{101},\,\,e_{110},\,\,e_{111}.
$$
In this basis the operators $T_{1n}$ and $T_{2n}$ act as follows
if the measures $\mu_{\alpha_{st}}$ are invariant ($T_{kn}$ acts on the space $H_{1k}\otimes H_{1n}\otimes H_{2k}\otimes H_{2n}$):
\begin{equation}
\label{Tinv2} T_{1n}:e_{ijl}\to e_{i,j+1,l},\,\, T_{2n}:e_{ijl}\to
e_{i,j+i,l+1},\,\, T_{kn}:e_{ijlr}\to e_{i,j+i,l,r+l}
\end{equation}
and as follows if the measure is not invariant:
\begin{equation}
\label{T2} T_{1n}:e_{ijl}\to \sqrt{\frac{ \alpha_{1n}(j) }{
\alpha_{1n}(j+1)} } e_{i,j+1,l},\,\, T_{2n}:e_{ijl}\to \sqrt{
\frac{
\alpha_{1n}(j)\alpha_{2n}(l)}{\alpha_{1n}(j+i)\alpha_{2n}(l+1)}}
e_{i,j+i,l+1},
\end{equation}
\begin{equation}
\label{T_kn} T_{kn}:e_{ijlr}\to \sqrt{ \frac{
\alpha_{1n}(j)\alpha_{2n}(r)}{\alpha_{1n}(j+i)\alpha_{2n}(r+l)}} e_{i,j+i,l,r+l}.
\end{equation}
Using (\ref{Tinv2}) we  have the following transformation of indices
of the basis $e_{ijl}$ under the action of $T_{1n}$ and $T_{2n}$:
$$
\left|
\begin{smallmatrix}
      &0  &  1&  2&  3&  4&  5&  6& 7 \\
  ijl &000&001&010&011&100&101&110&111\\
T_{1n}&010&011&000&001&110&111&100&101\\
T_{2n}&001&000&011&010&111&110&101&100
\end{smallmatrix}
\right|,
\quad\text{i.e.},\quad
 \left|\begin{smallmatrix}
  ijl      &0  &  1&  2&  3&  4&  5&  6& 7 \\
T_{1n}&2&3&0&1&6&7&4&5\\
T_{2n}&1&0&3&2&7&6&5&4
\end{smallmatrix} \right|.
$$
So, the operators $T_{1n}$ and $T_{2n}$ have the following forms:
$$
T_{1n}=\left(\begin{smallmatrix}
1&0\\
0&1
\end{smallmatrix} \right)\otimes
\left(\begin{smallmatrix}
0&a_n^{-1}\\
a_n&0
\end{smallmatrix} \right)\otimes
\left(\begin{smallmatrix}
1&0\\
0&1
\end{smallmatrix} \right)=
\left(\begin{smallmatrix}
0&0&a^{-1}&  0&  0&  0& 0&0\\
0&0&0&a^{-1}&  0&  0& 0&0\\
a& 0& 0&  0&0&  0& 0&0\\
0&a& 0& 0&  0&0&  0& 0\\
0& 0& 0&0&0&0&a^{-1}&0\\
0&0& 0& 0&0&0&0&a^{-1}\\
0&  0& 0&0&a& 0& 0&  0\\
0&  0& 0&0&0&a & 0& 0
\end{smallmatrix} \right),\,\,
$$
$$
T_{2n}=
\left(\begin{smallmatrix}
1&0&0&0\\
0&1&0&0\\
0&0&0&a_n^{-1}\\
0&0&a_n&0
\end{smallmatrix} \right)
\otimes
\left(\begin{smallmatrix}
0&b_n^{-1}\\
b_n&0
\end{smallmatrix} \right)=
\left(\begin{smallmatrix}
0&b^{-1}&0
&  0&  0&  0& 0&0\\
b&0&0& 0
&  0&  0& 0&0\\
0& 0& 0&  b^{-1}&0&  0& 0&0\\
0&0& b& 0&  0&0&  0& 0\\
 0&0& 0& 0&0&0&0&a^{-1}b^{-1}\\
 0&0& 0& 0&0&0& a^{-1}b&0\\
 0&0& 0& 0&0&ab^{-1}&0&0\\
0&0& 0& 0&ab &0&0&0
\end{smallmatrix} \right),
$$
where
\begin{equation}
\label{8.ab}
 a=a_n=\sqrt{\frac{ \alpha_{1n}(0) }{ \alpha_{1n}(1)} },\quad
b=b_n=\sqrt{\frac{ \alpha_{2n}(0) }{ \alpha_{2n}(1)} }.
\end{equation}
To calculate $T_{1n}$ and $T_{2n}$ we use the following formulas
\begin{equation}
\label{A-otimes- B}
\left(\begin{smallmatrix}
a&b\\
c&d
\end{smallmatrix} \right)\otimes
\left(\begin{smallmatrix}
1&0\\
0&1
\end{smallmatrix} \right)=
\left(\begin{smallmatrix}
a&0&b&0\\
0&a&0&b\\
c&0&d&0\\
0&c&0&d
\end{smallmatrix} \right),\quad
\left(\begin{smallmatrix}
1&0\\
0&1
\end{smallmatrix} \right)\otimes
\left(\begin{smallmatrix}
a&b\\
c&d
\end{smallmatrix} \right)=
\left(\begin{smallmatrix}
a&b&0&0\\
c&d&0&0\\
0&0&a&b\\
0&0&c&d
\end{smallmatrix} \right).
\end{equation}
\begin{rem}
\label{r.A-otimes-B}
The latter formulas are particular case of the formulas below.
 For $A=(a_{ij})_{ij}\in {\rm Mat}(n,{\bf k})$ and
$B=(b_{rs})_{rs}\in {\rm Mat}(m,{\bf k})$ we have in the basis
$e_i\otimes f_s$
\begin{equation*}
 A\otimes B=(a_{ij}b_{rs})_{((i,j),(r,s))}\in
{\rm End}({\bf k}^n\otimes{\bf k}^m).
\end{equation*}
Indeed, we get
$$
Ae_j=\sum_{i=1}^na_{ij}e_i,\quad Bf_s=\sum_{r=1}^mb_{rs}f_r,
$$
therefore, we have
$$
(A\otimes B)e_i\otimes f_s=Ae_i\otimes Bf_s=\sum_{i=1}^n\sum_{r=1}^ma_{ij}b_{rs}e_i\otimes f_r.
$$
\end{rem}
\subsection{Irreducibility $m=2$, $p=2$}
For two sequence $a=(a_n)_n$ and $b=(b_n)_n$ of positive numbers we say that they are {\it equivalent} are write $a_n\sim b_n$ if $C_1a_n\leq b_n\leq C_2a_n$ for all $n\in{\mathbb N}$.
Recall (see Notations before  Lemma~\ref{1in<2>}) that
$c_{kn}=2\sqrt{\alpha_{kn}(0)\alpha_{kn}(1)}$ and $
S_{kn}^L(\mu_\alpha)$ are defined as follows  (see (\ref{S^L(kk)}) and (\ref{S^L(kn)}))
$$
S_{11}^L(\mu_\alpha)\!=\!
\sum_{n=2}^\infty\left(1\!-\!c_{1n}\right),\,\,
S_{12}^L(\mu_\alpha)\!=\!
\sum_{n=3}^\infty\alpha_{2n}(1)\left(1\!-\!c_{1n}\right),\,\,
 S_{22}^L(\mu_\alpha)\!=\!
\sum_{n=3}^\infty\left(1\!-\!c_{2n}\right).
$$
\begin{thm}
Representation $T^{R,\mu,2}$ is irreducible if and only if\\
1) $(\mu_\alpha)^{L_{I+E_{12}}}\perp \mu_\alpha  \Leftrightarrow   S_{12}^L(\mu_\alpha)=\infty,$\\
2) $\mu_\alpha^2\perp \mu_{inv}^2\Leftrightarrow S_{22}^L(\mu_\alpha)=\infty.$
\end{thm}
Let $p=2$. To approximate $x_{1k}$ and $x_{2k}$ we use
the following expressions (see (\ref{Reg=otimes(r)}))
\begin{align}
\label{61}
 T_{kn}-t_{kn}(2)&=T_{kn}(1)\otimes
(T_{kn}(2)-t_{kn}(2))+ (T_{kn}(1)-I)\otimes t_{kn}(2)\\
&=T_{kn}(1)\otimes T_{kn}^c(2)+\hat{T}_{kn}(r)\otimes t_{kn}(2), \nonumber\\
\label{62}
T_{kn}-t_{kn}(1)&=(T_{kn}(1)-t_{kn}(1))\otimes
T_{kn}(2)+ t_{kn}(1)\otimes (T_{kn}(2)-I)\\
&=T_{kn}^c(r)\otimes T_{kn}(2)+t_{kn}(1)\otimes \hat{T}_{kn}(r),\nonumber\\
\label{63}
T_{kn}(1)-I&=x_{1k}\otimes
(T_{\alpha_{1n}}-I),\quad T_{kn}(2)-I=x_{2k}\otimes
(T_{\alpha_{2n}}-I),
\end{align}
  where  $T_{kn}=T_{kn}(1)\otimes T_{kn}(2)$ (see (\ref{Tkn(1)(2)})) and we set
\begin{equation}
\label{t(kn)(r)}
t_{kn}(r)=(T_{kn}(r){\bf 1},{\bf 1}),\quad
T_{kn}^c(r)=T_{kn}(r)-t_{kn}(r),\quad
\hat{T}_{kn}(r)=T_{kn}(r)-1.
\end{equation}
{\bf Notation}. Recall that we denote by  $x_{kn}$ the operator
$x_{kn}=\left(\begin{smallmatrix}
0&0\\
0&1
\end{smallmatrix}\right)=E_{11}$
on the space $H_{kn}=L^2({\mathbb F}_2,\mu_{\alpha_{kn}})$, see Remark~\ref{r.Idea-irr}.

We find  condition of the approximation of the operator
$x_{12}$ (respectively $x_{1k}$ and
$x_{2k},\,k>2$) by the following linear combinations
$$
\sum_{n}t_n(T_{2n}-t_{2n}(2))\,\big(\text{respectively by}\,
\sum_{n}t_n(T_{kn}-t_{kn}(2))\,\text{and}\,
\sum_{n}t_n(T_{kn}-t_{kn}(1))\big).
$$
We need the following lemma.
\begin{lem}
\label{l.(D(l,c)m,m)=}
For fixed $\lambda,\,\,c,\,\,\mu\in {\mathbb R}^m$ and the matrix $D_m(\lambda,c)$ defined as follows:
\begin{equation}
\label{D(lambda,c)}
D_m(\lambda,c)= \left(
\begin{array}{cccc}
1+\lambda_1+c_1^2&1+c_1c_2          &...&1+c_1c_m\\
1 +c_2c_1          &1+\lambda_2+c_2^2&...&1+c_2c_m\\
&&...&\\
1 +c_mc_1          &          1+c_mc_2 &...&1+\lambda_m+c_m^2
\end{array}
\right)
\end{equation}
we have
\begin{equation}
\label{D(l,c)=G()}
(D^{-1}_m(\lambda,c)\mu,\mu)=
\frac{\Gamma(f_m)+\Gamma(f_m,g_m)+\Gamma(f_m,h_m)}{1+\Gamma(g_m)+\Gamma(h_m)+\Gamma(g_m,h_m)}
\end{equation}
\begin{equation}
\label{D(l,c)mumu}
=\frac{\sum_{k=1}^m\frac{\mu_k^2}{\lambda_k}+\sum_{1\leq k<n\leq
m}\frac{(\mu_k-\mu_n)^2}{\lambda_k\lambda_n}+\sum_{1\leq k<n\leq
m}\frac{(c_n\mu_k-c_k\mu_n)^2}{\lambda_k\lambda_n}
}{
1+\sum_{k=1}^m\frac{1}{\lambda_k}+\sum_{k=1}^m\frac{c_k^2}{\lambda_k}+
\sum_{1\leq k<n\leq
m}\frac{(c_k-c_n)^2}{\lambda_k\lambda_n}
},
\end{equation}
where
\begin{equation}
\label{f,g,h}
f_m=\Big(\frac{\mu_k}{\sqrt{\lambda_k}}\Big)_{k=1}^m,\quad
g_m=\Big(\frac{1}{\sqrt{\lambda_k}}\Big)_{k=1}^m,\quad
h_m=\Big(\frac{c_k}{\sqrt{\lambda_k}}\Big)_{k=1}^m.
\end{equation}
\end{lem}
\begin{pf} We have for $m=2$
$$
D_2(\lambda,c)= \left(
\begin{array}{cc}
1+\lambda_1+c_1^2&1+c_1c_2 \\
1 +c_2c_1          &1+\lambda_2+c_2^2
\end{array}
\right),
$$
$$
F_2(\lambda,c):={\rm
det}\,D_2(\lambda,c)\!=\!1\!+\!\lambda_1\!+\!\lambda_2\!+\!c_1^2\!+\!c_2^2
\!+\!(\lambda_1+c_1^2)(\lambda_2+c_2^2)\!-\!(1+2c_1c_2+c_1^2c_2^2)
$$
$$
=\lambda_1\lambda_2\Big(1\!+\!\frac{1}{\lambda_1}\!+\!\frac{1}{\lambda_2}\!+\!
\frac{c_1^2}{\lambda_1}\!+\!\frac{c_2^2}{\lambda_2}\!+\!
\frac{(c_1-c_2)^2}{\lambda_1\lambda_2}\Big)\!=\!\lambda_1\lambda_2(1+\Gamma(g_2)+\Gamma(h_2)+\Gamma(g_2,h_2)),
$$
where
$$
g_2=\Big(\frac{1}{\sqrt{\lambda_1}},\frac{1}{\sqrt{\lambda_2}}\Big),\quad
h_2=\Big(\frac{c_1}{\sqrt{\lambda_1}},\frac{c_2}{\sqrt{\lambda_2}}\Big).$$
In the general case, we show that
\begin{equation}
\label{detD(l,c)}
F_m(\lambda,c):={\rm det}\,D_m(\lambda,c)=\prod_{k=1}^m\lambda_k
\Big(
1+\sum_{k=1}^m\frac{1}{\lambda_k}+\sum_{k=1}^m\frac{c_k^2}{\lambda_k}+
\sum_{1\leq k<n\leq
m}\frac{(c_k-c_n)^2}{\lambda_k\lambda_n}
\Big).
\end{equation}
We prove (\ref{detD(l,c)}) first for $m=3$:
$$
F_3(\lambda,c)=
\left|
\begin{array}{ccc}
1+\lambda_1+c_1^2&1+c_1c_2           &1+c_1c_3\\
1 +c_2c_1          &1+\lambda_2+c_2^2&1+c_2c_3\\
1 +c_3c_1          &1+c_3c_2         &1+\lambda_3+c_3^2
\end{array}
\right|=
$$
$$
\prod_{k=1}^3\lambda_k\Big(
1+\sum_{k=1}^3\frac{1+c_k^2}{\lambda_k}+
\sum_{1\leq k<n\leq
3}\frac{(c_k-c_n)^2}{\lambda_k\lambda_n}\Big)=\lambda_1\lambda_2\lambda_3+\lambda_2\lambda_3(1+c_1^2)+
$$
$$
\lambda_1\lambda_3(1+c_2^2)+
\lambda_1\lambda_2(1+c_3^2)+\lambda_3(c_1-c_2)^2+\lambda_2(c_1-c_3)^2+\lambda_1(c_2-c_3)^2.
$$
To prove the latter equality we show that the appropriate derivatives for the both sides of the relation coincide.
Indeed, they are equal respectively:
$$
\frac{\partial F_3(\lambda, c)}{\partial \lambda_i}\mid_{\lambda=0}=(c_j-c_k)^2,\quad
\frac{\partial^2 F_3(\lambda, c)}{\partial \lambda_i\lambda_j}\mid_{\lambda=0}=1+c_k^2,\quad
\frac{\partial^3 F_3(\lambda, c)}{\partial \lambda_1\lambda_2\lambda_3}\mid_{\lambda=0}=1,
$$
where $i,j,k$ is the cyclic permutations of the indices  $1,2,3$. In the general case, we have for both sides of the equation (\ref{detD(l,c)})
 {\small
$$
\frac{\partial^{m-2} F_m(\lambda, c)}{\partial \lambda_1\dots \lambda_{m-2}}\mid_{\lambda=0}\!=\!(c_{m-1}-c_m)^2,\,\,
\frac{\partial^{m-1} F_m(\lambda, c)}{\partial \lambda_1\dots\lambda_{m-1}}\mid_{\lambda=0}\!=\!1+c_m^2,\,\,
\frac{\partial^m     F_m(\lambda, c)}{\partial \lambda_1\dots\lambda_m}\mid_{\lambda=0}\!=\!1
$$
}
and the corresponding cyclic permutations of the indices. This proves (\ref{detD(l,c)}).

Further, for $m=2$ we have
$$
D_2^{-1}(\lambda,c)= \frac{1}{F_2(\lambda,c)}\left(
\begin{array}{cc}
1+\lambda_2+c_2^2&-(1+c_1c_2) \\
-(1 +c_2c_1)          &1+\lambda_1+c_1^2
\end{array}
\right)
\quad\text{and}\quad
%
%
(D_2^{-1}(\lambda,c)\mu,\mu)=
$$
$$
(F_2(\lambda,c))^{-1}\Big[
(1+\lambda_2+c_2^2)\mu_1^2
-2(1+c_1c_2)\mu_1\mu_2+
(1+\lambda_1+c_1^2)\mu_2^2
\Big]
$$
$$
=(F_2(\lambda,c))^{-1}\Big[
(\mu_1-\mu_2)^2+(c_2\mu_1-c_1\mu_2)^2+\lambda_2\mu_1^2+\lambda_1\mu_2^2
\Big]=
$$
$$
\frac{
\lambda_1\lambda_2\Big(\frac{\mu_1^2}{\lambda_1}+\frac{\mu_2^2}{\lambda_2}+
\frac{(\mu_1-\mu_2)^2}{\lambda_1\lambda_2}+
\frac{(c_2\mu_1-c_1\mu_2)^2}{\lambda_1\lambda_2}
\Big)
}
{
\lambda_1\lambda_2\Big(1+\frac{1}{\lambda_1}+\frac{1}{\lambda_2}+
\frac{c_1^2}{\lambda_1}+\frac{c_2^2}{\lambda_2}+
\frac{(c_1-c_2)^2}{\lambda_1\lambda_2}\Big)
}=
\frac{\Gamma(f_2)+\Gamma(f_2,g_2)+\Gamma(f_2,h_2)}{1+\Gamma(g_2)+\Gamma(h_2)+\Gamma(g_2,h_2)}.
$$
By analogy, we get (\ref{D(l,c)mumu}) for the genaral $m$.
\qed\end{pf}

\begin{lem}
\label{l.sigma1k-m=2}
 We have $ x_{1k}{\bf 1}\in
\langle\left[T_{kn}-t_{kn}(2)\right]{\bf 1}\mid n> k\rangle
\Leftrightarrow \Sigma_{1m}\to\infty $ where
\begin{equation}
\label{appr-P1}
\Sigma_{1m}=
(D^{-1}_m(\lambda,c)\mu,\mu)=
\frac{\Gamma(f_m)+\Gamma(f_m,g_m)+\Gamma(f_m,h_m)}{1+\Gamma(g_m)+\Gamma(h_m)+\Gamma(g_m,h_m)},
\end{equation}
the vectors $f_m,\,g_m,h_m$ are defined by (\ref{f,g,h}) and
$$
\lambda_n\!=\!\frac{2(1-c_{1n}^2+1-c_{2n}^2+1-c_{1n}^2c_{2n}^2)}{(1-c_{2n})^2},\,\,
c_n\!=\!c_{1n},\,\,\mu_n\!=\!-\frac{\sqrt{2}(1+c_{2n})(1-c_{1n})}{(1-c_{2n})}.
$$
\end{lem}
\begin{pf}
We have $T_{kn}=T_{kn}(1)\otimes T_{kn}(2)$ where an operator
$T_{kn}(1)$ acts on $H_{1k}\otimes H_{1n}$,  $T_{kn}(2)$ acts on
$H_{2k}\otimes H_{2n}$ and they are defined by
\begin{equation}
\label{Tkn(1)(2)}
 T_{kn}(1)\!=\!
\left(\begin{smallmatrix}
1&0&0&0\\
0&1&0&0\\
0&0&0&a_{n}^{-1}\\
0&0&a_{n}&0
\end{smallmatrix}\right),\,
T_{kn}(2)\!=\!
\left(\begin{smallmatrix}
1&0&0&0\\
0&1&0&0\\
0&0&0&b_{n}^{-1}\\
0&0&b_{n}&0
\end{smallmatrix}\right)
,\,a_n\!=\!\sqrt{\frac{\alpha_{1n}(0)}{\alpha_{1n}(1)}},\,
b_n\!=\!\sqrt{\frac{\alpha_{2n}(0)}{\alpha_{2n}(1)}}.
\end{equation}
Using (\ref{61}) we get
$$
T_{kn}-t_{kn}(2)=T_{kn}(1)\otimes
(T_{kn}(2)-t_{kn}(2))+(T_{kn}(1)-I)\otimes t_{kn}(2).
$$
Set $$ a_n^{(r)}:=(T_{\alpha_{rn}}-I){\bf
1},\,\,r=1,2\quad\text{and}\quad b=(b_n),\,\,
b_n=t_{kn}(2)Ma_n^{(1)}.
$$
Take $t=(t_n)_n$ such that $(t,b)=\sum_{n}t_nt_{kn}(2)Ma_n^{(1)}=1$.
Set
$$
f_n\!=\!\big[T_{kn}(1)\otimes (T_{kn}(2)-t_{kn}(2))\big]{\bf 1},\,\,
$$
$$
g_n\!=\!\big[(T_{kn}(1)-I)\otimes t_{kn}(2)\big]{\bf 1}\!=\!\big[x_{1k}\otimes
(T_{\alpha_{1n}}-I)\otimes t_{kn}(2)\big]{\bf 1},
$$
$$
g_n^c=\big[g_n-x_{1k}\otimes Ma_{n}^{(1)}\otimes t_{kn}(2)\big] {\bf
1}= \big[x_{1k}\otimes\big(T_{\alpha_{1n}}-c_{1n}\big)\otimes
t_{kn}(2)\big]{\bf 1}.
$$
We use the relation
\begin{equation}
\label{T-alpha-1n^c} T_{\alpha_{1n}}-1-Ma_n^{(1)}
=T_{\alpha_{1n}}-1-(c_{1n}-1)=T_{\alpha_{1n}}-c_{1n}.
\end{equation}
Using (\ref{61}--\ref{63}) we have $(T_{kn}-t_{kn}(2)){\bf
1}=f_n+g_n$ therefore, we get
$$
\Vert \big[ \sum_{n}t_n(T_{kn}-t_{kn}(2))-x_{1k}\big]{\bf 1}
\Vert^2= \Vert\sum_{n}t_n(f_n+g_n)-x_{1k}\otimes
\sum_{n}t_nt_{kn}(2)Ma_n^{(1)}{\bf 1} \Vert^2
$$
$$
=\Vert \sum_{n}t_nT_{kn}(1)\otimes (T_{kn}(2)-t_{kn}(2)){\bf 1}
+x_{1k}\otimes
\sum_{n}t_nt_{kn}(2)\big(T_{\alpha_{1n}}-c_{1n}\big){\bf 1}
\Vert^2
$$
$$
=\Vert \sum_{n}t_n(f_n+g_n^c)\Vert^2
=\sum_{n,m}t_nt_m(h_n,h_m)=\sum_{n,m}t_nt_m a_{nm}=(At,t),
$$
where
$$
h_n=f_n+g_n^c,\quad A=(a_{nm})_{n,m},\quad a_{nm}=(h_n,h_m).
$$
We use the following estimation
for a positively definite
operator $A$ acting on a space $H$ and a vector $b\in H$ (see \cite{Kos04}):
$$
\min_{t\in H}\Big((At,t)\mid (t,b)=1\Big)=(A^{-1}b,b)^{-1}.
$$
The minimum is reached  for
$x=A^{-1}b\left((A^{-1}b,b)\right)^{-1}$.

We calculate $a_{nm}=(h_n,h_m)$ and we show that $b_n=-\frac{1+c_{2n}}{2}(1-c_{1n})$ and
\begin{equation}
\label{a(nn)}
a_{nn}\!=
\!1\!-\!\big((1\!+\!c_{2n})/2\big)^2\big(1\!+\!c_{1n}^2)/2,
\end{equation}
\begin{equation}
\label{a(nm)}
a_{nm}\!=\!(1\!+\!c_{1n}c_{1m})(1\!-\!c_{2n})(1\!-\!c_{2m})/8.
\end{equation}
Since $(f_n,g_n^c)=0$ we get $\Vert f_n+ g_n^c\Vert^2=\Vert f_n\Vert^2+ \Vert g_n^c\Vert^2$.
 Indeed,
$$
(f_n,g_n^c)=([T_{kn}(1)\otimes (T_{kn}(2)-t_{kn}(2))]{\bf 1},
x_{1k}\otimes\big(T_{\alpha_{1n}}-c_{1n}\big)\otimes
t_{kn}(2){\bf 1})
$$
$$
= t_{kn}(2)(T_{kn}(1)x_{1k}\otimes
\big(T_{\alpha_{1n}}-c_{1n}\big){\bf 1},{\bf 1})((T_{kn}(2)-t_{kn}(2)){\bf 1},{\bf 1})=0.
$$
We have
\begin{equation}
\label{a(nn)'}
a_{nn}\!
=\!\Vert f_n\Vert^2+ \Vert g_n^c\Vert^2\!=\!
\Vert (T_{kn}(2)-t_{kn}(2)) {\bf 1}\Vert^2\!+\!t_{kn}^2(2)\Vert x_{1k}{\bf 1}\Vert^2
\Vert(T_{\alpha_{1n}}-c_{1n}) {\bf 1}\Vert^2.
\end{equation}
For general  $f,g\in L^1(X,\mu)\cap L^2(X,\mu)$ we use the following relation:
$$
(f-Mf,g-Mg)=(f,g)-Mf\overline{Mg}\quad\text{where}\quad Mf=\int_Xf(x)d\mu(x).
$$
In what follows we use the fact that $c_{kn}=2\sqrt{\alpha_{kn}(0)\alpha_{kn}(1)}\in (0,\,1]$. We assume that
$\alpha_{rk}(s)=1/2,\,\, r=1,2,\,\,s\in{\mathbb F}_2$. Obviously,
\begin{equation}
\label{64} Ma_n^{(r)}=c_{rn}-1,\quad
t_{kn}(r)=\alpha_{1k}(0)+c_{rn}\alpha_{1k}(1)=2^{-1}(1+c_{rn}),\,\,\,r=1,2,\,\,k>2.
\end{equation}
Indeed, we get for $r=1$
$$
t_{kn}(1)\!=\!(T_{kn}(1){\bf 1},{\bf 1})\!=\!\left(\!\!\left(\begin{smallmatrix}
1&0&0&0\\
0&1&0&0\\
0&0&0&a_{n}^{-1}\\
0&0&a_{n}&0
\end{smallmatrix}\!\!\right)\!{\bf 1},{\bf 1}\!\!\right)\!\!\stackrel{(\ref{Tkn(1)(2)})}{=}\!\!((1,1,a_{n}^{-1},a_{n}),(1,1,1,1))\!
\stackrel{(\ref{(f,g)})}{=}\!
$$
$$
\alpha_{1k}(0)\alpha_{1n}(0)\!+\!\alpha_{1k}(0)\alpha_{1n}(1)\!+\!
\sqrt{\frac{\alpha_{1n}(1)}{\alpha_{1n}(0)}}\alpha_{1k}(1)\alpha_{1n}(0)\!+\!
\sqrt{\frac{\alpha_{1n}(0)}{\alpha_{1n}(1)}}\alpha_{1k}(1)\alpha_{1n}(1)
$$
$$
\!=\!\alpha_{1k}(0)+c_{1n}\alpha_{1k}(1).
$$
Further, we have
$$
\Vert (T_{kn}(2)-t_{kn}(2)) {\bf 1}\Vert^2=1-t_{kn}^2(2)=1-\Big(\frac{1+c_{2n}}{2}\Big)^2,
$$
$$
\Vert (T_{2n}(2)-t_{2n}(2)) {\bf 1}\Vert^2=1-c_{2n}^2,
$$
\begin{equation}
\label{Mx_{kn}1}
\Vert x_{1k}{\bf 1} \Vert^2=\Vert \left(\begin{smallmatrix}
0&0\\
0&1
\end{smallmatrix}\right)
\left(\begin{smallmatrix}
1\\
1
\end{smallmatrix}\right)
\Vert^2=\Vert (0,1) \Vert^2=\alpha_{1k}(1)=1/2,
\end{equation}
\begin{equation}
\label{M(T-c)1}
\Vert(T_{\alpha_{1n}}-c_{1n}) {\bf 1}\Vert^2=1-c_{1n}^2.
\end{equation}
Using (\ref{a(nn)'})  we get  (\ref{a(nn)}):
\begin{align}
\label{a(nn)1}
a_{nn}\!=\!1\!-\!t_{kn}^2(2)+t_{kn}^2(2)\alpha_{1k}(1)(1\!-\!c_{1n}^2)\!=\!
1-t_{kn}^2(2)\frac{1+c_{1n}^2}{2}\!=\\
\!1\!-\!\big((1\!+\!c_{2n})/2\big)^2\big(1\!+\!c_{1n}^2)/2.\nonumber
\end{align}
Since $(f_n,g_m^c)=(g_n^c,f_m)=(g_n^c,g_m^c)=0$, for $n\not=m$, we get
$$
a_{nm}=(f_n+g_n^c,f_m+g_m^c)=(f_n,f_m)+(f_n,g_m^c)+(g_n^c,f_m)+(g_n^c,g_m^c)=
(f_n,f_m).
$$
Indeed,
\begin{align*}
 (f_n,g_m^c)=(T_{kn}(1)\otimes (T_{kn}(2)-t_{kn}(2)){\bf 1},
x_{1k}\otimes\big(T_{\alpha_{1m}}-c_{1m}\big)\otimes
t_{km}(2){\bf 1}),\\
= t_{km}(2)(T_{kn}(1){\bf 1},
x_{1k}\otimes\big(T_{\alpha_{1m}}-c_{1m}\big){\bf
1})((T_{kn}(2)-t_{kn}(2)){\bf 1},{\bf 1})\\
=C_1((T_{kn}(2)-t_{kn}(2)){\bf 1},{\bf 1})=0,\,\,\,
(g_n^c,f_m)=C_2((T_{km}(2)-t_{km}(2)){\bf 1},{\bf 1})=0,\\
(g_n^c,g_m^c)=t_{kn}(2)t_{km}(2)\Vert x_{1k}{\bf
1}\Vert^2((T_{\alpha_{1n}}-c_{1n}){\bf 1},
(T_{\alpha_{1m}}-c_{1m}){\bf 1})=0,
\end{align*}
where $C_k$ are some constants. Finally, we get
\begin{equation}
\label{a(nm)(1)}
a_{nm}=(f_n,f_m)=(T_{kn}(1){\bf 1},T_{km}(1){\bf 1})((T_{kn}(2)-t_{kn}(2)){\bf 1},
(T_{km}(2)-t_{km}(2)){\bf 1})
\end{equation}
$$
=(T_{kn}(1){\bf 1},T_{km}(1){\bf 1})\big[(T_{kn}(2){\bf 1},T_{km}(2){\bf 1})
-t_{kn}(2)t_{km}(2)\big].
$$
Using  (\ref{63}) we get
$$
(T_{kn}(1){\bf 1},T_{km}(1){\bf 1})=\big([I+x_{1k}\otimes (T_{\alpha_{1n}}-I)]{\bf 1},
[I+x_{1k}\otimes (T_{\alpha_{1m}}-I)]{\bf 1}\big)=1+
$$
$$
({\bf 1},x_{1k}\otimes (T_{\alpha_{1m}}-I){\bf 1})+(x_{1k}\otimes (T_{\alpha_{1n}}-I){\bf 1},{\bf 1})
+(x_{1k}\otimes (T_{\alpha_{1n}}-I){\bf 1},x_{1k}\otimes (T_{\alpha_{1m}}-I){\bf 1})
$$
$$
=1+\alpha_{1k}(1)[(c_{1m}-1)+(c_{1n}-1)+(c_{1m}-1)(c_{1n}-1)]=
1+\alpha_{1k}(1)(c_{1n}c_{1m}-1)
$$
$$
=\alpha_{1k}(0)+\alpha_{1k}(1)c_{1n}c_{1m}=(1+c_{1n}c_{1m})/2
$$
for $\alpha_{1k}(0)=\alpha_{1k}(1)=1/2$. Finally, we get
$$
a_{nm}=\frac{1+c_{1n}c_{1m}}{2}\Big(\frac{1+c_{2n}c_{2m}}{2}-\frac{1+c_{2n}}{2}
\frac{1+c_{2m}}{2}\Big)=\frac{1+c_{1n}c_{1m}}{8}(1-c_{2n})(1-c_{2m}).
$$
This proves (\ref{a(nm)}). In addition we have
$$
b_n=t_{kn}(2)Ma_n^{(1)}=-(1+c_{2n})(1-c_{1n})/2
,\quad t_{kn}(2)=(1+c_{2n})/2.
$$
We shall estimate $(A^{-1}b,b)$ for an operator $A=(a_{nm})_{n,m}$ defined as follows:
\begin{equation*}
a_{nn}\!=\!1\!-\!\Big(\frac{1+c_{2n}}{2}\Big)^2\frac{1+c_{1n}^2}{2},\quad
a_{nm}\!=\!(1\!+\!c_{1n}c_{1m})(1\!-\!c_{2n})(1\!-\!c_{2m})/8.
\end{equation*}
 We have $A=DD_m(\lambda,c)D$ where
$$
D={\rm diag}(d_n)_{n=1}^m,\quad d_n=(1-c_{2n})(2\sqrt{2})^{-1},
\quad c=(c_n)_n,\quad c_n=c_{1n}.
$$
Finally, we get
$$
(A^{-1}b,b)=(D^{-1}D^{-1}_m(\lambda,c)D^{-1}b,b)=(D^{-1}_m(\lambda,c)D^{-1}b,D^{-1}b)=(D^{-1}_m(\lambda,c)\mu,\mu)
$$
where $\mu=D^{-1}b$. Lemma~\ref{l.(D(l,c)m,m)=} finish the proof.
Since $1+\lambda_n+c_n^2=\frac{a_{nn}}{d_n^2}$  we get
$$
\lambda_n=\frac{a_{nn}-(1+c_n^2)d_n^2}{d_n^2}=
\frac{1-\big(\frac{1+c_{2n}}{2}\big)^2\frac{1+c_{1n}^2}{2}
-\frac{(1+c_{1n}^2)(1-c_{2n})^2}{8}
}{
\frac{(1-c_{2n})^2}{8}}
$$
$$
=\frac{
8-(1+c_{1n}^2)[(1+c_{2n})^2+(1-c_{2n})^2]
}{(1-c_{2n})^2}
$$
$$
=\frac{2[4-(1+c_{1n}^2)(1+c_{2n}^2)]}
{(1-c_{2n})^2}
=\frac{2(1-c_{1n}^2+1-c_{2n}^2+1-c_{1n}^2c_{2n}^2)}{(1-c_{2n})^2}.
$$
We have
$$
\Gamma(g)=\sum_{n=2}^\infty\frac{1}{\lambda_n}=\sum_n\frac{(1-c_{2n})^2}{2(1-c_{1n}^2+1-c_{2n}^2+1-c_{1n}^2c_{2n}^2)}.
$$
Since
$$
b_n=-\frac{(1+c_{2n})(1-c_{1n})}{2},\quad\text{we get}\quad \mu_n=\frac{b_n}{d_n}=-\frac{\sqrt{2}(1+c_{2n})(1-c_{1n})}{(1-c_{2n})}
$$
therefore,
$$
\Gamma(f)=\sum_{n=2}^\infty\frac{\mu_n^2}{\lambda_n}=
\sum_n\frac{(1+c_{2n})^2(1-c_{1n})^2}{2(1-c_{1n}^2+1-c_{2n}^2+1-c_{1n}^2c_{2n}^2)}.
$$
\qed\end{pf}

\begin{lem}
\label{l.x(12)-in-A}
We have $ x_{12}{\bf 1}\in
\langle\left[T_{2n}-c_{2n}\right]{\bf 1}\mid n> k\rangle$ if
$\Sigma_{12}=\infty $ where
\begin{equation}
\label{x(12)-in-A}
\Sigma_{12}=\sum_{n=3}^\infty\frac{(1-c_{1n})^2c_{2n}^2}
{(1-c_{1n})c_{2n}^2+1-c_{2n}}.
\end{equation}
\end{lem}
\begin{pf}
Using (\ref{63}) we get $T_{2n}=T_{2n}(1)\otimes T_{2n}(2)$ so,  $$T_{2n}-c_{2n}=(I+x_{12}(T_{1n}-I))\otimes T_{\alpha_{2n}}-c_{2n}=
x_{12}(T_{1n}-I)\otimes T_{\alpha_{2n}}+T_{\alpha_{2n}}-c_{2n}.
$$
If we chose $t=(t_n)_{n=3}^m$ such that $\sum_{n=3}^mt_nM\xi_n=1$ we get
{\small
$$
\Vert\Big[\sum_{n=3}^mt_n(T_{2n}-c_{2n})\!-\!x_{12}\Big]{\bf 1}\Vert=
\Vert x_{12}\Big(
\sum_{n=3}^mt_n(T_{1n}-I)\otimes T_{\alpha_{2n}}\!-\!I
\Big){\bf 1}+\sum_{n=3}^mt_n(T_{\alpha_{2n}}\!-\!c_{2n}){\bf 1}\Vert^2
$$
}
$$
\Vert\sum_{n=3}^mt_n \big[x_{12}(\xi_n-M\xi_n)+\eta_n\big]
\Vert^2=\sum_{n,k=3}^mt_nt_k(h_n,h_k)=\sum_{n=3}^mt_n^2(h_n,h_n),
$$
since $(h_n,h_k)=0$ for $n\not=k$, where $h_n=x_{12}(\xi_n-M\xi_n)+\eta_n$. Indeed, we show that
\begin{equation}
\label{(h,h)=0}
(h_n,h_k)\!=\!0,\,\,\text{for}\,\,\, n\not=k,\,\,\,\,\,\, (h_n,h_n)\!\sim\!
2(1-c_{1n})c_{2n}^2\!-\!(1-c_{1n})^2c_{2n}^2\!+\!1\!-\!c_{2n}^2.
\end{equation}
We have
$$
(h_n,h_n)=(x_{12}(\xi_n-M\xi_n)+\eta_n,x_{12}(\xi_n-M\xi_n)+\eta_n)=
$$
$$
(x_{12}{\bf 1},{\bf 1})\Vert\xi_n-M\xi_n\Vert^2+\Vert\eta_n\Vert^2+2(x_{12}{\bf 1},{\bf 1})(\xi_n-M\xi_n,\eta_n)=
$$
$$
\alpha_{12}(1)\big[2(1-c_{1n})-(1-c_{1n})^2c_{2n}^2\big]+1-c_{2n}^2+
2\alpha_{12}(1)(c_{1n}-1)(1-c_{2n}^2)=
$$
$$
2\alpha_{12}(1)(1-c_{1n})c_{2n}^2+1-c_{2n}^2-\alpha_{12}(1)(1-c_{1n})^2c_{2n}^2
$$
since
$$
(\xi_n\!-\!M\xi_n,\eta_n)\!=\!(\big[(T_{\alpha_{1n}}\!-\!1)T_{\alpha_{2n}}-(c_{1n}-1)c_{2n}\big]{\bf 1},(T_{\alpha_{2n}}-c_{2n}){\bf 1})=((T_{\alpha_{1n}}\!-\!1){\bf 1},{\bf 1})
$$
$$
\times(T_{\alpha_{2n}}{\bf 1},(T_{\alpha_{2n}}-c_{2n}){\bf 1})-(c_{1n}-1)c_{2n}
((T_{\alpha_{2n}}-c_{2n}){\bf 1},{\bf 1})=(c_{1n}-1)(1-c_{2n}^2).
$$
Finally, we get
$$
\min_{t\in{\mathbb R}^{m-2}}\Big(\sum_{n=3}^mt_n^2(h_n,h_n)\mid\sum_{n=3}^mt_nM\xi_n=1\Big)\to 0
\Leftrightarrow \sum_{n=3}^\infty\frac{\vert M\xi_n\vert^2}{(h_n,h_n)}\sim
\Sigma_{12}=\infty.
$$
We use the following equivalence
$$
\sum_{n=3}^\infty\frac{\vert M\xi_n\vert^2}{(h_n,h_n)}=\sum_{n=3}^\infty\frac{(1-c_{1n})^2c_{2n}^2}
{2\alpha_{12}(1)(1-c_{1n})c_{2n}^2+1-c_{2n}^2-\alpha_{12}(1)(1-c_{1n})^2c_{2n}^2}
$$
$$
\sim
\sum_{n=3}^\infty\frac{(1-c_{1n})^2c_{2n}^2}
{(1-c_{1n})c_{2n}^2+1-c_{2n}}=\Sigma_{12}.
$$
\qed\end{pf}
To find the conditions of the approximation of $x_{2k}$ it is sufficient to interchange $c_{1n}$ and $c_{2n}$ in
Lemma~\ref{l.sigma1k-m=2}. Namely, by analogy we have the following lemma.
\begin{lem}
\label{l.sigma2k-m=2} We have $ x_{2k}{\bf 1}\in
\langle\left[T_{kn}-t_{kn}(1)\right]{\bf 1}\mid n> k\rangle
\Leftrightarrow \Sigma_{2m}\to\infty $ where
\begin{equation}
\label{appr-P2}
\Sigma_{2m}:=
(D^{-1}_m(\lambda,c)\mu,\mu)=
\frac{\Gamma(f_m^{(2)})+\Gamma(f_m^{(2)},g_m^{(2)})+\Gamma(f_m^{(2)},h_m^{(2)})}{1+\Gamma(g_m^{(2)})+\Gamma(h_m^{(2)})+\Gamma(g_m^{(2)},
h_m^{(2)})},
\end{equation}
the vectors $f_m^{(2)},\,g_m^{(2)},\,h_m^{(2)}$ are defined by (\ref{f,g,h}) and
$$
\lambda_n=\frac{2(1-c_{1n}^2+1-c_{2n}^2+1-c_{1n}^2c_{2n}^2)}{(1-c_{1n})^2},\,\,
c_n=c_{2n},\,\, \mu_n=-\frac{\sqrt{2}(1+c_{1n})(1-c_{2n})}{(1-c_{1n})}.
$$
\end{lem}

%
{\small
Finally, to approximate $x_{1k}$ or $x_{2k}$ it is sufficient to have respectively $\Sigma_{1m}\to\infty $
or $\Sigma_{2m}\to\infty $ where
$$
\Sigma_{1m}=
\frac{\Gamma(f_m)+\Gamma(f_m,g_m)+\Gamma(f_m,h_m)}{1+\Gamma(g_m)+\Gamma(h_m)+\Gamma(g_m,h_m)}\sim
$$
$$
\frac{1+\sum_{k=1}^3\Gamma(x_k)+\sum_{1\leq k<r\leq 3}\Gamma(x_k,x_r)+\Gamma(x_1,x_2,x_3)}
{1+\Gamma(x_2)+\Gamma(x_3)+\Gamma(x_2,x_3)}
$$
(where $x_1=f_m,\,\,\,\,x_2=g_m,\,\,\,x_3=h_m$)
$$
=\frac{{\rm det} [I+\gamma(f_m,g_m,h_m)]} {{\rm det}
[I+\gamma(g_m,h_m)]}\sim \frac{{\rm det}
[I+\gamma(F^m,G^m,H^m_1)]} {{\rm det} [I+\gamma(G^m,H^m_1)]},
$$
$$
\Gamma(g)=\sum_n\frac{1}{\lambda_n}=\sum_n\frac{(1-c_{2n})^2}{2(1-c_{1n}^2+1-c_{2n}^2+1-c_{1n}^2c_{2n}^2)},\quad
\Gamma(h)=\sum_n\frac{c_{1n}^2}{\lambda_n}\leq \Gamma(g),
$$
$$
\Gamma(f)=\sum_n\frac{\mu_n^2}{\lambda_n}=
\sum_n\frac{(1+c_{2n})^2(1-c_{1n})^2}{2(1-c_{1n}^2+1-c_{2n}^2+1-c_{1n}^2c_{2n}^2)},
$$
and
$$
\Sigma_{2m}=
\frac{\Gamma(f_m^{(2)})+\Gamma(f_m^{(2)},g_m^{(2)})+\Gamma(f_m^{(2)},h_m^{(2)})}
{1+\Gamma(g_m^{(2)})+\Gamma(h_m^{(2)})+\Gamma(g_m^{(2)},h_m^{(2)})}\sim
$$
$$
\frac{1+\sum_{k=1}^3\Gamma(x_k)+\sum_{1\leq k<r\leq 3}\Gamma(x_k,x_r)+\Gamma(x_1,x_2,x_3)}
{1+\Gamma(x_2)+\Gamma(x_3)+\Gamma(x_2,x_3)}
$$
(where $x_1=f_m^{(2)},\,\,\,\,x_2=g_m^{(2)},\,\,\,x_3=h_m^{(2)}$)
$$
=\frac{{\rm det} [I+\gamma(f_m^{(2)},g_m^{(2)},h_m^{(2)})]}{{\rm
det} [I+\gamma(g_m^{(2)},h_m^{(2)})]} \sim\frac{{\rm det}
[I+\gamma(G^m,F^m,H^m_2)]} {{\rm det} [I+\gamma(F^m,H^m_2)]},
$$
$$
\Gamma(g^{(2)})\!=\!\sum_n\frac{1}{\lambda_n}=\sum_n\frac{(1-c_{1n})^2}
{2(1\!-\!c_{1n}^2\!+\!1\!-\!c_{2n}^2\!+\!1\!-\!c_{1n}^2c_{2n}^2)},
\,\,
\Gamma(h^{(2)})\!=\!\sum_n\frac{c_{2n}^2}{\lambda_n}\leq \Gamma(g^{(2)}),
$$
$$
\Gamma(f^{(2)})=\sum_n\frac{\mu_n^2}{\lambda_n}=
\sum_n\frac{(1+c_{1n})^2(1-c_{2n})^2}{2(1-c_{1n}^2+1-c_{2n}^2+1-c_{1n}^2c_{2n}^2)}.
$$
Since $1<1+c_{rn}\leq 2$ finally, we get
\begin{equation}
\label{appr-(F,G)1}
\Sigma_{1m}\simeq
\frac{\Gamma(F^m)+\Gamma(F^m,G^m)+\Gamma(F^m,H^m_1)}{1+\Gamma(G^m)+\Gamma(H^m_1)+\Gamma(G^m,H^m_1)},
\end{equation}
\begin{equation}
\label{appr-(F,G)2} \Sigma_{2m}\simeq
\frac{\Gamma(G^m)+\Gamma(G^m,F^m)+\Gamma(G^m,H^m_2)}{1+\Gamma(F^m)+\Gamma(H^m_2)+\Gamma(F^m,H^m_2)},
\end{equation}
 where we denote
\begin{equation}
\label{F,G=}
F_n\!=\!\frac{1-c_{1n}}{(1-c_{1n}+1-c_{2n}+1-c_{1n}c_{2n})^{1/2}},\,\,
G_n\!=\!\frac{1-c_{2n}}{(1-c_{1n}+1-c_{2n}+1-c_{1n}c_{2n})^{1/2}},
\end{equation}
$$
F=(F_n)_n,\,\,G=(G_n)_n,\,\, H_1=(H_n^1)_n,\quad H_2=(H_n^2)_n,\,\, H^1_n=G_nc_{1n},\,\, H^2_n=F_nc_{2n},
$$
\begin{equation}
\label{F^m,G^m=}
F^m=(F_n)_{n=2}^m,\quad G^m=(G_n)_{n=2}^m,\quad H^m_1=(H^1_n)_{n=2}^m,\quad H^m_2=(H^2_n)_{n=2}^m.
\end{equation}
}
\begin{lem}
\label{l.1k+2k=!}
If $S_{11}^L(\mu_\alpha)+S_{22}^L(\mu_\alpha)=\infty$ then
\begin{equation}
\label{f+g=infty}
\Gamma(f)+\Gamma(g)=\infty,\quad\Gamma(f^{(2)})+\Gamma(g^{(2)})=\infty
\quad\text{and}\quad \Gamma(F)+\Gamma(G)=\infty.
\end{equation}
\end{lem}
\begin{pf}
Since $c_{kn}=2\sqrt{\alpha_{kn}(0)\alpha_{kn}(1)}\in (0,1]$ we conclude that
\begin{equation}
\label{sigma1}
\Gamma(f)\sim \Gamma(g^{(2)})\sim \Gamma(F)=\sum_{n=3}^\infty\frac{(1-c_{1n})^2}{(1-c_{1n}+1-c_{2n}+1-c_{1n}c_{2n})}
\end{equation}
$$
\sim
\sum_{n=3}^\infty\frac{(1-c_{1n})^2(1+c_{2n})^2}{(1-c_{1n}+1-c_{2n}+1-c_{1n}c_{2n})},
$$
\begin{equation}
\label{sigma2}
\Gamma(g)\sim \Gamma(f^{(2)})\sim \Gamma(G)=\sum_{n=3}^\infty\frac{(1-c_{2n})^2}{(1-c_{1n}+1-c_{2n}+1-c_{1n}c_{2n})}
\end{equation}
$$
\sim
\sum_{n=3}^\infty\frac{(1+c_{1n})^2(1-c_{2n})^2}{(1-c_{1n}+1-c_{2n}+1-c_{1n}c_{2n})}.
$$
If $\Gamma(f)+\Gamma(g)<\infty$ or $\Gamma(f^{(2)})+\Gamma(g^{(2)})<\infty$ we get $\Gamma(F)+\Gamma(G)<\infty$ and
$$
\sum_{n=3}^\infty\frac{(1+c_{1n})^2(1-c_{2n})^2}{(1\!-\!c_{1n}\!+\!1\!-\!c_{2n}+1\!-\!c_{1n}c_{2n})}+
\sum_{n=3}^\infty\frac{(1\!-\!c_{1n})^2(1+c_{2n})^2}{(1\!-\!c_{1n}\!+\!1\!-\!c_{2n}\!+\!1\!-\!c_{1n}c_{2n})}<\infty
$$
therefore, $\Sigma<\infty$ where
$$
\Sigma:=\sum_{n=3}^\infty\frac{[
(1+c_{1n})(1-c_{2n})+(1-c_{1n})(1+c_{2n})]^2}
{4(1-c_{1n}+1-c_{2n}+1-c_{1n}c_{2n})}=
$$
$$
\sum_{n=3}^\infty\frac{(1-c_{1n}c_{2n})^2}
{(1-c_{1n}+1-c_{2n}+1-c_{1n}c_{2n})}.
$$
Finally, $\Gamma(F)+\Gamma(G)+\Sigma<\infty$ hence,
$$
\infty>\sum_n\frac{(1-c_{1n}+1-c_{2n}+1-c_{1n}c_{2n})^2}{(1-c_{1n}+1-c_{2n}+1-c_{1n}c_{2n})}=
\sum_n(1-c_{1n}+1-c_{2n}+1-c_{1n}c_{2n})
$$
$$
>S_{11}^L(\mu_\alpha)+S_{22}^L(\mu_\alpha)=\infty.
$$
This contradiction proves (\ref{f+g=infty}).
\qed\end{pf}
\begin{rem}
\label{op-in-W-alg} We have proved the convergence
$\sum_{n=N_1}^{N_2}t_n[T_{kn}-t_{kn}(r)]\to x_{rk}$ for $r=1,2$ only on
the vector $f={\bf 1}$.The same argument holds for the total set
of vectors of the form $f=\otimes_{k=1}^nf_k\otimes {\bf
1}\otimes{\bf 1}\cdots$ in the space $L^2(X^2,\mu^2)=H_{12}
\otimes_{k=3}^\infty\big(H_{1k}\otimes H_{2k}\big)$. Hence,
$x_{rk}\in {\mathfrak A}^2$. In what follows we will use the
same arguments.
\end{rem}
  It is useful
to use the analogy with the case of the field ${\bf k}={\mathbb
R}$.
\begin{rem}
\label{r.combi-R1} In the case of the field ${\bf k}={\mathbb R}$
the generators $A_{kn}$ and $A_{2n}$ of the corresponding
one-parameter groups has the following form \cite{Kos04}:
$$
A_{kn}=x_{1k}D_{1n}+x_{2k}D_{2n},\quad A_{2n}=x_{12}D_{1n}+D_{2n}.
$$
If we are able to approximate the variable $x_{2n}$,  but not the
$x_{1n}$. It is reasonable to use the following expressions
\begin{equation}
\label{combi-R1}
 A_{kn}-x_{2k}A_{2n}=(x_{1k}-x_{12}x_{2k})D_{1n}
\end{equation}
in order to approximate first the expression $x_{1k}-x_{12}x_{2k}$ by linear
combinations  $\sum_nt_n(x_{1k}-x_{12}x_{2k})D_{1n}^2$. Further, we
can approximate  the variable $x_{12}$ by
$\sum_nt_k(x_{1k}-x_{12}x_{2k})$. After we can approximate the
variable $x_{1k}$ (see also the details in {\rm \cite{Kos02.3}}).
\end{rem}

Let $x_{2k}={\rm diag}(0,1)\in {\mathfrak A}^2$. We {\it try to guess} an
analogue of the expression $A_{kn}-x_{2k}A_{2n}$. The analogue of $A_{kn}$ is
$T_{kn}-I=C(T_{kn})-p$, by  Remark~\ref{r.gen:C(T)-p}. So, the
analogue of $A_{kn}-x_{2k}A_{2n}$ is
$T_{kn}-I-x_{2k}\otimes(T_{2n}-I).$ We shall  use the following
combinations:
\begin{equation}
\label{}
 T_{kn}-I-x_{2k}\otimes(T_{2n}-I).
\end{equation}
For $k<n$ set
\begin{equation}
\label{x=P,t=(T-I)}
\tau_{kn}=(T_{\alpha_{kn}}-I).
\end{equation}
Using (\ref{T2n(1)-1.0}) we get
$$
T_{kn}=(I+x_{1k}\tau_{1n})(I+x_{2k}\tau_{2n})=I+x_{1k}\tau_{1n}+x_{2k}\tau_{2n}+x_{1k}x_{2k}\tau_{1n}\tau_{2n},
$$
$$
T_{2n}=(I+x_{12}\tau_{1n})(I+\tau_{2n})=I+x_{12}\tau_{1n}+\tau_{2n}+x_{12}\tau_{1n}\tau_{2n},
$$
therefore, we get
\begin{equation}
\label{x-x.x(1)}
 T_{kn}-I-x_{2k}\otimes(T_{2n}-I)=(x_{1k}-x_{12}x_{2k})\tau_{1n}+(x_{1k}x_{2k}-x_{12}x_{2k})\tau_{1n}\tau_{2n}.
\end{equation}
\begin{lem}
\label{x(1k)-x(12)x(2k)} We have
$$
(x_{1k}-x_{12}x_{2k}){\bf 1}\in \langle\left[T_{kn}-I-x_{2k}(T_{2n}-I)\right]{\bf 1}\mid n> k\rangle
$$
if and only if $\Delta(f'_{m},g'_{m})
\to\infty$ where
\begin{equation}
\label{1k-12.2k}
\Delta(f_m',g_m')=\frac{\Gamma(f_m')+\Gamma(f'_m,g'_m)}{1+\Gamma(g'_m)},\quad
\end{equation}
and $f_m',\,\,g_m'\in{\mathbb R}^{m-2}$ are defined as follows:
\begin{equation}
\label{f_m',g_m'}
f_m'=(\sqrt{1-c_{1n}})_{n=3}^m,\quad g_m'=(\sqrt{1-c_{1n}}(1-c_{2n}))_{n=3}^m.
\end{equation}
\end{lem}
\begin{pf}
Set $b_n=M(T_{\alpha_{1n}}-I){\bf 1}=M\tau_{1n}{\bf 1}=c_{1n}-1$,
$$
f_n=[(x_{1k}-x_{12}x_{2k})\tau_{1n}+(x_{1k}x_{2k}-x_{12}x_{2k})\tau_{1n}\tau_{2n}]{\bf 1},
$$
$$
f_n^c=[(x_{1k}-x_{12}x_{2k})(\tau_{1n}-M\tau_{1n}{\bf 1})+(x_{1k}x_{2k}-x_{12}x_{2k})\tau_{1n}\tau_{2n}]{\bf 1}
$$
$$
=[(x_{1k}-x_{12}x_{2k})(T_{\alpha_{1n}}-c_{1n})-(x_{1k}x_{2k}-x_{12}x_{2k})\otimes(T_{\alpha_{1n}}-I)\otimes(T_{\alpha_{2n}}-I)]{\bf 1}.
$$
Take $t=(t_n)_{n=3}^N$ such that $\sum_{n=3}^Nt_nb_n=1$, then we get
$$
\Vert\sum_{n=3}^Nt_nf_n-(x_{1k}-x_{12}x_{2k}){\bf 1}\Vert^2=\Vert\sum_{n=3}^Nt_nf_n^c\Vert^2=
\sum_{n=3}^Nt_nt_m(f_n^c,f_m^c)=(At,t),
$$
where
$$
A=(a_{nm})_{n,m=3}^N,\quad a_{nm}=(f_n^c,f_m^c)_{n,m=3}^N.
$$
As before, we use the following estimation
for a positively definite
operator $A$ acting in a space $H$ and a vector $b\in H$:
$$
\min_{t\in H}\Big((At,t)\mid (t,b)=1\Big)=(A^{-1}b,b)^{-1}.
$$
We show that
\begin{equation}
\label{a_{nn}2}
a_{nn}=(2-c_{2n})(1-c_{1n}^2)/2+(1-c_{1n})(1-c_{2n}),
\end{equation}
\begin{equation}
\label{a_{nm}2}
a_{nm}=d_nd_m,\,\,n\not=m,\quad\text{where}\quad d_n=(1-c_{1n})(1-c_{2n})/2.
\end{equation}
Set
$$
h=(x_{1k}-x_{12}x_{2k}),\quad g=(x_{1k}x_{2k}-x_{12}x_{2k}),\,\,
$$
$$
\xi_n=(T_{\alpha_{1n}}-c_{1n}){\bf 1},\quad
\eta_n=(T_{\alpha_{1n}}-I)\otimes(T_{\alpha_{2n}}-I){\bf 1}.
$$
We can suppose that $\mu_{\alpha_{1k}}=\mu_{\alpha_{2k}}$ are invariant measures for  all $k\leq n_0$ and  fixed $n_0$. This does not change
the equivalence class of the measure $\mu$ hence, the equivalence class of the representation.
Then using (\ref{Mx_{kn}1}) we get
$$
a_{nn}=(f_n^c,f_n^c)=\Vert h\Vert^2\Vert\xi_n\Vert^2+\Vert g\Vert^2\Vert\eta_n\Vert^2+
2(g,h)(\xi_n,\eta_n),
$$
$$
\Vert h\Vert^2= \Vert (x_{1k}-x_{12}x_{2k}){\bf 1}\Vert^2=
((x_{1k}^2-2x_{12}x_{1k}x_{2k}+x_{12}^2x_{2k}^2){\bf 1},{\bf 1})=
1/2,
$$
$$
\Vert g\Vert^2=\Vert (x_{1k}x_{2k}-x_{12}x_{2k}){\bf 1}\Vert^2=
(x_{1k}^2x_{2k}^2-2x_{12}x_{1k}x_{2k}^2+x_{12}^2x_{2k}^2){\bf 1},{\bf 1})=
1/4,
$$
$$
(g,h)=((x_{1k}-x_{12}x_{2k}){\bf 1},(x_{1k}x_{2k}-x_{12}x_{2k}){\bf 1})
$$
$$
=((x_{1k}^2x_{2k}-x_{12}x_{1k}x_{2k}-x_{12}x_{1k}x_{2k}^2+x_{12}^2x_{2k}^2){\bf 1},{\bf 1})=
1/4,
$$
$$
\Vert\xi_n\Vert^2=1-c_{1n}^2,\quad \Vert\eta_n\Vert^2=4(1-c_{1n})(1-c_{2n}),\quad (\xi_n,\eta_n)=(1-c_{1n}^2)(1-c_{2n}).
$$
Finally, we get (\ref{a_{nn}2}). Indeed, we have
$$
a_{nn}=(1-c_{1n}^2)/2+(1-c_{1n})(1-c_{2n})+(1-c_{1n}^2)(1-c_{2n})/2
$$
$$
=(2-c_{2n})(1-c_{1n}^2)/2+(1-c_{1n})(1-c_{2n}).
$$
We show that $(\xi_n,\eta_n)=-(1-c_{1n}^2)(1-c_{2n})$. Indeed, we get
$$
(\xi_n,\eta_n)=((T_{\alpha_{1n}}-c_{1n}){\bf 1},
(T_{\alpha_{1n}}-I)\otimes(T_{\alpha_{2n}}-I)]{\bf 1})=
((T_{\alpha_{1n}}-c_{1n}){\bf 1},(T_{\alpha_{1n}}-I){\bf 1})\times
$$
$$
((T_{\alpha_{2n}}-I){\bf 1},{\bf 1})=\Vert (T_{\alpha_{1n}}-c_{1n}){\bf 1}\Vert^2((T_{\alpha_{2n}}-I){\bf 1},{\bf 1})=
-(1-c_{1n}^2)(1-c_{2n}).
$$
Further, since $(\xi_n,\xi_m)=(\xi_n,\eta_m)=0$  for $n\not=m$ we get
$$
a_{nm}=(f_n^c,f_m^c)=(h\xi_n+g\eta_n,h\xi_m+g\eta_m)=(g,g)(\eta_n,\eta_m)=d_nd_m/4,
$$
where $d_n=\frac{(1-c_{1n})(1-c_{2n})}{2}$. This proves (\ref{a_{nm}2}).
We use the fact that $(g,g)=1/4$  and
$$
(\eta_n,\eta_m)=((T_{\alpha_{1n}}-I)\otimes(T_{\alpha_{2n}}-I){\bf 1},
(T_{\alpha_{1m}}-I)\otimes(T_{\alpha_{2m}}-I){\bf 1})
$$
$$
=((T_{\alpha_{1n}}-I)\otimes
(T_{\alpha_{2n}}-I){\bf 1},{\bf 1})
(T_{\alpha_{1m}}-I)(T_{\alpha_{2m}}-I){\bf 1},{\bf 1})
$$
$$
=(1-c_{1n})(1-c_{2n})(1-c_{1m})(1-c_{2m}).
$$
Since $a_{nm}$ is a product $a_{nm}=d_nd_m/4$ we can use  the particular case of Lemma~\ref{l.(D(l,c)m,m)=} for $c=0$.
to calculate $(A^{-1}b,b)$. We have $A=DD_m(\lambda)D$ where $D={\rm diag}(D_n)_{n=1}^m$ and
$D_m(\lambda)$ is defined by (\ref{D(lambda,c)}).
Finally, we get if we set $D^{-1}b=\mu$
$$
(A^{-1}b,b)=
(D^{-1}_m(\lambda,c)\mu,\mu)\stackrel{(\ref{D(l,c)=G()})}{=}
\frac{\Gamma(f_m)+\Gamma(f_m,g_m)}{1+\Gamma(g_m)}=\Delta(f_m,g_m),
$$
where $f_m=(\mu_k/\sqrt{\lambda_k})_{k=1}^m,\,\,g_m=\Big(1/\sqrt{\lambda_k})_{k=1}^m$ (see (\ref{f,g,h})).
To calculate $\Delta(f_m,g_m)$
we have $\lambda_n=\frac{a_{nn}}{d_n^2}-1=\frac{a_{nn}-d_n^2}{d_n^2}$ and $\mu_n=\frac{b_n}{d_n}=-\frac{2}{(1-c_{2n})}$ therefore,
$$
\lambda_n=\frac{(2-c_{2n})(1-c_{1n}^2)/2+(1-c_{1n})(1-c_{2n})-(1-c_{1n})^2(1-c_{2n})^2/4}{(1-c_{1n})^2(1-c_{2n})^2/4}
$$
$$
=\frac{2(2-c_{2n})(1+c_{1n})+4(1-c_{2n})-(1-c_{1n})(1-c_{2n})^2}
{(1-c_{1n})(1-c_{2n})^2}
$$
$$
=\frac{2(1+c_{1n})(2-c_{2n})+(1-c_{2n})[4-(1-c_{1n})(1-c_{2n})]
}{(1-c_{1n})(1-c_{2n})^2}.
$$
Finally, we get $\Gamma(g)=\sum_{n=3}^\infty  \frac{1}{\lambda_n}$
and $\Gamma(f)=\sum_{n=3}^\infty \frac{\mu_n^2}{\lambda_n}$ or
$$
\Gamma(g)=
\sum_{n=3}^\infty \frac{(1-c_{1n})(1-c_{2n})^2}{f(x_n,y_n)},\quad
\Gamma(f)=
4\sum_{n=3}^\infty \frac{(1-c_{1n})}{f(x_n,y_n)},
$$
where $f(x_n,y_n)=2(2+x_n)(1+y_n)+y_n(4-x_ny_n)$ and $x_n=1-c_{1n},$ \\$y_n=1-c_{2n}$.

Since $4\leq f(x_n,y_n)\leq 16$ we conclude that $\Gamma(f)\sim \Gamma(f')$ and $\Gamma(g)\sim
\Gamma(g')$ where  $f'$ and $g'$ are defined as follows (see (\ref{f_m',g_m'}))
$$
f'=(\sqrt{1-c_{1n}})_{n=3}^\infty,\quad g'=(\sqrt{1-c_{1n}}(1-c_{2n}))_{n=3}^\infty.
$$
$$\text{Hence,}\,\,\,\Gamma(f')=\sum_{n=3}^\infty(1-c_{1n})=S_{11}^L(\mu)=\infty,\quad
\Gamma(g')=\sum_{n=3}^\infty(1-c_{1n})(1-c_{2n})^2.
$$
$$
\text{Finally, we get}\quad\lim_m\Delta(f_m,g_m)\sim \lim_m\Delta(f_m',g_m')=
\frac{\Gamma(f')+\Gamma(f',g')}{1+\Gamma(g')}.
$$
\qed\end{pf}


{\bf The schema of the proof of irreducibility for $m=2$}.
Recall some notations (see (\ref{appr-(F,G)1})--(\ref{F^m,G^m=}))
$$
\Sigma_{12}=\sum_{n=2}^\infty\frac{(1-c_{1n})^2c_{2n}^2}{(1-c_{1n})c_{2n}^2+1-c_{2n}},\quad
d_n=1-c_{1n}+1-c_{2n}+1-c_{1n}c_{2n},
$$
$$
\Gamma(F^m)= \Vert
F^m\Vert^2=\!\sum_{n=3}^m\frac{(1-c_{1n})^2}{d_n},\quad
\Gamma(G^m)= \Vert
G^m\Vert^2=\!\sum_{n=3}^m\frac{(1-c_{2n})^2}{d_n},
$$
$$
\Gamma(H_1^m)=\Vert
H_1^m\Vert^2=\sum_{n=3}^m\frac{(1-c_{2n})^2c_{1n}^2}{d_n},\,\,
\Gamma(H_2^m)=\Vert
H_2^m\Vert^2=\sum_{n=3}^m\frac{(1-c_{1n})^2c_{2n}^2}{d_n},
$$
\begin{equation*}
\Sigma_{1m}\simeq
\frac{\Gamma(F^m)+\Gamma(F^m,G^m)+\Gamma(F^m,H^m_1)}{1+\Gamma(G^m)+\Gamma(H^m_1)+\Gamma(G^m,H^m_1)},
\end{equation*}
\begin{equation*}
\Sigma_{2m}\simeq
\frac{\Gamma(G^m)+\Gamma(G^m,F^m)+\Gamma(G^m,H^m_2)}{1+\Gamma(F^m)+\Gamma(H^m_2)+\Gamma(F^m,H^m_2)},
\end{equation*}
$$
\Delta(f,g)=\frac{\Gamma(f)+\Gamma(f,g)}{1+\Gamma(g)},\quad
f=(f_n)_{n\in{\mathbb N}},\,\,g=(g_n)_{n\in{\mathbb N}},
$$
$$
\Gamma(f')=\sum_{n=3}^\infty(1-c_{1n}),\quad
\Gamma(g')=\sum_{n=3}^\infty(1-c_{1n})(1-c_{2n})^2,
$$
$$
S_{11}^L(\mu)\!=\!\sum_{n=2}^\infty(1-c_{1n}),\,\, S_{22}^L(\mu)\!=\!\sum_{n=3}^\infty(1-c_{2n}),\,\, S_{12}^L(\mu)\!=\!\sum_{n=3}^\infty\alpha_{2n}(1)(1-c_{1n}).
$$
To prove the irreducibility,  consider different cases.

{\bf Case 1.} Let $\Sigma_{12}=\infty$ then by
Lemma~\ref{l.x(12)-in-A} we conclude that $x_{12}\in
{\mathfrak A}^2$. By (\ref{Reg=otimes(r)})) and (\ref{T2n(1)-1.0}) we get
\begin{equation}
\label{T_2n=} T_{2n}\!=\!T_{2n}(1)\otimes
T_{2n}(2),\,\,\,T_{2n}(1)\!=\!I+x_{12}\otimes (T_{\alpha_{1n}}-I),\,\,
T_{1n}=T_{\alpha_{1n}},\,\,T_{2n}(2)\!=\! T_{\alpha_{2n}},
\end{equation}
\begin{equation}
\label{T_kn=}
T_{kn}\!=\!T_{kn}(1)\otimes T_{kn}(2),\,\,
T_{kn}(1)\!=\!I+x_{1k}\otimes (T_{\alpha_{1n}}-I),\,\,
T_{kn}(2)\!=\!I+x_{2k}\otimes (T_{\alpha_{2n}}-I).
\end{equation}
%
%
\begin{rem}
\label{r.x(kn)-in-A}
We note that

(a) if $x_{12}\in {\mathfrak A}^2$ then $T_{\alpha_{2n}}\in {\mathfrak A}^2$ for $n\geq 3$;

(b) if $x_{1k}\in {\mathfrak A}^2$ then $x_{2k}\in {\mathfrak A}^2$ for $k\geq 3$;

(c) if $x_{12},\,x_{2k}\in {\mathfrak A}^2$ for $n\geq 3$ then $x_{1k}\in {\mathfrak A}^2$ for $k\geq 3$.

To prove (a) it is sufficient to use (\ref{T_2n=}) and $T_{1n}=T_{\alpha_{1n}}\in {\mathfrak A}^2$. Since $x_{12},\,\,T_{1n}\in {\mathfrak A}^2$ then $T_{2n}^{-1}(1)=T_{2n}(1)\in {\mathfrak A}^2$ therefore,  $T_{\alpha_{2n}}=T_{2n}(2)=T_{2n}^{-1}(1)T_{2n}\in {\mathfrak A}^2$.

(b) Since $T_{1k},\,T_{kn}\in {\mathfrak A}^2$, $T_{1n}=T_{\alpha_{1n}}$ using (\ref{T_kn}) we conclude that  $T_{kn}(1)^{-1}=T_{kn}(1) \in {\mathfrak A}^2$
therefore,  $T_{kn}(1)^{-1}T_{kn}=T_{kn}(2)\in {\mathfrak A}^2$
and finally, $x_{2k}\in {\mathfrak A}^2$ by Lemma~\ref{x(2k)-in-A} that is an analogue of Lemma~\ref{l.P2}. For an arbitrary $p$ we can use the relation
$T_{kn}(1)^{-1}=T_{kn}(1)^{p-1}\in{\mathfrak A}^2$.

(c) By (a) we get  $T_{\alpha_{2n}}\in {\mathfrak A}^2$. Since $x_{2k}\in {\mathfrak A}^2$ we conclute that $T_{kn}(2)\in {\mathfrak A}^2$ therefore, $T_{kn}(1)\in {\mathfrak A}^2$ hence, $x_{1k}\in {\mathfrak A}^2$ by Lemma~\ref{l.P2}.
\end{rem}
\begin{lem}
\label{x(2k)-in-A} We have $x_{2k}{\bf 1}\!\in\!
\langle(T_{kn}(2)\!-\!I){\bf 1}\mid n\!>\!k\rangle $ if and only if
$S^L_{22}(\mu)\!=\!\infty.$
\end{lem}

\begin{rem}
\label{r.x(12),(A)'=}
Since $x_{12},T_{1n},T_{2n+1}\in{\mathfrak A}^2,\,\,n\geq 2$ we
conclude by Remark~\ref{r.x(kn)-in-A} (a) that
$T_{\alpha_{2n}}\in{\mathfrak
A}^2,n>2$. Finally, we get
$x_{12},\,\,T_{\alpha_{1n}},T_{\alpha_{2n}}\in{\mathfrak A}^2$ for
all $n\geq 3$. This family of operators is commuting and has
common simple spectrum \cite{Ber86} therefore, the von Neumann
algebra
$$
L^\infty(x_{12},T_1,T_2):=L^\infty\left(\begin{smallmatrix}
x_{12}&T_{\alpha_{13}}&\cdots&T_{\alpha_{1n}}&\cdots\\
                          &T_{\alpha_{23}}&\cdots&T_{\alpha_{2n}}&\cdots\\
\end{smallmatrix}\right)\ni
 f\left(\begin{smallmatrix}
 x_{12}&T_{\alpha_{13}}&\cdots&T_{\alpha_{1n}}&\cdots\\
       &T_{\alpha_{23}}&\cdots&T_{\alpha_{2n}}&\cdots\\
\end{smallmatrix}\right),
$$
 generated by this family is maximal abelian subalgebra in
${\mathfrak A}^2$ and consists of all $L^\infty$ functions, i.e.,
bounded  operator-valued functions depending on the variables
$(x_{12},T_{\alpha_{1n}},T_{\alpha_{2n}},\,\,n\geq 3)$.  Therefore, $({\mathfrak
A}^2)'\subset L^\infty(x_{12},T_1,T_2)'= L^\infty(x_{12},T_1,T_2)$
hence, any operator $f=f(x_{12},T_{\alpha_{1n}},T_{\alpha_{2n}},\,\,n\geq 3)$ from $({\mathfrak A}^2)'$  belongs to
$L^\infty(x_{12},T_1,T_2)$. Since $T_{12}\in {\mathfrak A}^2$ the
relation $[f,T_{12}]=0$ implies that the operator $f$ does not
depend on $x_{12}$. Finally, any operator from $({\mathfrak
A}^2)'$  is a function $f$ in commuting family
$(T_{\alpha_{1n}},T_{\alpha_{2n}},n\geq 3)$.
\begin{equation}
\label{x(12),(A)'=}
({\mathfrak A}^2)'\subset L^\infty\left(\begin{smallmatrix}
T_{\alpha_{13}}&\cdots&T_{\alpha_{1n}}&\cdots\\
T_{\alpha_{23}}&\cdots&T_{\alpha_{2n}}&\cdots\\
\end{smallmatrix}\right).
\end{equation}
\end{rem}
The commutation
$[f,T_{kk+1}]=0$ for all $k\geq 1$ implies, by
Lemma~\ref{l.Delta_1} and its analogue, that $f$ depends only on
the following expressions:
$$
\Delta^{1,r}_{s,\infty}\!=\!T_{\alpha_{1s}}^r\!\prod_{k=s+1}^\infty\!
p^{-1}C(T_{\alpha_{1k}}),\,\,
\Delta^{2,r}_{s,\infty}\!=\!T_{\alpha_{2s}}^r\!\prod_{k=s+1}^\infty\!
p^{-1}C(T_{\alpha_{2k}}),\,\,r\in{\mathbb F}_p\setminus\{0\},\,s\geq 3.
$$
But the latter expressions are well defined if and only if
$S^L_{11}(\mu)<\infty$ and $S^L_{22}(\mu)<\infty$ by
Lemma~\ref{4cond.} and its analogue.

{\bf Case 2.} Let $\Sigma_{12}<\infty$ and $\Gamma(G)<\infty$ then $\Gamma(F)=\infty$. This conditions are incompatible.
Indeed,  since $d_n=1-c_{1n}+1-c_{2n}+1-c_{1n}c_{2n}<3$ we get
\begin{equation}
\label{Gam(G)<infty}
\infty>\Gamma(G)=\sum_{n=3}^\infty\frac{(1-c_{2n})^2}{d_n}>\frac{1}{3}
\sum_{n=3}^\infty(1-c_{2n})^2.
\end{equation}
Therefore, $\lim_nc_{2n}=1$ hence, $\infty>\Sigma_{12}\sim \Vert F \Vert^2=\infty$, contradiction. Indeed, for an arbitrary $\varepsilon>0$ we have for sufficiently big $N$
$$
\infty>\Sigma_{12}=\sum_{n=3}^\infty\frac{(1-c_{1n})^2c_{2n}^2}
{(1-c_{1n})c_{2n}^2+1-c_{2n}}\geq
(1-\varepsilon)^2
\sum_{n=N}^\infty
\frac{(1-c_{1n})^2}
{(1-c_{1n})c_{2n}^2+1-c_{2n}}>
$$
$$
(1-\varepsilon)^2
\sum_{n=N}^\infty
\frac{(1-c_{1n})^2}
{1-c_{1n}+1-c_{2n}}
\sim \Vert F \Vert^2=\infty.
$$

{\bf Case 3.}  Let $\Sigma_{12}<\infty$ and $\Gamma(F)<\infty$,
then $\Gamma(H_2)<\infty$ and $\Gamma(G)=\infty$, by Lemma~\ref{l.1k+2k=!}. Therefore, $\Sigma_{2m}\to\infty$ hence,
$x_{2k}\in {\mathfrak A}^2$ for $k\geq 3$.
Further, using notations
$\tau_{1n}=T_{\alpha_{1n}}-I$ and $\tau_{2n}=T_{\alpha_{2n}}-I$  (see (\ref{x=P,t=(T-I)}))
 we get
\begin{equation}
\label{x-x.x}
T_{kn}-I-x_{2k}(T_{2n}-I)=(x_{1k}-x_{12}x_{2k})\tau_{1n}+(x_{1k}x_{2k}-x_{12}x_{2k})\tau_{1n}\tau_{2n}
\end{equation}
(see (\ref{x-x.x(1)})).  By Lemma~\ref{x(1k)-x(12)x(2k)}, we conclude that
$x_{1k}-x_{12}x_{2k}\in {\mathfrak A}^2$ if
$\Delta(f_m',g_m')$ $\to\infty$ where
\begin{equation}
\label{Delta(fm,gm)}
\Delta(f_m',g_m')=\frac{\Gamma(f_m')+\Gamma(f_m',g_m')}{1+\Gamma(g_m')},
\end{equation}
 \begin{equation}
\label{Gamma(fm,gm)}
\Gamma(f'_m)=\sum_{n=3}^m(1-c_{1n}),\quad  \Gamma(g'_m)=\sum_{n=3}^m(1-c_{1n})(1-c_{2n})^2.
\end{equation}
 The case (3) splits into two cases
(30), when $\Delta(f_m',g_m')\to\infty$ and case (31), when $\Vert f'-tg'\Vert<\infty$ for some
$t\in{\mathbb R}\setminus\{0\}$.

In the case {\bf 30}, by Lemma~\ref{x(1k)-x(12)x(2k)}, we conclude that
$$
x_{1k}-x_{12}x_{2k}\in {\mathfrak A}^2\quad\text{if}\quad \Delta(f_m',g_m')\to\infty.
$$

We use Lemmas~\ref{l.min=proj0} and \ref{l.min=proj} proved in \cite{Kos16Ar}:
\begin{lem}
\label{l.min=proj0}
Let $f=(f_k)_{k\in{\mathbb N}}$ and $g=(g_k)_{k\in{\mathbb N}}$ be two real vectors  such that
$\Vert f\Vert^2=\infty$  where $\Vert f\Vert^2=\sum_kf_k^2$.  Denote by $f_{(n)}$,
$g_{(n)}\in {\mathbb R}^n$ their projections to the subspace
${\mathbb R}^n$, i.e., $f_{(n)}=(f_k)_{k=1}^n,\quad  g_{(n)}=(g_k)_{k=1}^n$ and set
\begin{equation}
\label{Delta-to-infty}
\Delta(f_{(n)},g_{(n)})=\frac{\Gamma(f_{(n)})+\Gamma(f_{(n)},g_{(n)})}{\Gamma(g_{(n)})+1}
\quad\text{then}\quad \lim_{n\to\infty}\Delta(f_{(n)},g_{(n)})=\infty
\end{equation}
in the following cases:
\begin{eqnarray*}
(a)&\Vert g\Vert^2<\infty,\\
(b)&\Vert g\Vert^2=\infty,\quad\text{and}\quad  \lim_{n\to\infty}\frac{\Vert f_{(n)}\Vert}{\Vert g_{(n)}\Vert}=\infty,\\
(c)&\Vert f\Vert^2=\Vert g\Vert^2=\Vert f+s g\Vert^2=\infty,\quad\text{for all}\quad
     s\in {\mathbb R}\setminus\{0\}.
\end{eqnarray*}
\end{lem}
\begin{pf} Obviously $\lim_{n\to\infty}\Delta(f_{(n)},g_{(n)})=\infty$ if conditions (a)  or (b) hold.
The implication $ (c)\Rightarrow (\ref{Delta-to-infty})$ is based on the following lemma.
\qed\end{pf}
\begin{lem}
\label{l.min=proj}
Let $f=(f_k)_{k\in{\mathbb N}}$ and $g=(g_k)_{k\in{\mathbb N}}$ be two real vectors  such that
\begin{equation}
\label{norm=infty}
\Vert f\Vert^2=\Vert g\Vert^2=\Vert C_1 f+C_2 g\Vert^2=\infty,\quad\text{for all}\quad
(C_1,C_2)\in {\mathbb R}^2\setminus\{0\},
\end{equation}
\begin{equation}
\label{final.}
\text{then}\quad
\lim_{n\to\infty}\frac{\Gamma(f_{(n)},g_{(n)})}{\Gamma(g_{(n)})}=\infty\quad\text{and}\quad
\lim_{n\to\infty}\frac{\Gamma(f_{(n)},g_{(n)})}{\Gamma(f_{(n)})}=\infty
\end{equation}
\end{lem}
Obviously, $\Gamma(f')=\Vert f'\Vert^2=S^L_{11}(\mu)=\infty$. Consider the following cases:

{\it case (a)}, when $\Gamma(g')<\infty$ then $\Delta(f_m',g_m')\to\infty$,

{\it case (b)},  when $\Gamma(g')=\infty$ and $\Gamma(f_m')/\Gamma(g_m')\to\infty$ then  $\Delta(f_m',g_m')\to\infty$,

{\it case (c)}, when  $\Gamma(f_m')/\Gamma(g_m')\leq C$ and $\Vert C_1f'+C_2g'\Vert^2=\infty$ for all $(C_1,C_2)\in
{\mathbb R}^2\setminus\{0\}$ then $\Delta(f_m',g_m')\to\infty$ by Lemma~\ref{l.min=proj}.
\begin{rem}
\label{r.x(12)-in-A}
 We can  approximate $x_{12}$ by linear combination of $x_{1k}-x_{12}x_{2k}$ due to
 Lemma~\ref{x_12-in-A} if $\sigma^{(0)}_{12}(\mu)=\infty$ (see (\ref{l.x_12-in-A})).
 The divergence  $ \sigma^{(0)}_{12}(\mu)=\infty$
 follows from the inequality $\alpha_{kn}(0)\alpha_{kn}(1)\leq 1/4$
 based on  the relation $(1-x)x\leq 1/4$ for $x\in[0,1]$ and the divergence
 $ \sigma_{12}(\mu)=\sum_{n=3}^\infty\alpha^2_{2n}(1)=\infty$ which
 follows from Lemma~\ref{l.*}.
 The convergence  $\sum_n(1-c_{1n})^2<\infty$ follows from the fact that
$\Gamma(F)<\infty$ (see (\ref{Gam(G)<infty})).
\end{rem}
 \begin{lem}
 \label{l.*}
 Let $ S_{12}^L(\mu)\!=\!\infty$ and $\sum_n(1\!-\!c_{1n})^2\!<\!\infty$ then
 $\sigma_{12}(\mu):=$\\$\sum_n\alpha^2_{2n}(1)\!=\!\infty.$
 \end{lem}
 \begin{pf}
Using Cauchy-Schwarz inequality $\mid(x,y)\mid\leq \Vert x\Vert\cdot
\Vert y\Vert$ for $x,y\in {\mathbb R}^m$ we get
$$
\Big(\sum_{n=1}^m
\alpha_{2n}(1)(1-c_{1n})
\Big)^2\leq
\Big(\sum_{n=1}^m
\alpha_{2n}^2(1)
\Big)
\Big(\sum_{n=1}^m
(1-c_{1n})^2
\Big),\quad \text{for all}\quad
m\in {\mathbb N}.
$$
Therefore,
$$
\Big(\sum_{n=1}^\infty \alpha_{2n}^2(1)\Big)\geq
\Big(\sum_{n=1}^\infty (1-c_{1n})^2\Big)^{-1}
\Big(\sum_{n=1}^\infty\alpha_{2n}(1)(1-c_{1n})\Big)^2=\infty.
$$
\qed\end{pf}
Finally, in the cases $(a),\,(b)$ or $(c)$ we can approximate $x_{1k}-x_{12}x_{2k}$. Then, by Remark~\ref{r.x(12)-in-A} we can approximate $x_{12}$.  Therefore,  we can
approximate all the  variables $x_{1n},x_{2n+1},\,\,n\geq2$ and the proof is completed.

{\bf 31}.  In the opposite case $(\overline a) \cap (\overline
b)\cap (\overline c)$, i.e., in the case
$\Gamma(f_m')/\Gamma(g_m')\leq C$ and
\begin{equation}
\label{not(a,b,c)} \Vert
f'-tg'\Vert^2=\sum_n(1-c_{1n})[1-t(1-c_{1n})]^2<\infty
\end{equation}
for some $t\in{\mathbb R }\setminus\{0\}$ consider again the
expressions $T_{kn}-I-x_{2k}(T_{2n}-I)$ see (\ref{x-x.x(1)})
\begin{equation}
\label{x-x.x(2)}
T_{kn}-I-x_{2k}(T_{2n}-I)=(x_{1k}-x_{12}x_{2k})\tau_{1n}+(x_{1k}x_{2k}-x_{12}x_{2k})\tau_{1n}\tau_{2n}.
\end{equation}
\begin{lem}
\label{l.a+beta b}
 We can approximate by $\sum_{n=3}^mt_n[T_{kn}-I-x_{2k}(T_{2n}-I)]$ the following expression:
 \begin{equation}
\label{a+beta b}
(x_{1k}-x_{12}x_{2k})+\beta(x_{1k}x_{2k}-x_{12}x_{2k})=x_{1k}-(1-s)x_{12}x_{2k}-sx_{1k}x_{2k}
\end{equation}
where
$s:=-\beta\in [0,\,1]$ and $\beta$ is defined as follows:
\begin{equation}
 \label{beta=}
\beta=\lim_m\beta^{(3)}_m,\quad \beta^{(3)}_m=-\Big(\sum_{n=3}^m\frac{1-c_{1n}}{1+c_{1n}}\Big)^{-1}
\sum_{n=3}^m\frac{(1-c_{1n})(1-c_{2n})}{1+c_{1n}}.
\end{equation}
\end{lem}
\begin{pf}
First, we show that $\sum_{n=3}^mt_n\tau_{1n}\to I$ for an appropriate $t=(t_n)_{n=3}^m$. Second, we show that
$\lim_{m\to\infty}\sum_{n=3}^mt_n\tau_{1n}\tau_{2n}=\beta$.   Indeed, set $b_n:=M\tau_{1n}{\bf 1}$ and
$b_n^{(3)}:=M\tau_{1n}\tau_{2n}{\bf 1}$ then we get
\begin{equation}
\label{}
b_n\!=\!M(T_{1n}-I){\bf 1}\!=\!(c_{1n}-1),\quad
b_n^{(3)}\!=\!M(T_{1n}-I)(T_{\alpha_{2n}}-I){\bf 1}\!=\!(1-c_{1n})(1-c_{2n}).
\end{equation}
By Lemma~\ref{l.P2} we conclude that $\sum_{n=3}^mt_n\tau_{1n}\to I$  if and only if $S^L_{11}(\mu)=\sum_{n=2}^\infty(1-c_{1n})=\infty$. Indeed, if we set $t=(t_n)_{n=3}^m,\,\,b=(b_n)_{n=3}^m$
we get
$$
\sum_{n=3}^mt_n\tau_{1n}\to I \Leftrightarrow
\min_{t\in {\mathbb R}^{m-2}}\Big(\sum_{n=3}^mt_n^2(1-c_{1n})^2\mid (t,b)=1\Big)=
$$
$$
\Big(\sum_{n=3}^m\frac{b_n^2}{a_n}=\sum_{n=3}^m\frac{(1-c_{1n})^2}
{1-c_{1n}^2}=\sum_{n=3}^m\frac{1-c_{1n}}{1+c_{1n}}\Big)^{-1}\to 0
\Leftrightarrow \sum_{n=2}^\infty(1-c_{1n})=\infty,
$$
where
\begin{equation}
\label{a(n),t(n)=}
a_n\!=\!\Vert (T_{1n}\!-\!c_{1n}){\bf 1}\Vert^2=1\!-\!c_{1n}^2, \,\,
t_n\!=\!\frac{b_n}{a_n}\Big(\sum_{n=3}^m\frac{b_n^2}{a_n}\Big)^{-1}\!=\!-\frac{1}{1+c_{1n}}
\Big(\sum_{n=3}^m\frac{1-c_{1n}}{1+c_{1n}}\Big)^{-1}.
\end{equation}
Further, we get
$$
\Vert
\Big[\sum_{n=3}^mt_n\tau_{1n}\tau_{2n}-\beta_m^{(3)}\Big]{\bf 1}
\Vert^2=
\Vert
\Big[\sum_{n=3}^mt_n\tau_{1n}\tau_{2n}-\sum_{n=3}^mt_n(1-c_{1n})(1-c_{2n}) \Big]{\bf 1}
\Vert^2=
$$
$$
\Vert
\sum_{n=3}^mt_n[\tau_{1n}\tau_{2n}-(1-c_{1n})(1-c_{2n})]{\bf 1}
\Vert^2=\sum_{n=3}^mt_n^2\Vert [\tau_{1n}\tau_{2n}-(1-c_{1n})(1-c_{2n})]{\bf 1}\Vert^2
$$
$$
=
\sum_{n=3}^mt_n^2[
2(1-c_{1n})2(1-c_{2n})-(1-c_{1n})^2(1-c_{2n})^2]\sim
\sum_{n=1}^mt_n^2(1-c_{1n})(1-c_{2n}),
$$
since we have $3<4-(1-x)(1-y)\leq 4$ for $x,y\in (0,\,1]$ and
$$
\Vert [\tau_{1n}\tau_{2n}-(1-c_{1n})(1-c_{2n})]{\bf 1}\Vert^2=
\Vert (T_{1n}-I)(T_{\alpha_{2n}}-I){\bf 1}\Vert^2-(1-c_{1n})^2(1-c_{2n})^2=
$$
$$
2(1-c_{1n})2(1-c_{2n})-(1-c_{1n})^2(1-c_{2n})^2=(1-c_{1n})(1-c_{2n})[4-(1-c_{1n})(1-c_{2n})].
$$
We show that $\sum_{n=1}^mt_n^2(1-c_{1n})(1-c_{2n})\to 0$. Indeed, we have
$$
\sum_{n=3}^mt_n^2(1-c_{1n})(1-c_{2n})=\Big(\sum_{n=3}^m\frac{1-c_{1n}}{1+c_{1n}}\Big)^{-2}
\sum_{n=3}^m\frac{(1-c_{1n})(1-c_{2n})}{1+c_{1n}} \leq
$$
$$
4\Big(\sum_{n=3}^m(1-c_{1n})\Big)^{-2}\sum_{n=3}^m(1-c_{1n})=4\Big(\sum_{n=3}^m(1-c_{1n})\Big)^{-1}\to 0.
$$
Obviously, the sequence $\beta_m^{(3)}$  defined by (\ref{beta=}) is bounded $\beta_m^{(3)}\in [-1,\,0]$ for all $m\in {\mathbb N}$ therefore, there exists a subsequence having the limit $\beta\in [-1,\,0]$.
\qed\end{pf}

We have to study two cases: (310), when $s\not=1$ and (311), when  $s=1$.

{\bf 310.} When $s\not=1$ the proof of the irreducibility is finished.
Indeed, we conclude  by  Lemma~\ref{x12inA}  (we shall prove this lemma below) that $x_{12}\in
{\mathfrak A}^2$ if $\sigma_{12}^{(s)}(\mu)=\infty$.
The divergence $\sigma_{12}^{(s)}(\mu)=\infty$ follows from
$\Sigma_k\alpha_{2k}^2(1)=\infty$ and the estimation
$$
\alpha_{rk}(0)\alpha_{rk}(1)\leq 1/4,\quad \alpha_{2k}(0)+(1-s)^2\alpha_{2k}(1)\leq 1.
$$
The divergence $\Sigma_k\alpha_{2k}^2(1)=\infty$ follows from
$\Gamma(F)=\Sigma_n(1-c_{1n})^2<\infty$ and Lemma~\ref{l.*}.

Further, $x_{2k}\in {\mathfrak A}^2$
hence, we get using (\ref{a+beta b})
$x_{1k}-sx_{1k}x_{2k}=x_{1k}(1-sx_{2k})=x_{1k}\otimes
\left(\begin{smallmatrix}
1&0\\
0&1-s
\end{smallmatrix}\right)\in {\mathfrak A}^2$. Finally, we conclude that $x_{1k}\in {\mathfrak
A}^2$ since
$$
(1-sx_{2k})^{-1}=\left(\begin{smallmatrix}
1&0\\
0&1-s
\end{smallmatrix}\right)^{-1}=\left(\begin{smallmatrix}
1&0\\
0&(1-s)^{-1}
\end{smallmatrix}\right)=(I-x_{2k})+(1-s)^{-1}x_{2k}\in {\mathfrak
A}^2.
$$
Now we get $x_{1k},\,\,x_{2k+1}\in {\mathfrak A}^2$ for $k\geq 2$, this finish the proof in the
case $s\not=1$.

We prove Lemma~\ref{x12inA} to finish the case $s\not=1$ before passing to the case $s=1$.
To approximate $x_{12}$, we correct a little bit the expression $x_{1k}-sx_{1k}x_{2k}$.
\begin{lem}
\label{optim(t1,t2)}
For $s\in {\mathbb R}$ we have
$$
\min_{(t_1,t_2)\in{\mathbb R}^2}\Vert (x_{1k}-sx_{1k}x_{2k}+t_1 +
t_2x_{2k}){\bf 1} \Vert^2=\Vert(x_{1k}-Mx_{1k}{\bf
1})(1-sx_{2k}){\bf 1} \Vert^2
$$
$$
=\alpha_{1k}(0)\alpha_{1k}(1)[\alpha_{2k}(0)+
(1-s)^2\alpha_{2k}(1)].
$$
\end{lem}
We see that $(x_{1k}-Mx_{1k}{\bf
1})(1-sx_{2k})=x_{1k}-sx_{1k}x_{2k}-Mx_{1k}{\bf 1}+s(Mx_{1k}{\bf
1})x_{2k}$ hence, minimum we have for $t_1=-Mx_{1k}{\bf
1},\,\,t_2=s(Mx_{1k}{\bf 1})$.
\begin{pf}
We note that the distance $d(f_{n+1};\langle f_1,...,f_n\rangle)$ of
the vector $f_{n+1}$ in a Hilbert space  $H$ from the hyperplane
$\langle f_1,...,f_n\rangle$ generated by vectors $f_1,...,f_n$ may
be calculated in terms of the {\it Gram determinants}
$\Gamma(f_1,f_2,...,f_k)$ corresponding to the set of vectors
$f_1,f_2,...,f_k$ (see \cite{Gan58}):
\begin{equation}
\label{8.dist}
 d^2(f_{n+1};\langle f_1,...,f_n\rangle)=
\min_{t=(t_k)\in{\mathbb R}^n}\Vert f_{n+1}+\sum_{k=1}^nt_kf_k
\Vert^2=
\frac{\Gamma(f_1,f_2,...,f_{n+1})}{\Gamma(f_1,f_2,...,f_n)},
\end{equation}
where the Gram determinant is defined by
$\Gamma(f_1,f_2,...,f_n)={\rm det}\,\gamma(f_1,f_2,...,f_n)$ and
$\gamma(f_1,f_2,...,f_n)$ is the {\it Gram matrix}
\index{Gram!determinant}
\index{Gram!matrix}
$$
\gamma(f_1,f_2,...,f_n)= \left(\begin{array}{cccc}
(f_1,f_1)&(f_1,f_2)&...&(f_1,f_n)\\
(f_2,f_1)&(f_2,f_2)&...&(f_2,f_n)\\
         &         &...&         \\
(f_n,f_1)&(f_n,f_2)&...&(f_n,f_n)
\end{array}\right).
$$
Let us denote $f_0=x_{1k}(1-sx_{2k}){\bf 1},\,\,f_1={\bf 1},\,\,
f_2=x_{2k}{\bf 1}$. We have
$$
\min_{(t_1,t_2)\in{\mathbb R}^2}\Vert f_0+t_1f_1+t_2f_2 \Vert^2=
\frac{\Gamma(f_0,f_1,f_2)}{\Gamma(f_1,f_2)}.
$$
Since we have for operators $x_{1k}$ and $x_{2k}$ (acting on the
spaces $H_{1k}$ and $H_{2k}$ respectively) the same expressions: $
\begin{array}{cc}
\left(\begin{smallmatrix}
0&0\\
0&1
\end{smallmatrix}\right)
\end{array}
$ and $1-sx_{2k}=
\begin{array}{cc}
\left(\begin{smallmatrix}
1&0\\
0&1
\end{smallmatrix}\right)
- \left(\begin{smallmatrix}
0&0\\
0&s
\end{smallmatrix}\right)=
\left(\begin{smallmatrix}
1&0\\
0&1-s
\end{smallmatrix}\right)
\end{array}
$ (to be more precise we write)
$$
x_{1k}=\begin{array}{cc}
\left(\begin{smallmatrix}
0&0\\
0&1
\end{smallmatrix}\right)
\end{array}\otimes \begin{array}{cc}
\left(\begin{smallmatrix}
1&0\\
0&1
\end{smallmatrix}\right)
\end{array},\quad
1-sx_{2k}= \begin{array}{cc} \left(\begin{smallmatrix}
1&0\\
0&1
\end{smallmatrix}\right)
\end{array}\otimes
 \begin{array}{cc}
\left(\begin{smallmatrix}
1&0\\
0&1-s
\end{smallmatrix}\right)
\end{array},
$$
we get
{\small
$$
(f_0,f_0)=\Vert x_{1k}{\bf 1}\Vert^2\Vert(1-sx_{2k}){\bf 1}\Vert^2=
\alpha_{1k}(1)(\alpha_{2k}(0)+ (1-s)^2\alpha_{2k}(1)),
$$
$$
(f_0,f_1)\!=\!(x_{1k}(1\!-\!sx_{2k}){\bf 1},{\bf 1})\!=\!
(x_{1k}{\bf 1},{\bf 1})((1\!-\!sx_{2k}){\bf 1},{\bf 1})\!=\!
\alpha_{1k}(1)(\alpha_{2k}(0)\!+\!(1\!-\!s)\alpha_{1k}(1)),
$$
$$
(f_0,f_2)=(x_{1k}(1-sx_{2k}){\bf 1},x_{2k}{\bf 1})= (x_{1k}{\bf
1},{\bf 1})((1-sx_{2k}){\bf 1},x_{2k}{\bf
1})=\alpha_{1k}(1)(1-s)\alpha_{2k}(1),
$$
$$
(f_1,f_1)=1,\quad (f_1,f_2)=({\bf 1},x_{2k}{\bf
1})=\alpha_{2k}(1),\quad(f_2,f_2)=(x_{2k}{\bf 1},x_{2k}{\bf
1})=\alpha_{2k}(1).
$$
}
Finally, we conclude that
$$
\Gamma(f_0,f_1,f_2)= \left|\begin{array}{ccc}
(f_0,f_0)&(f_0,f_1)&(f_0,f_2)\\
(f_1,f_0)&(f_1,f_1)&(f_1,f_2)\\
(f_2,f_0)&(f_2,f_1)&(f_2,f_2)
\end{array}\right|=
$$
{\small
$$
\left|\begin{array}{ccc}
\alpha_{1k}(1)(\alpha_{2k}(0)\!+\!(1\!-\!s)^2\alpha_{2k}(1))
&\alpha_{1k}(1)(\alpha_{2k}(0)\!+\!(1\!-\!s)\alpha_{2k}(1))&\alpha_{1k}(1)(1\!-\!s)\alpha_{2k}(1)\\
\alpha_{1k}(1)(\alpha_{2k}(0)+ (1-s)\alpha_{2k}(1))&1              &\alpha_{2k}(1)\\
\alpha_{1k}(1)(1-s)\alpha_{2k}(1) &\alpha_{2k}(1) &\alpha_{2k}(1)
\end{array}\right|
$$
$$
\!=\!\alpha_{1k}(1)\alpha_{2k}(1) \left|\begin{array}{ccc}
\alpha_{2k}(0)\!+\!(1\!-\!s)^2\alpha_{2k}(1)
&\alpha_{2k}(0)+ (1-s)\alpha_{2k}(1)&(1-s)\\
\alpha_{1k}(1)(\alpha_{2k}(0)+ (1-s)\alpha_{2k}(1))&1              &1\\
\alpha_{1k}(1)(1-s)\alpha_{2k}(1) &\alpha_{2k}(1) &1
\end{array}\right|
$$
}
$$
=\alpha_{1k}(1)\alpha_{2k}(1) \left|\begin{array}{ccc}
\alpha_{2k}(0)+(1-s)^2\alpha_{1k}(0)\alpha_{2k}(1)
&\alpha_{2k}(0)&0\\
\alpha_{1k}(1)\alpha_{2k}(0)&\alpha_{2k}(0)             &0\\
\alpha_{1k}(1)(1-s)\alpha_{2k}(0) &\alpha_{2k}(1) &1
\end{array}\right|
$$
$$
=\alpha_{1k}(1)\alpha_{2k}(1) \left|\begin{array}{cc}
\alpha_{2k}(0)+(1-s)^2\alpha_{1k}(0)\alpha_{2k}(1)
&\alpha_{2k}(0)\\
\alpha_{1k}(1)\alpha_{2k}(0)&\alpha_{2k}(0)
\end{array}\right|
$$
$$
=\alpha_{1k}(1)\alpha_{2k}(1)\alpha_{2k}(0)\left(
\alpha_{2k}(0)+(1-s)^2\alpha_{1k}(0)\alpha_{2k}(1)-
\alpha_{1k}(1)\alpha_{2k}(0)\right)
$$
$$
=\alpha_{1k}(0)\alpha_{1k}(1)\alpha_{2k}(0)\alpha_{2k}(1)\left(
\alpha_{2k}(0)+(1-s)^2\alpha_{2k}(1)\right).
$$
 For $\Gamma(f_1,f_2)$ we have
$$
\Gamma(f_1,f_2)= \left|\begin{array}{cc}
(f_1,f_1)&(f_1,f_2)\\
(f_2,f_1)&(f_2,f_2)
\end{array}\right|=
\left|\begin{array}{cc}
1&\alpha_{2k}(1)\\
\alpha_{2k}(1) &\alpha_{2k}(1)
\end{array}\right|=\alpha_{2k}(0)\alpha_{2k}(1),
$$
hence, $\Gamma(f_0,f_1,f_2)(\Gamma(f_1,f_2))^{-1}=
\alpha_{1k}(0)\alpha_{1k}(1)\left(
\alpha_{2k}(0)+(1-s)^2\alpha_{2k}(1)\right)$ and
$$
\Vert(x_{1k}-Mx_{1k}{\bf 1})(1-sx_{2k}){\bf 1} \Vert^2=
\Vert(x_{1k}-Mx_{1k}{\bf 1})\Vert^2\Vert (1-sx_{2k}){\bf 1} \Vert^2
$$
$=\alpha_{1k}(0)\alpha_{1k}(1)\left(
\alpha_{2k}(0)+(1-s)^2\alpha_{2k}(1)\right).$ \qed\end{pf}
By Lemma~\ref{optim(t1,t2)} we have for optimal $t_1$ and $t_2$
$$
x_{1k}-sx_{1k}x_{2k}-(1-s)x_{12}x_{2k}+t_1+t_2x_{2k}=
(x_{1k}-Mx_{1k}{\bf 1})(1-sx_{2k})-(1-s)x_{12}x_{2k}.
$$
\begin{lem}
\label{x12inA} For $s\not=1$ we have
{\small
$$
-(1-s)x_{12}{\bf 1}\!\in\!\langle \left[(x_{1k}\!-\!Mx_{1k}{\bf
1})(1-sx_{2k})\!-\!(1-s)x_{12}x_{2k}\right]{\bf 1}\mid k\geq 3 \rangle
\Leftrightarrow \sigma_{12}^{(s)}(\mu)\!=\!\infty,
$$
}
$$
\text{\, where\,}\,\,\sigma_{12}^{(s)}(\mu):=
\sum_{k}\frac{\alpha_{2k}^2(1)}{
\alpha_{2k}(0)\alpha_{2k}(1)+\alpha_{1k}(0)\alpha_{1k}(1)(\alpha_{2k}(0)+
(1-s)^2\alpha_{2k}(1))}.
$$
\end{lem}
\begin{pf} We can procced as before. Let us denote
 $$
 \xi_k=x_{2k}{\bf 1},\quad\eta_k^s=(x_{1k}-Mx_{1k}{\bf
1})(1-sx_{2k}){\bf 1},\quad\text{then}\quad M\xi_k=\alpha_{2k}(1),
$$
$$
\Vert \xi_k-M\xi_k\Vert^2=\alpha_{2k}(0)\alpha_{2k}(1),\quad
\Vert\eta_k^s
\Vert^2=\alpha_{1k}(0)\alpha_{1k}(1)\left(\alpha_{2k}(0)+
(1-s)^2\alpha_{2k}(1)\right).
$$
If we take $(t_k)_k$ such that $\sum_{k=3}^{N+3}t_kM\xi_k=1$ we obtain
$$
\Vert \Big( \sum_{k=3}^{N+3}t_k \left[(x_{1k}-Mx_{1k}{\bf
1})(1-x_{2k})-(1-s)x_{12}x_{2k}\right]+(1-s)x_{12}
 \Big){\bf 1}\Vert^2=
$$
$$
\Vert \sum_{k=3}^{N+3}t_k[ \eta_k^s-(1-s)x_{12}(\xi_k-M\xi_k) ]
\Vert^2= \sum_{k=2}^{N+2}t_k^2\Vert
\eta_k^s-(1-s)x_{12}(\xi_k-M\xi_k)
 \Vert^2
$$
$$
=\sum_{k=3}^{N+3}t_k^2\left(\Vert \eta_k^s\Vert^2+\Vert (1-s)
x_{12}{\bf 1}\Vert^2 \Vert (\xi_k-M\xi_k)\Vert^2\right).
$$
Hence,
$$
\min_{t\in{\mathbb R}^N}\Big(\sum_{k=3}^{N+3}t_k^2\big(\Vert
\eta_k^s\Vert^2+\Vert (1-s)x_{12}{\bf 1}\Vert^2 \Vert
(\xi_k-M\xi_k)\Vert^2\big)
\mid
\sum_{k=3}^{N+3}t_kM\xi_k=1
\Big)
$$
$$
= \Big( \sum_{k=3}^{N+3}\frac{|M\xi_k|^2}{(1-s)^2\Vert x_{12}{\bf
1}\Vert^2\Vert \xi_k-M\xi_k \Vert^2+\Vert\eta_k^s \Vert^2}
\Big)^{-1}
$$
{\small
$$
= \Big( \sum_{k=3}^N\frac{\alpha_{2k}^2(1)}{
(1-s)^2\alpha_{12}(1)\alpha_{2k}(0)\alpha_{2k}(1)+\alpha_{1k}(0)\alpha_{1k}(1)
(\alpha_{2k}(0)+(1-s)^2\alpha_{2k}(1))} \Big)^{-1}\sim
$$
}
$$
 \Big( \sum_{k=3}^N\frac{\alpha_{2k}^2(1)}{
\alpha_{2k}(0)\alpha_{2k}(1)+\alpha_{1k}(0)\alpha_{1k}(1)
(\alpha_{2k}(0)+(1-s)^2\alpha_{2k}(1))}
\Big)^{-1},\text{\, if\,} s\not=1.
$$
\qed\end{pf}
{\bf 311.} When $s=1$ we get from (\ref{a+beta b})
$x_{1k}-x_{1k}x_{2k}$. The condition (\ref{not(a,b,c)})
\begin{equation}
\label{f-tg-in-l(2)} \Vert
f'-tg'\Vert^2=\sum_n(1-c_{1n})[1-t(1-c_{2n})]^2<\infty,\quad\text{for
some}\quad t\in{\mathbb R}\setminus\{0\},
\end{equation}
splits into two cases (3110), when $t\not=1$ and (3111), when $t=1$. We show that in the  first and the second case
we get respectively:
\begin{equation}
\label{t=(not,1)}
\sum_n(1-c_{1n})c_{2n}^2=\infty\quad\text{and}\quad\sum_n(1-c_{1n})\alpha_{2n}^2(1)=\infty.
\end{equation}
To approximate $x_{12}$, under the above conditions we use the following expression
(see (\ref{x-x.x(1)}))  in the first case
$$
T_{kn}-I-x_{2k}(T_{2n}-I)-(x_{1k}-x_{1k}x_{2k})(T_{1n}-I)=
(x_{1k}-x_{12}x_{2k})\tau_{1n}+
$$
\begin{equation}
\label{xx-xx}
(x_{1k}x_{2k}-x_{12}x_{2k})\tau_{1n}\tau_{2n}-
(x_{1k}-x_{1k}x_{2k})\tau_{1n}=
(x_{1k}x_{2k}-x_{12}x_{2k})\tau_{1n}(I+\tau_{2n}).
\end{equation}
In the second case, if  we multiply the latter expression by $T_{2n}=(I+x_{12}\tau_{1n})(I+\tau_{2n})$ we get
(see (\ref{x-xxx}))
\begin{equation}
(x_{1k}x_{2k}-x_{12}x_{2k})\tau_{1n}(I+\tau_{2n})
(I+x_{12}\tau_{1n})(I+\tau_{2n})=(x_{1k}x_{2k}-x_{12}x_{2k}(2x_{1k}-I))\tau_{1n}.
\end{equation}
Consider the case {\bf 3110}, when $t\not=1$.
\begin{lem}
\label{f,g-not-l(2)} Let $f,g\not\in l_2$ where  $f=(f_n)_{n\in{\mathbb N}},\,\,g=(g_n)_{n\in{\mathbb N}}$.
If for some $t\in{\mathbb R}$ holds $tf+(1-t)g\in l_2$ then such a $t$ is unique.
\end{lem}
\begin{pf}
Set $H(t)=tf+(1-t)g,\,\,t\in{\mathbb R}$. Suppose that $H(t_1),\,\,H(t_2)\in l_2$ for two different $t_1$ and $t_2$.Then we
get the contradiction, since for some $s\in {\mathbb R}$ holds
$f=sH(t_1)+(1-s)H(t_2)$ and by assumption we get $l_2\not\ni f=sH(t_1)+(1-s)H(t_2)\in l_2$.
We note that $s=(1-t_2)(t_1-t_2)^{-1}.$
\qed\end{pf}
\begin{rem}
\label{r.t-not=1} The condition (\ref{f-tg-in-l(2)}) for
$t\not=1$ implies the first condition of (\ref{t=(not,1)}).
Indeed, by Lemma~\ref{f,g-not-l(2)} we get  $\Vert f'-g'\Vert^2=\sum_n(1-c_{1n})c_{2n}^2=\infty$ for $t=1$.
\end{rem}
\begin{lem}
\label{l.xx-xx}
We have
$$
(x_{1k}x_{2k}-x_{12}x_{2k}){\bf 1}\in\langle
(x_{1k}x_{2k}-x_{12}x_{2k})\tau_{1n}(I+\tau_{2n}){\bf 1}\mid n>k\rangle
$$
if and only if $\Sigma_{12}^{(1)}:=\sum_{n}(1-c_{1n})c_{2n}^2=\infty$.
\end{lem}
\begin{pf} It is sufficient to show that $\sum_n t_n\big[(T_{\alpha_{1n}}-I)\otimes T_{\alpha_{2n}}\big]\to 1$ if and only if
$\Sigma_{12}^{(1)}=\infty$. Set $\xi_n=\big[(T_{\alpha_{1n}}-I)\otimes T_{\alpha_{2n}}\big]{\bf 1}$ and
$\xi_n^c=\xi_n-M\xi_n$, then
$$
M\xi_n=(c_{1n}-1)c_{2n},\quad \Vert\xi_n\Vert^2=2(1-c_{1n}),\quad \Vert\xi_n^c\Vert^2=\Vert\xi_n\Vert^2-\mid M\xi_n\mid^2.
$$
Indeed, we have
$$
\Vert\xi_n\Vert^2=
\Vert\big[(T_{\alpha_{1n}}-I)\otimes T_{\alpha_{2n}}\big]{\bf 1}\Vert^2 =
$$
$$
\Vert(T_{\alpha_{1n}}-I){\bf 1}\Vert^2=\Vert T_{\alpha_{1n}}{\bf 1}\Vert^2-2(T_{\alpha_{1n}}{\bf 1},{\bf 1})+1=
2(1-c_{1n}).
$$
Take $(t_n)_n$ such that $\sum_{n=2}^Nt_nM\xi_n=1$ then
{\small
$$
\Vert\Big(\sum_{n=2}^{N+2}t_n\big[(T_{\alpha_{1n}}-I)\otimes T_{\alpha_{2n}}\big]-I\Big){\bf 1}\Vert^2=
\Vert\Big(\sum_{n=2}^{N+2}t_n\big[(T_{\alpha_{1n}}-I)\otimes T_{\alpha_{2n}}\big]-\sum_{n=2}^{N+2}t_nM\xi_n\Big){\bf 1}\Vert^2
$$
}
$$
=\sum_{n=2}^{N+2}t_n^2\Vert\big[(T_{\alpha_{1n}}-I)\otimes T_{\alpha_{2n}}-M\xi_n\big]{\bf 1}\Vert^2
=\sum_{n=2}^{N+2}t_n^2\Vert\xi_n^c\Vert^2.
$$
Finally, we get
$$
\min_{t\in{\mathbb R}^N}\Big(\sum_{n=2}^{N+2}t_n^2\Vert\xi_n^c\Vert^2
\mid \sum_{n=2}^{N+2}t_nM\xi_n=1
\Big)=(\Sigma_{12,N}^{(1)})^{-1}\quad\text{where}
$$
$$
\Sigma_{12,N}^{(1)}:=\sum_{n=2}^{N+2}\frac{\mid M\xi_n\mid^2}{\Vert\xi_n^c\Vert^2}=
\sum_{n=2}^{N+2}\frac{\mid M\xi_n\mid^2}{\Vert\xi_n\Vert^2-\mid M\xi_n\mid^2}\sim
$$
$$
\sum_{n=2}^{N+2}\frac{\mid M\xi_n\mid^2}{\Vert\xi_n\Vert^2}=
\sum_{n=2}^{N+2}\frac{(1-c_{1n})^2c_{2n}^2}{2(1-c_{1n})}=\frac{1}{2}\sum_{n=2}^{N+2}(1-c_{1n})c_{2n}^2.
$$
\qed\end{pf}
 Now we get $x_{1k}x_{2k}-x_{12}x_{2k},\,\,x_{1k}-x_{1k}x_{2k}\in {\mathfrak A}^2$
hence, $x_{1k}-x_{12}x_{2k}\in{\mathfrak A}^2$. Using Lemma~\ref{x12inA} for $s=0$ we get
\begin{lem}
\label{x_12-in-A} We have $ x_{12}{\bf
1}\in\langle(x_{1k}-x_{12}x_{2k}){\bf
1}\mid n>k\rangle $ if and only if $\sigma_{12}^{(0)}(\mu)=\infty$
where
\begin{equation}
\label{l.x_12-in-A}
\sigma_{12}^{(0)}(\mu)=\sum_{k}\frac{\alpha_{2k}^2(1)}{
\alpha_{2k}(0)\alpha_{2k}(1)+\alpha_{1k}(0)\alpha_{1k}(1)}.
\end{equation}
\end{lem}
We use the following obvious implications:
$$
\Gamma(F)<\infty\stackrel{(\ref{Gam(G)<infty})}{\Rightarrow} \sum_n(1-c_{1n})^2<\infty\stackrel{\text{Lemma~\ref{l.*}}}{\Rightarrow} \sum_n\alpha_{2n}^2(1)=\infty
\Rightarrow \sigma_{12}^{(0)}(\mu)=\infty.
$$
By  Lemma~\ref{x_12-in-A},  we conclude
that $x_{12}\in {\mathfrak A}^2$ hence, $x_{1k},\,\,x_{2k+1}\in
{\mathfrak A}^2$ for all $k\geq 2$. This finish the proof in this
case.

Consider the case {\bf 3111}, when $t=1$. Since $p=2$ we get
$(I+\tau_{2n})^2=T_{2n}^2(2)=T_{\alpha_{2n}}^2=I$ and
$\tau_{1n}^2=-2\tau_{1n}$. Indeed, we get
$$
\tau_{1n}^2=(T_{\alpha_{1n}}-I)^2\!=\!T_{\alpha_{1n}}^2-2T_{\alpha_{1n}}+I=-2(T_{\alpha_{1n}}-I)=-2\tau_{1n}.
$$
Hence, we have
$$
(x_{1k}x_{2k}\!-\!x_{12}x_{2k})\tau_{1n}(I+\tau_{2n})(I+x_{12}\tau_{1n})(I+\tau_{2n})\!=\!
(x_{1k}x_{2k}-x_{12}x_{2k})\tau_{1n}(I+x_{12}\tau_{1n})
$$
\begin{equation}
\label{x-xxx}
\!=\![x_{1k}x_{2k}\!-\!x_{12}x_{2k}\!-\!2x_{12}(x_{1k}x_{2k}-x_{12}x_{2k})]\tau_{1n}
\!=\![x_{1k}x_{2k}\!-\!x_{12}x_{2k}(2x_{1k}\!-\!I)]\tau_{1n}.
\end{equation}
The condition $S_{11}^L(\mu)=\infty$ implies $\sum_{n=1}^mt_n\tau_{1n}\to I$ therefore,
$$
x_{1k}x_{2k}-x_{12}x_{2k}(2x_{1k}-I)\in {\mathfrak A}^2.
$$
Since
$x_{1k}-x_{1k}x_{2k}\in {\mathfrak A}^2$ we conclude finally, that
$x_{1k}-x_{12}x_{2k}(2x_{1k}-I)\in {\mathfrak A}^2$.
\begin{rem}
\label{r.x-xxx}
The condition (\ref{f-tg-in-l(2)}) for $t=1$ means that
$\sum_n(1-c_{1n})c_{2n}^2<\infty$. This implies that
$\sum_n(1-c_{1n})\alpha_{2n}^2(1)=\infty$. Indeed, otherwise, if
we suppose that \\$\sum_n(1\!-\!c_{1n})\alpha_{2n}^2(1)\!<\!\infty$
we obtain the contradiction with the condition
$S_{12}^L(\mu)=\sum_n(1-c_{1n})\alpha_{2n}(1)=\infty$. In fact,
since $c_{2n}^2=4\alpha_{2n}(0)\alpha_{2n}(1)$,  we get
$$
\infty>\sum_n(1-c_{1n})c_{2n}^2=4\sum_n(1-c_{1n})(\alpha_{2n}(1)-\alpha_{2n}^2(1)).
 $$
By Lemma~\ref{l.x(12)-in-A!} below, we conclude that $x_{12}\in {\mathfrak A}^2$. Since $x_{2k}\in {\mathfrak A}^2$ we conclude,  by Remark~\ref{r.x(kn)-in-A} (c), that  $x_{1k}\in  {\mathfrak A}^2$ hence, $x_{1k},\,x_{2k+1}\in  {\mathfrak A}^2$ for $k\geq 2$ and the proof is finished.
\end{rem}
\begin{lem}
\label{l.x(12)-in-A!} We have $ x_{12}{\bf
1}\in\langle[x_{1k}-x_{12}x_{2k}(2x_{1k}-I)]{\bf 1}\mid k>2\rangle
$ if and only if $\sum_n(1-c_{1n})\alpha_{2n}^2(1)=\infty$.
\end{lem}
\begin{pf}
Set $\eta_k=x_{1k}{\bf 1}$ and $\xi_k=x_{2k}(2x_{1k}-I){\bf 1}$ then
$$
M\eta_k=\alpha_{1k}(1),\,\,\, M\xi_k =\alpha_{2k}(1)(-\alpha_{1k}(0)+\alpha_{1k}(1)),\quad
$$
$$
\Vert \eta_k\Vert^2=\alpha_{1k}(1),\quad\Vert \xi_k\Vert^2=\alpha_{2k}(1)(\alpha_{1k}(0)+\alpha_{1k}(1))=\alpha_{2k}(1),
$$
since
$$
2x_{1k}-I=\left(\begin{smallmatrix}
0&0\\
0&2
\end{smallmatrix} \right)-\left(\begin{smallmatrix}
1&0\\
0&1
\end{smallmatrix} \right)=\left(\begin{smallmatrix}
-1&0\\
0&1
\end{smallmatrix} \right).
$$
Set $h_n=\eta_k- M\eta_k-x_{12}(\xi_{k}-M\xi_{k})$ then
$$a_n=\Vert h_n\Vert^2=\Vert \eta_k- M\eta_k-x_{12}(\xi_{k}-M\xi_{k})\Vert^2.
$$
We have
$$
\Vert \sum_{n=1}^mt_n[(x_{1k}- Mx_{1k})-x_{12}\xi_{k}]{\bf 1}-x_{12}{\bf 1}\Vert^2=
$$
$$
\Vert \sum_{n=1}^mt_n[(x_{1k}- Mx_{1k})-x_{12}(\xi_{k}-M\xi_{k})]{\bf 1}\Vert^2=
\sum_{n=1}^mt_n^2a_n.
$$
To calculate $a_n$ we get
$$
a_n=\Vert h_n\Vert^2=\Vert [(\eta_k- M\eta_k)-x_{12}(\xi_{k}-M\xi_{k})]{\bf 1}\Vert^2=
\Vert\eta_k- M\eta_k\Vert^2+
$$
$$
\Vert x_{12}(\xi_{k}-M\xi_{k}){\bf 1}\Vert^2-(x_{12}{\bf 1},{\bf 1})(\eta_k- M\eta_k,\xi_{k}-M\xi_{k})=
$$
$$
\Vert x_{1k}{\bf 1}\Vert^2-\mid Mx_{1k}{\bf 1}\mid^2+\frac{1}{2}[\Vert \xi_{k}\Vert^2-\mid M\xi_k\mid^2]-(\eta_k- M\eta_k,\xi_{k}-M\xi_{k})=
$$
$$
\alpha_{1k}(1)-\alpha_{1k}^2(1)+\frac{1}{2}[\alpha_{2k}(1)-\alpha_{2k}^2(1)(-\alpha_{1k}(0)+\alpha_{1k}(1))^2]-
(\eta_k- M\eta_k,\xi_{k}-M\xi_{k}).
$$
Since
$$
(\eta_k- M\eta_k,\xi_{k}-M\xi_{k})=(\eta_k,\xi_{k})-M\eta_kM\xi_{k}=(x_{1k}{\bf 1},x_{2k}(2x_{1k}-I){\bf 1})-M\eta_kM\xi_{k}
$$
$$
=(x_{2k}{\bf 1},x_{1k}{\bf 1})-M\eta_kM\xi_{k}=\alpha_{1k}(1)\alpha_{2k}(1)-\alpha_{1k}(1)\alpha_{2k}(1)(-\alpha_{1k}(0)+\alpha_{1k}(1))
$$
$$
=\alpha_{1k}(1)\alpha_{2k}(1)(1+\alpha_{1k}(0)-\alpha_{1k}(1))=2\alpha_{1k}(0)\alpha_{1k}(1)\alpha_{2k}(1),
$$
we conclude that
$$
a_k\!=\!\alpha_{1k}(1)\!-\!\alpha_{1k}^2(1)\!+\!\frac{1}{2}[\alpha_{2k}(1)\!-\!\alpha_{2k}^2(1)(\alpha_{1k}(0)-\alpha_{1k}(1))^2]-
2\alpha_{1k}(0)\alpha_{1k}(1)\alpha_{2k}(1)
$$
$$
\sim
a_k':=\alpha_{1k}(0)\alpha_{1k}(1)+\frac{1}{2}\alpha_{2k}(1)[1-4\alpha_{1k}(0)\alpha_{1k}(1)]
=\frac{c_{1k}^2}{4}+\frac{1}{2}\alpha_{2k}(1)[1-c_{1k}^2].
$$
Finally, we get
$$
\min_{t\in{\mathbb R}^m}\Big(\sum_{n=1}^mt_n^2a_n\mid
\sum_{n=1}^mt_nb_n=1\Big)=\Big(\sum_{n=1}^m\frac{b_n^2}{a_n}\Big)^{-1}\quad\text{where}\quad
b_n=\frac{1}{2}M\xi_k\quad \text{and}
$$
$$
\sum_{n=1}^\infty\frac{b_n^2}{a_n}\sim\sum_{n=1}^\infty\frac{b_n^2}{a_n'}=
\sum_{k=1}^\infty\frac{\frac{1}{4}
\alpha_{2k}^2(1)(-\alpha_{1k}(0)+\alpha_{1k}(1))^2}
{\alpha_{1k}(0)\alpha_{1k}(1)+\frac{1}{2}\alpha_{2k}(1)[1-4\alpha_{1k}(0)\alpha_{1k}(1)]}=
$$
$$
\sum_{k=1}^\infty\frac{ \frac{1}{4}  \alpha_{2k}^2(1))(1-c_{1k}^2)
}{ \frac{c_{1k}^2}{4}+\frac{1}{2}\alpha_{2k}(1)(1-c_{1k}^2) }\geq
\frac{1}{3}\sum_{k=1}^\infty\alpha_{2k}^2(1))(1-c_{1k}^2)\sim
\sum_{k=1}^\infty\alpha_{2k}^2(1)(1-c_{1k}),
$$
since $1\leq 1+c_{1k}< 2,\,\,\,c_{1k}^2\leq 1$ and $\alpha_{2k}(1)[1-c_{1k}^2]<1$.
We use the following relation for $x=\alpha_{1k}(1)$:
$$
(-\alpha_{1k}(0)+\alpha_{1k}(1))^2=(1-2x)^2=1-4x(1-x)=1-c_{1k}^2.
$$
\qed\end{pf}
{\bf Case 4.}  Let $\Sigma_{12}<\infty$ and $\Gamma(F)=\Gamma(G)=\infty$.
Condition $\Sigma_{12}<\infty$ implies
$\Gamma(H_2)<\infty$ hence, $\Sigma_{2m}\sim \Delta(G_m,F_m)$. We
have two cases:

(4I),$\,\,\,$ when   $\Delta(G_m,F_m)\to \infty$;
$$\quad{\rm (4II)},\,\, \text{when}\,\,
\Gamma(G_m)/\Gamma(F_m)
\leq C\quad \text{and}\,\,\,\Vert G-tF\Vert^2<\infty\,\,\text {for some}\,\,\, t\in{\mathbb R}\setminus\{0\}.$$
In the first case  (4I),  we can approximate $x_{2k}$ and we
are reduced to the case   (3) but some particular cases should
be considered in addition.

In the second case  (4II), we show that by linear combinations of
the expressions
$$T_{kn}-I=x_{1k}\tau_{1n}+x_{2k}\tau_{2n}+x_{1k}x_{2k}\tau_{1n}\tau_{2n}$$
we can approximate $x_{1k}+\beta^{(2)}_1x_{2k}-\beta^{(3)}_1x_{1k}x_{2k}$
or $\beta^{(1)}_2x_{1k}+x_{2k}-\beta^{(3)}_2x_{1k}x_{2k}$,
see Lemma~\ref{x+x+x}.

Case {\bf 4I.} We follow step by step the  case (3) with the
same notations, just replacing $3$ by $4$. We know that
$\Delta(G_m,F_m)\to \infty$ in two cases, due to
Lemmas~\ref{l.min=proj0} and \ref{l.min=proj}:
\par{\it case (b)}, when $\Gamma(G_m)/\Gamma(F_m)\to\infty$,
\par{\it case (c)}, when $\Gamma(G_m)/\Gamma(F_m)\leq C$ for all  $m>3$ and
 $ \Vert C_1F+C_2G\Vert^2=\infty$ for all $(C_1,C_2)\in {\mathbb R}^2\setminus\{0\}$.

 Consider the following expression  (see (\ref{x-x.x(2)}))
$$
T_{kn}-I-x_{2k}(T_{2n}-I)=
(x_{1k}-x_{12}x_{2k})\tau_{1n}+(x_{1k}x_{2k}-x_{12}x_{2k})\tau_{1n}\tau_{2n}.
$$

The case (4I) splits into two cases (40), when $\Delta(f_m',g_m')\to\infty$
and (41), when
$$
\Gamma(f_m')/\Gamma(g_m')\leq C\quad\text{and}\quad
\Vert f'-tg'\Vert^2=\sum_n(1-c_{1n})[1-t(1-c_{1n})]^2<\infty
$$
for some $t\in {\mathbb R}\setminus\{0\}$ (see (\ref{not(a,b,c)})),  as in  the  case (31).

In the case {\bf 40}, by Lemma~\ref{x(1k)-x(12)x(2k)}, we conclude that
$x_{1k}-x_{12}x_{2k}\in {\mathfrak A}^2$ since $\Delta(f_m',g_m')\to\infty$.
The case (40) splits into  two cases: (400), when $\sigma_{12}(\mu)=\sum_n\alpha_{2n}^2(1)$ $=\infty$ and the case
(401), when $\sigma_{12}(\mu)=\sum_n\alpha_{2n}^2(1)<\infty$. In the case (400) we can approximate $x_{12}$ and
the proof is finished.

In the case {\bf 401}, the condition $\sigma_{12}(\mu)=\sum_n\alpha_{2n}^2(1)<\infty$ implies  $\lim_nc_{2n}=0$ indeed,
$$
\lim_nc^2_{2n}=\lim_n4\alpha_{2n}(1)(1-\alpha_{2n}(1))=0.
$$
Since $\sum_n(1-c_{1n})\alpha_{2n}^2(1)< \sum_n\alpha_{2n}^2(1)<\infty$ we conclude that  $\sum_n(1-c_{1n})c_{2n}^2=\infty$. Indeed, use the fact that $S_{12}^L(\mu)=\sum_n(1-c_{1n})\alpha_{2n}(1)=\infty$  and consider the equality
$$
\sum_n(1-c_{1n})c_{2n}^2=4\sum_n(1-c_{1n})(\alpha_{2n}(1)-\alpha_{2n}^2(1)).
$$
\begin{ex}
\label{ex.401}
Let $1-c_{1n}=\frac{1}{n^\beta},\,\,\, c_{2n}=\frac{1}{n^\alpha}$ where $\alpha,\beta>0$. We show that the conditions of the divergence of the following series, which gives us the case (401),
{\small
$$
\sum_{n}(1-c_{1n})\!=\!\sum_{n}(1-c_{2n})\!=\!\sum_{n}(1-c_{1n})\alpha_{2n}(2)\!=\!
\sum_{n}(1-c_{1n})c_{2n}^2=\Vert F\Vert^2\!=\!\Vert G\Vert^2\!=\!\infty
$$
}
are as follows:
\begin{equation}
\label{D:=}
D=\{(\alpha,\beta)\in {\mathbb R}^2\mid 2\alpha+\beta\leq 1,\,\,2\alpha+2\beta> 1,\,\,4\alpha>1\}.
\end{equation}
\end{ex}
Indeed, we have
$$
S^L_{11}(\mu)=\sum_{n}(1-c_{1n})=\sum_n \frac{1}{n^\beta}=\infty\quad\text{for}\quad\beta\in(0,\,1],
$$
$$
S^L_{22}(\mu)=\sum_{n}(1-c_{2n})=\sum_n(1- \frac{1}{n^\alpha})=\infty\quad\text{for}\quad\alpha>0.
$$
To find $x=\alpha_{2n}(1)$, we use the identity $c_{2n}^2=4\alpha_{2n}(0)\alpha_{2n}(1)=4(1-x)x$ (see notation $c_{1n}$ before Lemma~\ref{1in<2>}). We have $x^2-x+\frac{c_{2n}^2}{4}$. The roots are as follows:
\begin{equation}
\label{x1,x2=}
x_1\!=\!\Big(1\!-\!\sqrt{1-c_{2n}^2}\Big)/2\!=\!c_{2n}^2
\Big(2\big(1\!+\!\sqrt{1\!-\!c_{2n}^2}\big)\Big)^{-1}
,\,\,x_2\!=\!\Big(1\!+\!\sqrt{1\!-\!c_{2n}^2}\Big)/2.
\end{equation}
Only the first root is suitable since $\alpha_{2n}(1)\to 0$.
We have $\alpha_{2n}(1)\sim c_{2n}^2/4$. Therefore,
$$
S^L_{12}(\mu)\!=\!\sum_{n}(1-c_{1n})\alpha_{2n}(1)\sim\sum_{n}(1-c_{1n})c_{2n}^2=
\sum_{n}\frac{1}{n^{2\alpha+\beta}}\!=\!\infty\,\,\text{for}\,\,2\alpha+\beta\leq 1,
$$
{\small
$$
\Sigma_{12}\!=\!\sum_{n=2}^\infty\frac{(1-c_{1n})^2c_{2n}^2}{(1-c_{1n})c_{2n}^2+1-c_{2n}}
\!\sim\!\sum_{n}(1\!-\!c_{1n})^2c_{2n}^2\!=\!
\sum_{n}\frac{1}{n^{2\alpha+2\beta}}<\infty\,\,\,\,\text{for}\,\,\,\,2\alpha\!+\!2\beta\!>\!1,
$$
}
$$
\sigma_{12}(\mu)=\sum_{n}\alpha^2_{2n}(1)\sim\sum_{n}c_{2n}^4=\sum_{n}\frac{1}{n^{4\alpha}}<\infty
\quad\text{for}\quad 4\alpha>1,
$$
that proves (\ref{D:=}). Further, we get
$$
\Vert F_m\Vert^2\sim \sum_{n=1}^m\frac{(1-c_{1n})^2}{1-c_{1n}+1-c_{2n}}\sim
\sum_{n=1}^m(1-c_{1n})^2=\sum_{n=1}^m\frac{1}{n^{2\beta}}\sim\frac{m^{1-2\beta}}{1-2\beta}
\to\infty,
$$
$$
\Vert G_m\Vert^2\sim \sum_{n=1}^m\frac{(1-c_{2n})^2}{1-c_{1n}+1-c_{2n}}\sim
\sum_{n=1}^m(1-c_{2n})^2=\sum_{n=1}^m(1-\frac{1}{n^{\alpha}})^2\sim m
\to\infty.
$$
Therefore,
$$\Delta(G_m,F_m)\to\infty\quad\text{since}\quad \Vert G_m\Vert^2/\Vert F_m\Vert^2\sim (1-2\beta)m^{2\beta}\to\infty.
$$
In addition we get for $f_m'=(\sqrt{1-c_{1n}})_{n=1}^m$ and $g_m'=(\sqrt{1-c_{1n}}(1-c_{1n}))_{n=1}^m$
$$
\Vert f_m'\Vert^2/\Vert g_m'\Vert^2=\sum_{n=1}^m(1-c_{1n})/\sum_{n=1}^m(1-c_{1n})(1-c_{2n})^2
\to \lim_n(1-c_{2n})^{-2}=1.
$$
For all $t\in {\mathbb R}$ we have
$\Vert g_m'-tf_m'\Vert^2= \sum_{n=1}^m(1-c_{1n})(1-c_{2n}-t)^{2}\to\infty.$
Indeed, for $t=1$ we get
$
\Vert g_m'-f_m'\Vert^2= \sum_{n=1}^m(1-c_{1n})c_{2n}^2\to\infty
$ hence, $\Vert g_m'-tf_m'\Vert^2\to \infty$ for all $t\in {\mathbb R}$. Therefore,
$\Delta(f_m',g_m')\to\infty$ by Lemma~\ref{l.min=proj} and we are in the case (401).
\begin{lem}
\label{l.final}
We have
\begin{equation}
\label{final}
x_{12}{\bf 1}\in \langle\left[T_{2n},x_{1n}-x_{12}x_{2n}\right]x_{2n}{\bf 1}\mid n\geq 3\rangle
\end{equation}
if and only if $\Sigma_{12}^{(2)}=\sum_{n}\alpha_{1n}(1)\alpha_{2n}(0)=\infty$.
\end{lem}
\begin{pf} Recall (see (\ref{8.ab})), that
$$
a_n=
\sqrt{\alpha_{1n}(0)/\alpha_{1n}(1)}
\quad
b_n=
\sqrt{\alpha_{2n}(0)/\alpha_{2n}(1)}.
$$
We show that
\begin{equation}
\label{[T,x-xx1]}
\left[T_{2n},x_{1n}-x_{12}x_{2n}\right]=2x_{12}\left(
\left(\begin{smallmatrix} 0&a_n^{-1}\\
0&0\\
\end{smallmatrix}\right)\otimes
\left(\begin{smallmatrix} 0&0\\
b_n&0\\
\end{smallmatrix}\right)-
\left(\begin{smallmatrix} 0&0\\
a_n&0\\
\end{smallmatrix}\right)\otimes
\left(\begin{smallmatrix} 0&b_n^{-1}\\
0&0\\
\end{smallmatrix}
\right)
\right),
\end{equation}
therefore,
\begin{equation}
\label{[T,x-xx2]}
\left[T_{2n},x_{1n}-x_{12}x_{2n}\right]x_{2n}=
-2x_{12}\left(\begin{smallmatrix} 0&0\\
a_n&0\\
\end{smallmatrix}\right)\otimes
\left(\begin{smallmatrix} 0&b_n^{-1}\\
0&0\\
\end{smallmatrix}
\right).
\end{equation}
Indeed, since
$$
x_{1n}=
\left(\begin{smallmatrix} 0&0\\
0&1\\
\end{smallmatrix}\right),\quad
T_{\alpha_{1n}}=
\left(\begin{smallmatrix} 0&a_n^{-1}\\
a_n&0\\
\end{smallmatrix}\right),\quad
T_{\alpha_{2n}}=
\left(\begin{smallmatrix} 0&b_n^{-1}\\
b_n&0\\
\end{smallmatrix}\right),
$$
we get
$$
[T_{\alpha_{1n}},x_{1n}]=\left(\begin{smallmatrix} 0&a_n^{-1}\\
-a_n&0\\
\end{smallmatrix}\right),\quad
[T_{\alpha_{2n}},x_{2n}]=\left(\begin{smallmatrix} 0&b_n^{-1}\\
-b_n&0\\
\end{smallmatrix}\right).
$$
Using the (\ref{T(kn)=otimes}) and (\ref{T2n(1)-1.0}) we get
$$
T_{2n}=T_{2n}(1)\otimes T_{2n}(2)=\big(x_{12}(T_{\alpha_{1n}}-I)+I\big)\otimes T_{\alpha_{2n}}
$$
that implies (\ref{[T,x-xx1]}). Indeed, we have
$$
\left[T_{2n},x_{1n}-x_{12}x_{2n}\right]=
\left[\big(x_{12}(T_{\alpha_{1n}}-I)+I\big)\otimes T_{\alpha_{2n}},x_{1n}-x_{12}x_{2n}\right]=
$$
$$
x_{12}\Big([T_{\alpha_{1n}},x_{1n}]T_{\alpha_{2n}}-T_{\alpha_{1n}}
[T_{\alpha_{2n}},x_{2n}]\Big)=
$$
$$
x_{12}\left(
\left(\begin{smallmatrix} 0&a_n^{-1}\\
-a_n&0\\
\end{smallmatrix}\right)\otimes
\left(\begin{smallmatrix} 0&b_n^{-1}\\
b_n&0\\
\end{smallmatrix}\right)-
\left(\begin{smallmatrix} 0&a_n^{-1}\\
a_n&0\\
\end{smallmatrix}\right)\otimes
\left(\begin{smallmatrix} 0&b_n^{-1}\\
-b_n&0\\
\end{smallmatrix}
\right)
\right),
$$
$$
x_{12}\left(
\left(\begin{smallmatrix}
0&0&0&a_n^{-1} b_n^{-1}\\
0&0&a_n^{-1} b_n&0\\
0&-a_n b_n^{-1}&0&0\\
-a_nb_n&0&0&0\\
\end{smallmatrix}\right)-
\left(\begin{smallmatrix}
0&0&0&a_n^{-1} b_n^{-1}\\
0&0&-a_n^{-1} b_n&0\\
0&a_n b_n^{-1}&0&0\\
-a_nb_n&0&0&0\\
\end{smallmatrix}\right)
\right)=
$$
$$
2x_{12}
\left(\begin{smallmatrix}
0&0&0&0\\
0&0&a_n^{-1} b_n&0\\
0&-a_n b_n^{-1}&0&0\\
0&0&0&0\\
\end{smallmatrix}\right)=
2x_{12}\left(
\left(\begin{smallmatrix} 0&a_n^{-1}\\
0&0\\
\end{smallmatrix}\right)\otimes
\left(\begin{smallmatrix} 0&0\\
b_n&0\\
\end{smallmatrix}\right)-
\left(\begin{smallmatrix} 0&0\\
a_n&0\\
\end{smallmatrix}\right)\otimes
\left(\begin{smallmatrix} 0&b_n^{-1}\\
0&0\\
\end{smallmatrix}
\right)
\right).
$$
Now we show that
$\sum_{n=3}^mt_n\left(\begin{smallmatrix} 0&0\\
a_n&0\\
\end{smallmatrix}\right)\otimes
\left(\begin{smallmatrix} 0&b_n^{-1}\\
0&0\\
\end{smallmatrix}
\right)\to I$ if and only if \\ $\sum_{n}\alpha_{1n}(1)\alpha_{2n}(0)=\infty$.

Indeed, set $$
\Xi_n=\left(\begin{smallmatrix} 0&0\\
a_n&0\\
\end{smallmatrix}\right)\otimes
\left(\begin{smallmatrix} 0&b_n^{-1}\\
0&0\\
\end{smallmatrix}
\right)\quad\text{and}\quad
\xi_n=\left(\begin{smallmatrix} 0&0\\
a_n&0\\
\end{smallmatrix}\right)\otimes
\left(\begin{smallmatrix} 0&b_n^{-1}\\
0&0\\
\end{smallmatrix}
\right){\bf 1}.
$$
We get
$$
\Vert
\Big(\sum_{n=3}^mt_n\Xi_n-I\Big){\bf 1}
\Vert^2=
\Vert\sum_{n=3}^mt_n(\xi_n-M\xi_n)
\Vert^2=\sum_{n=3}^mt_n^2\Vert \xi_n-M\xi_n\Vert^2\to 0
$$
under the condition $\sum_{n=3}^mt_nM\xi_n=1$ if and only if
$\sum_n\frac{b_n^2}{a_n}\sim \Sigma_{12}^{(2)}=\infty$. Indeed, we have
{\small
$$
\sum_n\frac{b_n^2}{a_n}\!=\!\sum_n\frac{c_{1n}^2c_{2n}^2/16}
{\alpha_{1n}(0)\alpha_{2n}(1)\!-\!c_{1n}^2c_{2n}^2/16}\!\sim\!
\sum_n\frac{\alpha_{1n}(0)\alpha_{1n}(1)\alpha_{2n}(0)\alpha_{2n}(1)}
{\alpha_{1n}(0)\alpha_{2n}(1)}\!=\!\Sigma_{12}^{(2)}\!=\!
\infty,
$$
}
 where
$$
b_n=M\xi_n,\quad a_b=\Vert \xi_n-M\xi_n\Vert^2=
\Vert \xi_n\Vert^2-\vert M\xi_n\vert^2,
$$
$$
b_n=M\xi_n=\Big(\left(\begin{smallmatrix} 0&0\\
a_n&0\\
\end{smallmatrix}\right){\bf 1},{\bf 1}\Big)
\Big(
\left(\begin{smallmatrix} 0&b_n^{-1}\\
0&0\\
\end{smallmatrix}
\right){\bf 1},{\bf 1}\Big)=c_{1n}c_{2n}/4,
$$
$$
\Vert \xi_n\Vert^2=\Vert
\left(\begin{smallmatrix} 0&0\\
a_n&0\\
\end{smallmatrix}\right){\bf 1}
\Vert^2\Vert
\left(\begin{smallmatrix} 0&b_n^{-1}\\
0&0\\
\end{smallmatrix}
\right){\bf 1}
\Vert^2=\alpha_{1n}(0)\alpha_{2n}(1).
$$
\qed\end{pf}
The condition $\Sigma_{12}^{(2)}\!=\!\sum_{n}\alpha_{1n}(1)\alpha_{2n}(0)\!=\!\infty$ follows from two facts:

(a) $\lim_n\alpha_{2n}(0)\!=\!1$ since $\lim_nc_{2n}=0$ and

(b) $\lim_k\alpha_{1n_k}(1)=1/2$ since $\lim_kc_{1n_k}=1$ (see (\ref{x1,x2=})), that follows from
\begin{equation}
\label{c(1n)-to-1}
\Sigma_{12}\sim\sum_{n}(1-c_{1n})^2c_{2n}^2<\infty\quad\text{and}\quad
S^L_{12}\sim \sum_{n}(1-c_{1n})c_{2n}^2=\infty.
\end{equation}
 The first equivalence follows from $\lim_n[(1-c_{1n})c_{2n}^2+1-c_{2n}]=1$.
Indeed, the condition $1-c_{1n}\geq \varepsilon>0$ contradicts (\ref{c(1n)-to-1})
therefore, for some subsequence $(n_k)_k$ we have $\lim_kc_{1n_k}=1$ hence, $\Sigma_{12}^{(2)}\sim \sum_{n}\alpha_{1n}(1)>\sum_{k}\alpha_{1n_k}(1)=\infty.$

In the case {\bf 41}, when $\Vert f'-tg'\Vert^2<\infty$, by
Lemma~\ref{l.a+beta b}, we can approximate the following
expression:
\begin{equation*}
x_{1k}-(1-s)x_{12}x_{2k}-sx_{1k}x_{2k}.
\end{equation*}
We have two cases: (410), when $s\not=1$ and (411), when $s=1$.
The case (410) splits into two cases (4100), when
$\sum_k\alpha_{2k}^2(1)=\infty$ and (4101),  when
$\sum_k\alpha_{2k}^2(1)<\infty$.

In the case {\bf 4100} we can approximate $x_{12}$ and the proof
is finished.

In the case {\bf 4101}  the condition
$\sum_k\alpha_{2k}^2(1)<\infty$ implies that $\lim_nc_{2n}=0$
and we are reduced to the case (401).

In the case {\bf 411}, when $s=1$, we get $x_{1k}-x_{1k}x_{2k}\in {\mathfrak A}^2$ and we can consider the expression
$(x_{1k}x_{2k}-x_{12}x_{2k})\tau_{1n}(I+\tau_{2n})$, see (\ref{xx-xx}). The case
(411) splits into two cases: (4110), corresponding to the cases $t\not=1$ and $t=1$ in (\ref{f-tg-in-l(2)}), when
$\sum_n(1-c_{1n})c_{2n}^2=\infty$ and (4111), when
$\sum_n(1-c_{1n})c_{2n}^2<\infty$ hence, $\sum_n(1-c_{1n})\alpha_{2n}^2(1)=\infty$
(see cases (3110), (3111) and (\ref{t=(not,1)})).

In the  case {\bf 4110} we can approximate
$x_{1k}x_{2k}-x_{12}x_{2k}$, by Lemma~\ref{l.xx-xx}.
Since $x_{1k}-x_{1k}x_{2k}\in {\mathfrak
A}^2$ we get  $x_{1k}-x_{12}x_{2k}\in {\mathfrak A}^2$ hence, we
can approximate $x_{12}$, by Lemma~\ref{x_12-in-A}, when $\sum_k\alpha_{2k}^2(1)=\infty$. This finish the proof.

When $\sum_k\alpha_{2k}^2(1)<\infty$ we conclude that $\lim_nc_{2n}=0$
and we are in the case (401).

In the  case {\bf 4111},  as in the case (3111), we can use the expression
$x_{1k}-x_{12}x_{2k}(2x_{1k}-I)$ (see (\ref{x-xxx})).  By
Lemma~\ref{l.x(12)-in-A!}, we can approximate  $x_{12}$ since in this case
$\sum_n(1-c_{1n})\alpha_{2n}^2(1)=\infty$ (see
Remark~\ref{r.x-xxx})). Since $x_{12},\,\,x_{2n}\in {\mathfrak A}^2$, by Remark~\ref{r.x(kn)-in-A} (c), we conclude that $x_{1k}\in {\mathfrak A}^2$.
Finally, we have $x_{1k},\,\,x_{2k+1}\in {\mathfrak A}^2$ for all $k\geq 2$.

{\bf Case 4II.}  Let for some $t\in {\mathbb R}\setminus \{0\}$ holds $\Vert G-tF\Vert^2<\infty$ and $\Gamma(G_m)/\Gamma(F_m)$\\$
\leq C$. This means that
$$
\Vert G-tF\Vert^2  =\sum_{n=3}^\infty\frac{\mid (1-c_{2n})-t(1-c_{1n})\mid^2}{d_n},\quad
\frac{\sum_{n=3}^m(1-c_{2n})^2/d_n}{\sum_{n=3}^m(1-c_{1n})^2/d_n}\leq C
$$
where $d_n=1-c_{1n}+1-c_{2n}+1-c_{1n}c_{2n}$. Set $x_n=1-c_{1n},\,\,y_n=1-c_{2n}$, and
$$
d_1(x,y)=2x+2y-xy,
\quad d_2(x,y)=x+y,\quad x,y\in [0,\,1].
$$
\begin{lem}
\label{l.d1-sim-d3}
We have for $x,y\in [0,\,1]$
\begin{equation}
\label{d1-sim-d3}
d_2(x,y)\leq d_1(x,y)\leq 2d_2(x,y).
\end{equation}
\end{lem}
\begin{pf}
Indeed, since $x+y-xy=1-(1-x)(1-y)\in [0,1]$ we get
(\ref{d1-sim-d3})
$
x+y\leq 2x+2y-xy\leq 2(x+y).
$
\qed\end{pf}
Using the relations
$$
d_n=1-c_{1n}+1-c_{2n}+1-c_{1n}c_{2n}=
d_1(x_n,y_n),\quad
1-c_{2n}+1-c_{1n}c_{2n}=d_2(x_n,y_n)
$$
and Lemma~\ref{l.d1-sim-d3} we conclude that the following equivalences hold:
\begin{align}
\label{Si,G,F-sim1}
\Sigma_{12}=\sum_{n}\frac{x_n^2(1-y_n)^2}{x_n(1-y_n)^2+y_n},\quad
\Vert F\Vert^2\sim\sum_{n}\frac{x_n^2}{x_n+y_n},\quad\\
\label{Si,G,F-sim2}
\Vert G\Vert^2\sim\sum_{n}\frac{y_n^2}{x_n+y_n},\quad
\Vert G-tF\Vert^2\sim\sum_{n}\frac{(y_n-tx_n)^2}{x_n+y_n}.
\end{align}
We have to consider only the following three  possibilities:

(a)  the case when $1>x_n\geq \varepsilon>0$ for all $n\in {\mathbb N}$, the set of all limit points is $[\varepsilon,\,1]$;

(b)  the case when $\lim_nx_n=0$, the set of limit points is one point $0$;

(c) the intermediate case, when the set of all limit points is the segment $[0,\,1]$, in this case we have  ${\mathbb N}_0$ and ${\mathbb N}_1$ two infinite subsets of ${\mathbb N}$ such that
$$
x_n\geq C>0\quad  \forall n\in {\mathbb N}_0,\quad\text{and}\quad
\lim_{n\in {\mathbb N}_1}x_n=0.
$$
Consider the expression
$
T_{kn}-I=x_{1k}\tau_{1n}+ x_{2k}\tau_{2n}+x_{1k}x_{2k}\tau_{1n}\tau_{2n}.
$
\begin{lem}
\label{l.x+x+x}
We can approximate by linear combinations $\sum_{n=3}^mt_n(T_{kn}-I)$ the following expressions:
\begin{equation}
\label{x+x+x}
x_{1k}+\beta^{(2)}_1x_{2k}-\beta^{(3)}_1x_{1k}x_{2k},\quad\text{or}\quad
\beta^{(1)}_2x_{1k}+x_{2k}-\beta^{(3)}_2x_{1k}x_{2k}
\end{equation}
where
\begin{equation}
\label{beta(2,3)}
\beta^{(2)}_1=\lim_m\frac{\sum_{n=3}^m\frac{1-c_{2n}}{1+c_{1n}}}
{\sum_{n=3}^m\frac{1-c_{1n}}{1+c_{1n}}},\quad\quad
\beta^{(3)}_1=\lim_m\frac{\sum_{n=3}^m\frac{(1-c_{1n})(1-c_{2n})}{1+c_{1n}}}
{\sum_{n=3}^m\frac{1-c_{1n}}{1+c_{1n}}},
\end{equation}
\begin{equation}
\label{beta(2,3).2}
\beta^{(1)}_2=\lim_m\frac{\sum_{n=3}^m\frac{1-c_{1n}}{1+c_{2n}}}
{\sum_{n=3}^m\frac{1-c_{2n}}{1+c_{2n}}},\quad\quad
\beta^{(3)}_2=\lim_m\frac{\sum_{n=3}^m\frac{(1-c_{1n})(1-c_{2n})}{1+c_{2n}}}
{\sum_{n=3}^m\frac{1-c_{2n}}{1+c_{2n}}}.
\end{equation}
\end{lem}
\begin{pf}
Indeed,  to obtain the first expression or the second one in (\ref{x+x+x}) we use the fact that  $\sum_{n=3}^mt_n\tau_{1n}\to I$ or $\sum_{n=3}^mt_n\tau_{2n}\to I$ (see Remark~\ref{approx-I} and Lemma~\ref{1in<2>})
where $t_n$ are defined respectively by the following formulas (see (\ref{a(n),t(n)=})):
$$
t_n=-\frac{1}{1+c_{1n}}\Big(\sum_{n=3}^m\frac{1-c_{1n}}{1+c_{1n}}\Big)^{-1},\quad
t_n=-\frac{1}{1+c_{2n}}\Big(\sum_{n=3}^m\frac{1-c_{2n}}{1+c_{2n}}\Big)^{-1}.
$$
Further, we should proceed exactly as in the proof of Lemma~\ref{l.a+beta b}.
\qed\end{pf}
\begin{ex}
\label{ex.4IIa}
Let $x_n=C\in (0,\,1)$ for all $n\in {\mathbb N}$, then
$\Sigma_{12}<\infty$ if and only if $\sum_n(1-y_n)^2=\sum_nc_{2n}^2<\infty$. Indeed, we have
{\small
$$
\infty\!>\!\Sigma_{12}\!=\!\sum_{n}\frac{x_n^2(1-y_n)^2}{x_n(1-y_n)^2+y_n}\!>\!
\sum_{n}\frac{x_n^2(1-y_n)^2}{x_n+y_n}\!=\!
C^2\sum_{n}\frac{(1-y_n)^2}{C+y_n}\sim\sum_n(1-y_n)^2.
$$
}
We show that $\Vert F\Vert^2=\Vert G\Vert^2=\infty$ and
$\Vert G-tF\Vert^2<\infty$ for some $t\not=0$. Indeed, we have
$$
\Vert F\Vert^2\!\sim\!\sum_n\frac{x_n^2}{x_n+y_n}\sim
\sum_n\frac{C^2}{C+1}\!=\!\infty,\,\,
\Vert G\Vert^2\sim\sum_n\frac{y_n^2}{x_n+y_n}\sim
\sum_n\frac{1}{C+1}\!=\!\infty,
$$
$$
\Vert G-tF\Vert^2\sim\sum_n\frac{\vert y_n -tx_n\vert^2}{x_n+y_n}\sim
\sum_n\frac{\vert y_n -tC\vert^2}{C+1}<\infty\quad\text{for}\quad tC=1.
$$
\end{ex}
Since $1-y_n=c_{2n}$ we conclude that $\sum_nc_{2n}^2<\infty$ therefore, $\lim_nc_{2n}=0$ and finally, we conclude by Lemma~\ref{l.x+x+x} (see (\ref{beta(2,3)})),  Toeplitz theorem~\ref{t.Toeplitz} and Lemma~\ref{diverg}
that $\beta^{(2)}_1=\frac{1}{C}>1$ and $\beta^{(3)}_1=1$. Set $\beta:=\frac{1}{C}>1$ then we get
$$
x_{1k}+\beta x_{2k}-x_{1k}x_{2k}=x_{1k}(1-x_{2k})+\beta x_{2k}\in  {\mathfrak A}^2.
$$
{\it A regular matrix summability method} is a matrix transformation of a convergent sequence which preserves the limit.
\begin{thm}[Otto Toeplitz \cite{Toe1911}]
\label{t.Toeplitz}
An infinite matrix $(a_{i,j})_{i,j\in \mathbb {N}}$  with com\-plex-valued entries defines a regular summability method, i.e.,
 \begin{equation}
 \label{toepl,l(t)=l(s)}
 \lim_{n\to\infty}t_n=\lim_{n\to\infty}s_n\quad \text{where}\quad t_n:=\sum_{n=1}^{n_k}a_{kn}s_n
 \end{equation}
if and only if it satisfies all of the following properties:
\begin{eqnarray*}
 (I)& \lim_{i\to \infty }a_{i,j}=0\quad j\in \mathbb {N}& \text{(every column sequence converges to 0)},\\
(II)&\lim_{i\to \infty }\sum _{j=0}^{\infty }a_{i,j}=1& \text{(the row sums converge to 1)},\\
(III)& \sup _{i}\sum _{j=0}^{\infty }\vert a_{i,j}\vert <\infty & \text{(the absolute row sums are bounded)}.
\end{eqnarray*}
 \end{thm}
%
\begin{lem}[A particular case of the Toeplitz theorem]
\label{diverg} Let us have three sequences of real numbers
$(a_n),(b_n)$ and $(\alpha_n)$ with
$(a_n)>0,\,\,\sum_{n\in{\mathbb N}}a_n=\infty,\,\,
\sum_{k=1}^m\mid b_k\mid\left(\sum_{k=1}^ma_k\right)^{-1}\leq
C,\,\,m\in{\mathbb N},$ for some $C>0$ and
$\lim_n\alpha_n=\alpha\not=0.$ Set
\begin{equation}
\label{beta}
\beta_{m}=\sum_{k=1}^mb_k\Big(\sum_{k=1}^ma_k\Big)^{-1},\quad
\beta_{m}(\alpha]=\sum_{k=1}^m\alpha_kb_k\Big(\sum_{k=1}^m
a_k\Big)^{-1},
\end{equation}
\begin{equation*}
\beta_{m}[\alpha)=\sum_{k=1}^mb_k\Big(\sum_{k=1}^m\alpha_k
a_k\Big)^{-1},\quad
\beta_{m}(\alpha)=\sum_{k=1}^m\alpha_kb_k\Big(\sum_{k=1}^m\alpha_k
a_k\Big)^{-1}.
\end{equation*}
If the limit exists $\lim_m\beta_{m}=\beta\in {\mathbb R},$ then
the following limits  also exist and we have
\begin{equation}
\label{beta()} \lim_m\beta_{m}(\alpha]=\alpha\beta,\quad
\lim_m\beta_{m}[\alpha)=\alpha^{-1}\beta,\quad
\lim_m\beta_{m}(\alpha)=
\beta.
\end{equation}
\end{lem}
To prove that  $I-x_{2k}\in {\mathfrak A}^2$ we calculate $[T_{1k},x_{1k}]$. The operators
$x_{1k}$ and $T_{\alpha_{1k}}$ have the following form in $H_{1k}=L^2({\mathbb F}_2,\mu_{\alpha_{1k}})$ (see (\ref{T-alpha.p=2}), (\ref{X(kn)}) and (\ref{8.ab})):
$$
x_{1k}=
\left(\begin{smallmatrix} 0&0\\
0&1\\
\end{smallmatrix}\right),\quad
T_{1k}=\left(\begin{smallmatrix} 0& a_n^{-1}
\\
a_n &0
\end{smallmatrix}\right)\quad\text{where}\quad a_n=\sqrt{\alpha_{1k}(0)/\alpha_{1k}(1)}.
$$
We show that
\begin{equation}
\label{[T,x]^2=-1(2)}
 [T_{1k},x_{1k}]^2=-I,\quad  [T_{1k},x_{1k}(I-x_{2k})]^2=-(I-x_{2k}).
\end{equation}
 Indeed, we get
\begin{equation*}
[T_{1k},x_{1k}]=T_{1k}x_{1k}-x_{1k}T_{1k}=
\left(\begin{smallmatrix} 0& a_n^{-1}
\\
-a_n &0
\end{smallmatrix}\right).
\end{equation*}
This implies (\ref{[T,x]^2=-1(2)}) since
$
[T_{1k},x_{1k}]^2=
\left(\begin{smallmatrix}
 -1&0\\
0&-1
\end{smallmatrix}\right),
$ and $x_{2k}^2=x_{2k}$.

%
Finally, we get $x_{2k}\in {\mathfrak A}^2$ therefore, $x_{1k}-x_{1k}x_{2k}\in {\mathfrak A}^2$ and we can use the following expression (see (\ref{x-xxx}))
$$
[T_{kn}-I-x_{2k}(T_{2n}-I)-(x_{1k}-x_{1k}x_{2k})(T_{1n}-1)]T_{2n}
=
[x_{1k}x_{2k}-x_{12}x_{2k}(2x_{1k}-I)]\tau_{1n}.
$$
Since
$\sum_n(1-c_{1n})c_{2n}^2\!<\! \sum_nc_{2n}^2\!<\!\infty$ we conclude
by Remark~\ref{r.x-xxx} that\\  $\sum_n(1\!-\!c_{1n})\alpha_{2n}^2(1)=\infty$. By Lemma~\ref{l.x(12)-in-A!} we get $x_{12}\in {\mathfrak A}^2$. Since $x_{12},\,\,x_{2n}\in {\mathfrak A}^2$, by Remark~\ref{r.x(kn)-in-A} (c), we conclude that $x_{1k}\in {\mathfrak A}^2$.
Finally, we have
$x_{1k},\,\,x_{2k+1}\in {\mathfrak A}^2$ for all $k\geq 2$ and the proof of the irreducibility of the example is finished.

Consider now the general case {\bf (a)}.  Since $0<\varepsilon\leq x_n=1-c_{1n}< 1$ for all $n\in {\mathbb N}$ we conclude that for some subsequence we get $\lim_k(1-c_{1n_k})=C_1\in [\varepsilon,\,1].$ As in Example~\ref{ex.4IIa} we conclude that $\Sigma_{12}\sim\sum_n(1-y_n)^2<\infty.$ We can repeat then step by step the  proof of the irreducibility as it was done in the Example~\ref{ex.4IIa}.

The case {\bf (c)} is similar to the case (a). In this case we conclude that
{\small
$$
\infty\!>\!\Sigma_{12}\!=\!\!\!
\sum_{n\in {\mathbb N}_0}\frac{x_n^2(1-y_n)^2}{x_n(1-y_n)^2+y_n}
\!>\!\!\sum_{n\in {\mathbb N}_0}\frac{x_n^2(1-y_n)^2}{x_n+y_n}\!\geq\!\! C^2\!\!
\sum_{n\in {\mathbb N}_0}\frac{(1\!-\!y_n)^2}{1+y_n}\!\sim\!\!\sum_{n\in {\mathbb N}_0}(1\!-\!y_n)^2.
$$
}
Therefore, $\sum_{n\in {\mathbb N}_0}(1-y_n)^2<\infty$.  Set ${\mathbb N}^m_0={\mathbb N}_0\bigcap [1,\,m]$ and define  $\beta^{(2)}_0$ and $\beta^{(3)}_0$ as follows:
\begin{equation}
\label{beta(2,3)0}
\beta^{(2)}_0=\lim_m\frac{\sum_{n\in {\mathbb N}^m_0}\frac{1-c_{2n}}{1+c_{1n}}}
{\sum_{n\in {\mathbb N}^m_0}\frac{1-c_{1n}}{1+c_{1n}}},\quad\quad
\beta^{(3)}_0=\lim_m\frac{\sum_{n\in {\mathbb N}^m_0}\frac{(1-c_{1n})(1-c_{2n})}{1+c_{1n}}}
{\sum_{n\in {\mathbb N}^m_0}\frac{1-c_{1n}}{1+c_{1n}}}.
\end{equation}
Since $\lim_{n\in {\mathbb N}_0}c_{2n}=0$ we conclude that
$$
\beta^{(2)}_0=(\lim_{n\in {\mathbb N}_0} x_n)^{-1}=C_1^{-1}\in [1,\,C^{-1})\quad\text{and}\quad \beta^{(3)}_0=1.
$$
We repeat step by step the proof done in the Example~\ref{ex.4IIa} to the case {\bf (a)}.

We show that the case {\bf (b)} can not be realized. Indeed, let $\lim_nx_n=0$.
Since for some $t\in {\mathbb R}\setminus \{0\}$ holds
$$
\Vert G-tF\Vert^2  \sim\sum_{n=3}^\infty\frac{\mid y_n-tx_n\mid^2}{x_n+y_n}<\infty,\quad\text{so}\quad
$$
$$
0=\lim_n\frac{\mid y_n-tx_n\mid^2}{x_n+y_n}\geq \frac{1}{2}\lim_n\mid y_n-tx_n\mid^2=\frac{1}{2}(\lim_ny_n-t\lim_nx_n)^2=\frac{1}{2}(\lim_ny_n)^2.
$$
Therefore, $\lim_ny_n=0$. This contradicts with two conditions:
$$\Sigma_{12}<\infty\quad\text{ and}\quad \Vert F\Vert^2\sim\sum_{n}\frac{x_n^2}{x_n+y_n}=\infty.
$$
Indeed, fix some $\varepsilon>0$. For sufficiently big $N\in {\mathbb N}$ we get
$$
\infty>\Sigma_{12}>\sum_{n>N}\frac{x_n^2(1-y_n)^2}{x_n+y_n}\geq (1-\varepsilon)^2 \sum_{n>N}\frac{x_n^2}{x_n+y_n}\sim \Vert F\Vert^2=\infty.
$$


We give another proof of the irreducibility in the {\bf case 1}.

The case $\Sigma_{12}=\infty$, in fact, is included in the cases (2),
(3) and (4), we shall denote them respectively by (2*), (3*) and
(4*). Since $\Sigma_{12}=\infty$ we have $x_{12}\in {\mathfrak
A}^2$.

Case {\bf (2*)}. Let $\Gamma(G)<\infty$. Then $\Gamma(H_1)<\infty$
therefore, $\Sigma_{1m}\to\infty$ hence, $x_{1k}\in {\mathfrak
A}^2$ for $k\geq 3$. In addition $x_{12}\in {\mathfrak A}^2$.
Since $x_{1k},\,\,T_{\alpha_{1n}}\in {\mathfrak A}^2$ we conclude, by Remark~\ref{r.x(kn)-in-A}(b), that $x_{2k}\in
{\mathfrak A}^2,\,\,k>3$. Finally, $x_{1k},\,\,x_{2k+1}\in {\mathfrak
A}^2$ for $k\geq 2$.

Case {\bf (3*)}. Let $\Gamma(F)<\infty$. Then $\Gamma(H_2)<\infty$
therefore, $\Sigma_{2m}\to\infty$ hence, $x_{2k}\in {\mathfrak
A}^2$ for $k\geq 3$. As in the case {\bf (30)} $(a),\,\,
(b)\,\,\text{or}\,\,(c)$ we get $x_{1k}-x_{12}x_{2k}\in {\mathfrak
A}^2$.
Since $x_{12},\,x_{2k}\in {\mathfrak A}^2$ for $k\geq 3$ we
conclude that $x_{1k},\,\,x_{2k+1}\in {\mathfrak A}^2$ for $k\geq
2$. The proof is finished.

In the opposite case, i.e., $(\overline a) \cap (\overline b)\cap (\overline c)$,
by Lemma~\ref{l.a+beta b}, we can
approximate the following expression:
$x_{1k}-(1-s)x_{12}x_{2k}-sx_{1k}x_{2k}$. In the case $s\not=1$
since $x_{12},\,x_{2k}\in {\mathfrak A}^2$ we conclude that
$x_{1k}-sx_{1k}x_{2k}=x_{1k}(1-sx_{2k})\in {\mathfrak A}^2$ hence,
$x_{1k}\in {\mathfrak A}^2$ (see the case {\bf (310)}). The proof is finished.

If $s=1$ we get $x_{1k}-x_{1k}x_{2k}\in {\mathfrak A}^2$. Since $x_{12}\in {\mathfrak A}^2$ we conclude that $T_{\alpha_{2k}}\in {\mathfrak A}^2$ for $k
\geq 3$, by Remark~\ref{r.x(kn)-in-A} (a) therefore, (see (\ref{[T,x]^2=-1(2)}))
\begin{equation}
\label{[T,x]^2=-1(1)}
[T_{\alpha_{2k}},x_{1k}-x_{1k}x_{2k}]^2=(-x_{1k}[T_{\alpha_{2k}},x_{2k}])^2=-x_{1k}.
\end{equation}
At last, we have $x_{1k},\,\,x_{2k+1}\in {\mathfrak A}^2$ for $k\geq
2$ and the proof is finished.

Case {\bf (4*)}. Let $\Sigma_{12}=\Gamma(F)=\Gamma(G)=\infty$, then $x_{12}\in {\mathfrak A}^2$ and  $T_{\alpha_{2n}}\in {\mathfrak A}^2$ for $n\geq 3$.
 Using Lemma~\ref{l.x+x+x}, we can approximate by linear combinations $\sum_{n=3}^mt_n(T_{kn}-I)$ the following expressions:
$$
x_{1k}+\beta^{(2)}_1x_{2k}-\beta^{(3)}_1x_{1k}x_{2k},\quad\text{or}\quad
\beta^{(1)}_2x_{1k}+x_{2k}-\beta^{(3)}_2x_{1k}x_{2k}
$$
since one of two sequence
$$
\sum_{n=3}^m(1-c_{1n})\Big(\sum_{n=3}^m (1-c_{2n})\Big)^{-1}\quad\text{or}\quad
\sum_{n=3}^m(1-c_{2n})\Big(\sum_{n=3}^m (1-c_{1n})\Big)^{-1}
$$
should be bounded.

 Because of the symmetry between the first and the second rows, i.e., between variables $(x_{1k})_k$ and $(x_{2k})_k$, it is sufficient to consider the case when $x_{1k}+\beta^{(2)}_1x_{2k}-\beta^{(3)}_1x_{1k}x_{2k}\in {\mathfrak A}^2$   where $0<\beta^{(3)}_1\leq \beta^{(2)}_1<\infty$.  By (\ref{[T,x]^2=-1(2)}) we get
$$
[T_{\alpha_{1k}},x_{1k}(I-\beta^{(3)}_1x_{2k})+\beta^{(2)}_1x_{2k}]^2=
-(I-\beta^{(3)}_1x_{2k}),
$$
 therefore, $x_{2k}\in {\mathfrak A}^2$ for $k\geq 3$ when $\beta^{(3)}_1>0$. By  Remark ~\ref{r.x(kn)-in-A} (c), we get that $x_{1k}\in {\mathfrak A}^2$ for $k\geq 3$ and the proof is finished.

Let $\beta^{(2)}_1>\beta^{(3)}_1=0$, then $x_{1k}+\beta^{(2)}_1x_{2k}\in {\mathfrak A}^2$. We prove the following
\begin{lem}
\label{l.ABD-irr}
 The von Neumann algebra $C_n$ generated by operators $T_{\alpha_{1n}},\,\,T_{\alpha_{2n}}$ and $x_{1n}+\beta x_{2n}$ is irreducible in the space $H_n:=L^2({\mathbb
F}_2,\mu_{\alpha_{1n}})\otimes L^2({\mathbb
F}_2,\mu_{\alpha_{2n}})$ for $\beta \in (0,\,1]$.
\end{lem}
\begin{pf}
Using (\ref{A-otimes- B}), (\ref{T1nT2n}),  (\ref{X(kn)}) and Remark~\ref{otimes=1}  we get
\begin{equation}
\label{ABD-irr}
 T_{\alpha_{1n}}\!\!=\!
 \left(\begin{smallmatrix}
 0  &0&a_n^{-1}&0\\
 0  &0&0&a_n^{-1}\\
 a_n&0&0&0     \\
 0  &a_n&0&0
\end{smallmatrix}\right),\,\,
T_{\alpha_{2n}}\!\!=\!\!
 \left(\begin{smallmatrix}
 0&b_n^{-1}&0&0\\
 b_n&0&0&0\\
 0&0&0&b_n^{-1}\\
 0&0& b_n&0
\end{smallmatrix}\right),\,\,x_{1n}\!+\!\beta x_{2n}\!\!=\!\!
 \left(\begin{smallmatrix}
 0&0&0&0\\
0&\beta &0&0\\
 0&0&1&0\\
 0&0&0&1+\beta
\end{smallmatrix}\right)
\end{equation}
where $a_n,\,b_n$ are defined by (\ref{8.ab}). Indeed, we have
$$
 T_{\alpha_{1n}}\otimes I=
\left(\begin{smallmatrix}
 0&a_n^{-1}\\
a_n&0
\end{smallmatrix}\right)\otimes
 \left(\begin{smallmatrix}
 1&0\\
 0&1
\end{smallmatrix}\right)=
 \left(\begin{smallmatrix}
 0  &0&a_n^{-1}&0\\
 0  &0&0&a_n^{-1}\\
 a_n&0&0&0     \\
 0  &a_n&0&0
\end{smallmatrix}\right),
$$
$$
 I\otimes T_{\alpha_{2n}}=
\left(\begin{smallmatrix}
 1&0\\
 0&1
\end{smallmatrix}\right)\otimes
\left(\begin{smallmatrix}
 0&b_n^{-1}\\
b_n&0
\end{smallmatrix}\right) =
 \left(\begin{smallmatrix}
 0&b_n^{-1}&0&0\\
 b_n&0&0&0\\
 0&0&0&b_n^{-1}\\
 0&0& b_n&0
\end{smallmatrix}\right),
$$
$$
x_{1n}\otimes I+\beta I\otimes x_{2n}=
 \left(\begin{smallmatrix}
 0&0\\
 0&1
\end{smallmatrix}\right)\otimes
 \left(\begin{smallmatrix}
 1&0\\
 0&1
\end{smallmatrix}\right)+\beta
 \left(\begin{smallmatrix}
 1&0\\
 0&1
\end{smallmatrix}\right)\otimes
 \left(\begin{smallmatrix}
 0&0\\
 0&1
\end{smallmatrix}\right)=
$$
$$
 \left(\begin{smallmatrix}
 0&0&0&0\\
0&0 &0&0\\
 0&0&1&0\\
 0&0&0&1
\end{smallmatrix}\right)+\beta
 \left(\begin{smallmatrix}
 0&0&0&0\\
0&1 &0&0\\
 0&0&0&0\\
 0&0&0&1
\end{smallmatrix}\right)=
 \left(\begin{smallmatrix}
 0&0&0&0\\
0&\beta &0&0\\
 0&0&1&0\\
 0&0&0&1+\beta
\end{smallmatrix}\right).
$$
In the case $\beta\in (0,\,1)$ the commutant $(x_{1n}\!+\!\beta x_{2n})'$ consists of all diagonal operators $D(\lambda)={\rm diag}(\lambda_1,\dots, \lambda_4)$ since eigenvalues of $x_{1n}\!+\!\beta x_{2n}$ are distinct. The commutation relation
$[D, T_{\alpha_{1n}}\otimes I]=0$ implies  $\lambda_1=\lambda_3,\,\,\lambda_2=\lambda_4$. The commutation relation $[D, I\otimes T_{\alpha_{2n}}, I\otimes T_{\alpha_{2n}}]=0$ implies  $\lambda_1=\lambda_2,\,\,\lambda_3=\lambda_4$. Finally, we get $D(\lambda)=\lambda I$.
In the case $\beta=1$ the commutant $(x_{1n}\!+\!\beta x_{2n})'$ consists of all  operators of the form
$$
D= \left(\begin{smallmatrix}
 \lambda_1&0&0&0\\
0& \lambda_2 &b&0\\
 0&c& \lambda_3&0\\
 0&0&0& \lambda_4
\end{smallmatrix}\right).
$$
The commutation relation
$[D, T_{\alpha_{1n}}\otimes I]=0$ implies $b=c=0,\,\,\lambda_1=\lambda_3,\,\,\lambda_2=\lambda_4$.
The commutation relation $[D, I\otimes T_{\alpha_{2n}}, I\otimes T_{\alpha_{2n}}]=0$ implies  $\lambda_1=\lambda_2,\,\,\lambda_3=\lambda_4$. Hence, in this case we get $D=\lambda I$.
\qed\end{pf}
The irreducibility of the representation in the case  $\beta^{(2)}_1>\beta^{(3)}_1=0$
follows from the fact that von Neumann algebra
${\mathfrak A}=(T_{12},x_{12})''\otimes_{n=3}^\infty C_n$
is irreducible since  the commutant ${\mathfrak A}'$ is trivial by Lemma~\ref{l.ABD-irr}. Indeed, we have
${\mathfrak A}'=(T_{12},x_{12})'\otimes_{n=3}^\infty C_n'.$

When $\beta^{(2)}_1=\beta^{(3)}_1=0$ we get $x_{1n}\in {\mathfrak A}^2 $.
By Remark ~\ref{r.x(kn)-in-A} (b), we conclude  that $x_{2k}\in {\mathfrak A}^2$
for $k\geq 3$ and the proof is finished.

\newpage
{\it Acknowledgements.} {The authors expresses his deep gratitude to
the Max Planck Institute for Mathematics for the financial grant and the hospitality in 2016--2017}.

\end{document}